\documentclass[oneside,english]{amsart}
\usepackage[T1]{fontenc}
\usepackage[latin9]{inputenc}
\usepackage{amsmath}
\usepackage[pdftex]{graphicx}
\usepackage{microtype}

\makeatletter
\numberwithin{equation}{section}
\numberwithin{figure}{section}

\makeatother

\usepackage{babel}
\begin{document}

\title[\tiny{Higher order tangents and Higher order Laplacians on Sierpinski Gasket Type Fractals}]{Higher order tangents and Higher order Laplacians on Sierpinski Gasket Type Fractals}
%    author one information
\author{Shiping Cao}
\address{School of Physics, Nanjing University, Nanjing, 210093, P.R. China.}
\curraddr{} \email{shipingcao@hotmail.com}
\thanks{}

%    author two information
\author{Hua Qiu$^*$}
\address{Department of Mathematics, Nanjing University, Nanjing, 210093, P. R. China.}
\curraddr{} \email{huaqiu@nju.edu.cn}
\thanks{$^*$ Corresponding author. }
\thanks{The research of the second author was supported by the Nature Science Foundation of China, Grant 11471157.}

%    For articles to be published after 1 January 2010, you may use
%    the following version:
%\subjclass[2010]{Primary }


\date{}

\dedicatory{}
\begin{abstract}

We study higher order tangents and higher order Laplacians on p.c.f. self-similar sets with fully symmetric structures, such as $D3$ or $D4$ symmetric fractals. Firstly, let $x$ be a vertex point in the graphs that approximate the fractal,  we prove that for any $f$ defined near $x$, the higher oder weak tangent of $f$ at $x$, if exists, is the uniform limit of local multiharmonic functions that agree with $f$ in some sense near $x$. Secondly, we prove that the higher order Laplacian on a fractal can be expressible as a renormalized uniform limit of higher order graph Laplacians on the graphs that approximate the fractal. The main technical tool is the theory of local multiharmonic functions and local monomials analogous to $(x-x_0)^j/j!$. The results in this paper are closely related to the theory of local Taylor approximations, splines and entire analytic functions. Some of our results can be extended to general p.c.f. fractals. In Appendix of the paper, we provide a recursion algorithm for the exact calculations of the boundary values of the monomials for $D3$ or $D4$ symmetric fractals, which is shorter and more direct than the previous work on the Sierpinski gasket. 
\end{abstract}
\maketitle

\section{Introduction }

Analysis on post critically finite (p.c.f.) fractals has been well developed since Kigami's original papers {[}Ki1,Ki2{]} on the direct analytic construction of the Laplacian on the familiar Sierpinski gasket $\mathcal{SG}$ (see Fig. 1.1), which now has been viewed as the ``post child'' for the theory. Recently, there are several works in connection with the differential calculus on fractals that involve derivatives, tangents, multiharmonic functions, higher order Laplacians, analogous to the theory of analysis on manifolds. Please see {[}BSSY, CQ, DRS, DSV, NSTY, S1-S3, SU,T2{]} and the references therein. 

\begin{figure}[h]
\centering
\includegraphics[width=6cm]{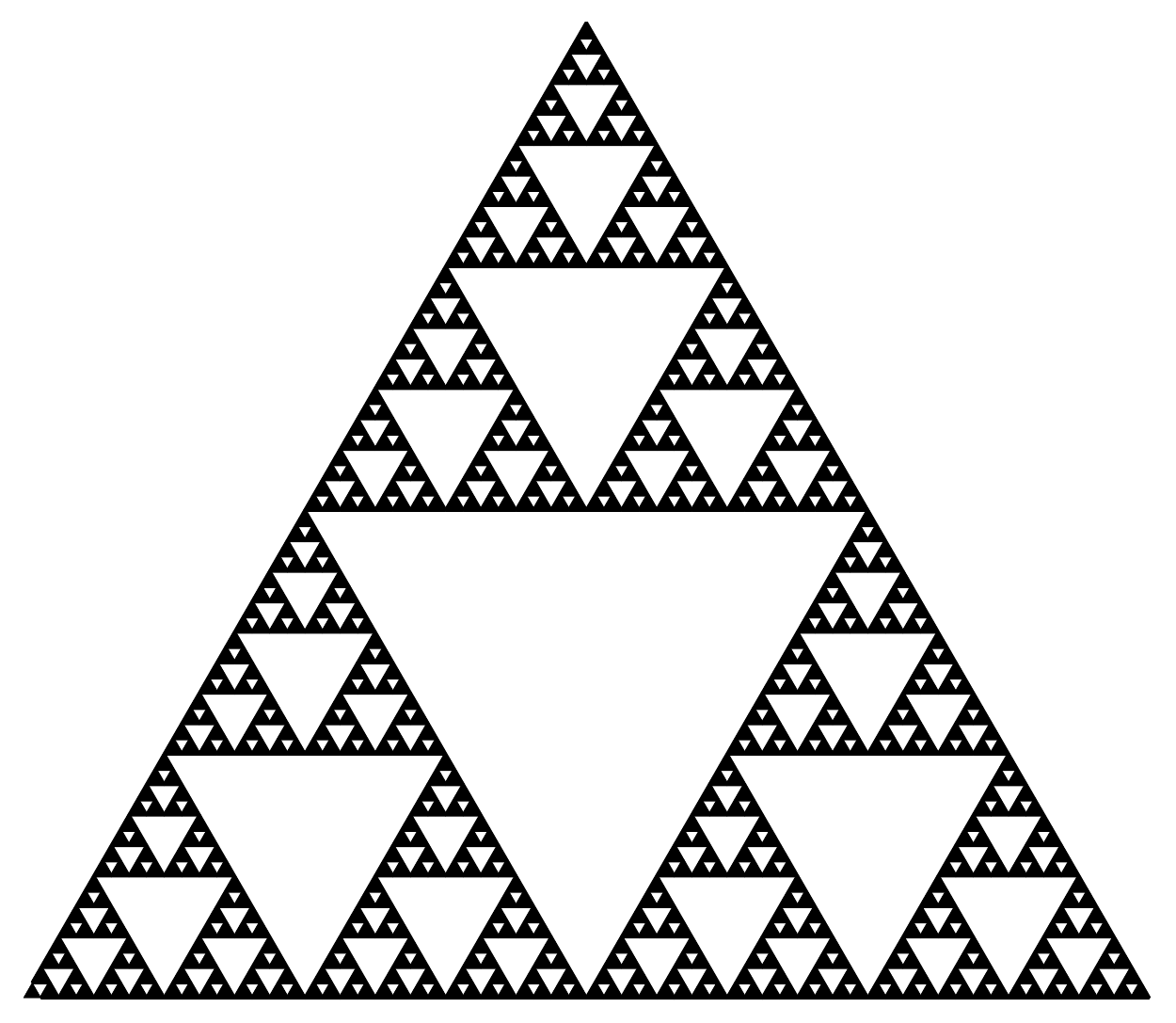}
\begin{center}
 \caption{The Sierpinski gasket $\mathcal{SG}$.}
 \end{center}
\end{figure}

In {[}S3{]}, Strichartz developed a theory of derivatives and gradients on a class of p.c.f. fractals with nondegenerate harmonic structure, which consummates the theory of normal derivatives and Gauss-Green's integration formula, and could be used to get an analogous theory of local Taylor approximations at vertices on fractals. There are also some other works concerned the gradients and tangents on fractals from different points of views, please see  {[}Ki3{]}, {[}Ku{]} and {[}T2{]}. In Teplyaev's work {[}T2{]} one could find a discussion on the relations between the different definitions and results of Kigami {[}Ki3{]}, Kusuoka {[}Ku{]}, Teplyaev and Strichartz on this topic.

In the end of Strichartz's paper {[}S3{]}, he post sevaral open problems that should be solved to complete the story of local Taylor approximations. Two of them are as follows.

\textbf{Question 1.} \textit{For a smooth function $f$ defined near a vertex $x$, could the higher order tangents $T_n(f)$ at $x$ be expressible as limits of local multiharmonic functions that agree with $f$ in a suitable sense near $x$?}

\textbf{Question 2.} \textit{For a smooth function $f$, is it possible to express $\Delta_\mu^n f$ as a uniform limit of a pointwise formula in terms of linear combinations of the values of $f$ at vertices approaching $x$?}
 
The main goal of this paper is to answer the above two questions. 

We will mainly focus on the $D3$ symmetric fractals, i.e., those fractals whose boundary consists of $3$ points and all structures posses full $D3$ symmetry. The results can be extended to fully symmetric p.c.f. fractals with suitable modification. 

Some of our results are extended to general p.c.f. fractals. 

This paper is organized as follows. In Section 2, we collect some notations and facts about Laplacians and derivatives on general p.c.f. fractals, most of which was introduced in {[}Ki5, S4{]}. In Section 3, we introduce the theory of local monomials which form a basis of local multiharmonic functions near a vertex $x$, using which we give a positive answer of Question 1. This theory will also play a key role for solving Question 2. Then in Section 4, the main part of this paper, we focus on Question 2 and prove a pointwise formula for the higher order Laplacians. In both  Section 3 and Section 4 we only consider those $D3$ symmetric fractals, since the full symmetric structures could provide many advantages for our discussion.  Then in Section 5, we turn to consider the general p.c.f. fractals to see to what extent can we extend our previous results. In Section 6, we show how to extend the previous results to those p.c.f. self-similar sets with full symmetry,  including some typical $D4$ symmetric fractals, such as the tetrahedral Sierpinski gasket and the Vicsek sets. Finally, in Appendix, we provide a recursion algorithm for the exact calculations of the boundary values of the monomials introduced in Section 3,  for some typical $D3$ or $D4$ symmetric fractals, which are important for results in Section 3. This algorithm is an improved version of the one developed in {[}NSTY{]}, which presents to be more direct and shorter.

\section{Notations}

We recall some  standard notations and results on Kigami's Laplacian and Stirchartz's derivatives on p.c.f. fractals, which are the necessary background of this paper. Please refer to {[}Ki1-Ki2, Ki5, S4{]} for any unexplained notion.

Let $(K, N, \{F_l\}_{0\leq l<N})$ be a \textit{ p.c.f. self-similar structure}. That is, there is a finite set of contractive continuous injections $\{F_l\}_{0\leq l<N}$ on some metric space, with a compact invariant set $K$ satisfying $K=\bigcup_{0\leq l<N}F_lK$. We define $W_m$ as the space of  \textit{words} 
$w=w_1\cdots w_m$ of length $|w|=m$, taking values from the alphabet $\{0,\dots,N-1\}$, then we could denote 
$F_w=F_{w_1}\circ\cdots\circ F_{w_m}$ and call $F_wK$ a  $m$ level \textit{cell} of $K$. The term ``p.c.f.'' means that $K$ is connected, and there is a finite set $V_0\subseteq K$ called the \textit{boundary} of $K$ such that $F_wK\cap F_{w'}K\subset F_w V_0\cap F_{w'}V_0$ for different $w$ and $w'$ with the same length. Moreover, each element in the  boundary set $V_0$ is required to be the fixed point of one of the mapping of $\{F_l\}_{0\leq l<N}$. Without loss of generality, we write $V_0=\{q_0, \cdots,q_{N_0-1}\}$ for $N_0\leq N$, and $F_l q_l=q_l$ for $l<N_0$.

Let $G_0$ denote the complete graph on $V_0$. We approximate $K$ by a sequence of \textit{graphs} $G_m$ with vertices $V_m$ and edge relation $x\sim_my$ defined by inductively applying the contractive mappings of $\{F_l\}$ to $G_0$. Let $V_*=\bigcup_{m\geq 0} V_m$ be the collection of all vertices of $K$.

Suppose there is a regular \textit{harmonic structure} on $(K, N, \{F_l\}_{0\leq l<N})$. Thus there is a sequence of \textit{renormalized  graph energies}   $\mathcal{E}_m$ on $G_m$ with
$$\mathcal{E}_m(f,g)=\sum_{x\sim_my}c_{xy}(f(x)-f(y))(g(x)-g(y))$$
for functions $f,g$ defined on $V_m$, satisfying the self-similar identity 
$$\mathcal{E}_m(f,g)=\sum_{l=0}^{N-1}r_l^{-1}\mathcal{E}_{m-1}(f\circ F_l, g\circ F_l),$$
where $c_{xy}$ are the $m$-level \textit{conductances} on graph $G_m$, and $\{r_l\}_{0\leq l<N}$ are the \textit{renormalization factors} satisfying $0<r_l<1$. For $0\leq i,j< N_0$, we use $c_{ij}$ to denote the $0$-level conductances on graph $G_0$. Obviously, for  $x\sim_my$, we have $c_{xy}=r_w^{-1} c_{ij}$, where $w$ is the word of length $m$ such that $x=F_wq_i$, $y=F_wq_j$, with $r_w=r_{w_1}\cdots r_{w_m}$.  Furthermore, if we denote $\mathcal{E}_m(f)=\mathcal{E}_m(f,f)$, then the restriction of $\mathcal{E}_m$ to $G_{m-1}$ equals $\mathcal{E}_{m-1}$, which means, if $f$ is defined on $G_{m-1}$, then for all extension $f'$ of $f$ to $G_m$, the one $\tilde{f}$ that minimize $\mathcal{E}_m$ satisfies $\mathcal{E}_m(\tilde{f})=\mathcal{E}_{m-1}(f)$. Hence the sequence $\{\mathcal{E}_m(f)\}$ is monotone increasing as $m$ goes to infinity for any function $f$ defined on $K$, and thus we could define
$$\mathcal{E}(f)=\lim_{m\rightarrow\infty}\mathcal{E}_m(f),$$
and by polarization identity, 
$$\mathcal{E}(f,g)=\lim_{m\rightarrow\infty}\mathcal{E}_m(f,g).$$
The domain $dom\mathcal{E}$ of $\mathcal{E}$ consists of continuous functions $f$ such that $\mathcal{E}(f)<\infty$.  The self-similar identity for graph energy becomes 
$$\mathcal{E}(f,g)=\sum_{w\in W_m}r_w^{-1}\mathcal{E}(f\circ F_w, g\circ F_w).$$

A function $h$ is \textit{harmonic} if it minimizes the energy from level $m$ to level $m+1$ for each $m$. All the harmonic functions form a $N_0$-dimensional space, denoted by $\mathcal{H}_0$, and hence any values on the boundary can uniquely determine a harmonic function on $K$. In fact, for every $0\leq l<N_0$, there is linear map $M_l: \mathcal{H}_0\rightarrow\mathcal{H}_0$ be defined by $M_lh=h\circ F_l$. We call $M_l$ the $l$-th \textit{harmonic extension matrix}. 

Let $\mu$ be the \textit{self-similar measure} with a set of probability weights $\{\mu_l\}$ on $K$, satisfying 
$$\mu(A)=\sum_{0\leq l<N}{\mu_l}\mu(F_l^{-1}A)$$
or equivalently,
$$\int_K fd\mu=\sum_{0\leq l<N}\mu_l\int_K f\circ F_ld\mu.$$
For $w\in W_m$, we denote $\mu_w=\mu_{w_1}\cdots\mu_{w_m}$ the measure of $F_wK$.

The \textit{graph Laplacian} $\Delta_m$ on $G_m$ is defined to be 
$$\Delta_mf(x)=\sum_{y\sim_mx}c_{xy}(f(y)-f(x))$$
for $x\in V_m\setminus V_0$. The \textit{Laplacian} with respect to $\mu$ on $K$ is defined as the renormalized limit
$$\Delta_\mu f(x)=\lim_{m\rightarrow\infty}\tilde{\Delta}_m f(x),$$
where $\tilde{\Delta}_mf(x)=(\int_K\psi_x^md\mu)^{-1}\Delta_mf(x)$. (We avoid the notation $\tilde{\Delta}_{\mu,m}$ without causing any confusion.) Here $\psi_x^m$ is a \textit{tent function} which is  harmonic on each  $m$-level cell taking value $1$ at $x$ and $0$ at other vertices in $V_m$.
More precisely, $f\in dom\Delta_\mu$ and $\Delta_\mu f=g$ means $f$ and $g$ are continuous and the above limit converges to $g$ uniformly on $V_*\setminus V_0$. There is an equivalent definition called \textit{weak formulation}, which says that for $f\in dom\mathcal{E}$ and continuous function $g$, $f\in dom\Delta_\mu$ with $\Delta_\mu f=g$ if and only if 
$$\mathcal{E}(f,v)=-\int_K gvd\mu$$
holds for all $v\in dom_0\mathcal{E}$, where $dom_0\mathcal{E}$ means those functions in $dom\mathcal{E}$ that vanishes on the boundary $V_0$.

There is a scaling identity
$$\Delta_\mu(f\circ F_w)=r_w\mu_w(\Delta_\mu f)\circ F_w.$$

The space of \textit{multiharmonic functions} (solutions of $\Delta_\mu^n h=0$ for some $n$) on fractals, analogous to \textit{polynomials} on the unit interval plays an important role in describing the approximation behavior of smooth functions, such as in the theory of Taylor approximations {[}S3{]}, splines {[}SU{]}, and power series expansions {[}NSTY{]}. Let $\mathcal{H}_n$ denote the collection of $(n+1)$-harmonic functions, the solutions of $\Delta_\mu^{n+1} h=0$, which is of $(n+1)N_0$-dimension.  

The following is the \textit{Gauss-Green's formula},
$$\mathcal{E}(f,g)=-\int_K \Delta_\mu f gd\mu+\sum_{q_l\in V_0}{\partial_n f(q_l) g(q_l)},$$ which connect the Laplacian $\Delta_\mu$ with the important concept of normal derivative.

We would not want to involve the general theory of derivatives for general p.c.f. fractals. 
In the rest of this section, we restrict our attention to the \textit{$D3$ symmetric fractals}. Here $D3$ symmetry means that all structures are invariant under any homeomorphism of $K$. In this case, $N_0=3$ and we could choose all $c_{ij}=1$. Now all the harmonic extension matrices $M_l$ only differ by permutations, and we must have $r_0=r_1=r_2$ and $\mu_0=\mu_1=\mu_2$. We denote them by $r$ and $\mu$ for simplicity respectively. We denote $\rho$ the value of $r\mu$. (In the next two sections, we actually need $r_l\mu_l=\rho$ for all $0\leq l < N$.)  It is easy to verify that $1$ is the largest eigenvalue and $r$ is the second large eigenvalue of the matrix $M_l$ for $l=0,1,2$.  We denote the \textit{third eigenvalue} by $\lambda$. Here we require that the matrix $M_l$ to be nondegenerate. 
  The familiar \textit{Sierpinski gasket} $\mathcal{SG}$ is a typical example, which is an invariant set generated by $3$ contractive mappings with fixed points $q_0,q_1,q_2$ the vertices of a triangle with contraction ratio $1/2$. For $\mathcal{SG}$, $r=3/5$, $\mu=1/3$, and $\rho=\lambda=1/5$. Two more examples are the \textit{level $3$ Sierpinski gasket} $\mathcal{SG}_3$ and the \textit{hexagasket} $\mathcal{HG}$. Here $\mathcal{SG}_3$ is an invariant set of six contractions of ratio $1/3$ as shown in Fig. 2.1, which has $r=7/15$, $\lambda=1/15$, $\mu=1/6$ and $\rho=7/90$. While $\mathcal{HG}$, which is also named as Star of David, is generated by six mappings with simultaneously rotating and contracting by a ratio of 1/3 as shown in Fig. 2.2, having
$r=3/7, \lambda=1/7$, $\mu=1/6$ and $\rho=1/14$. Please refer to {[}S4{]} for detailed information.

\begin{figure}[h]
\centering
\includegraphics[width=6cm]{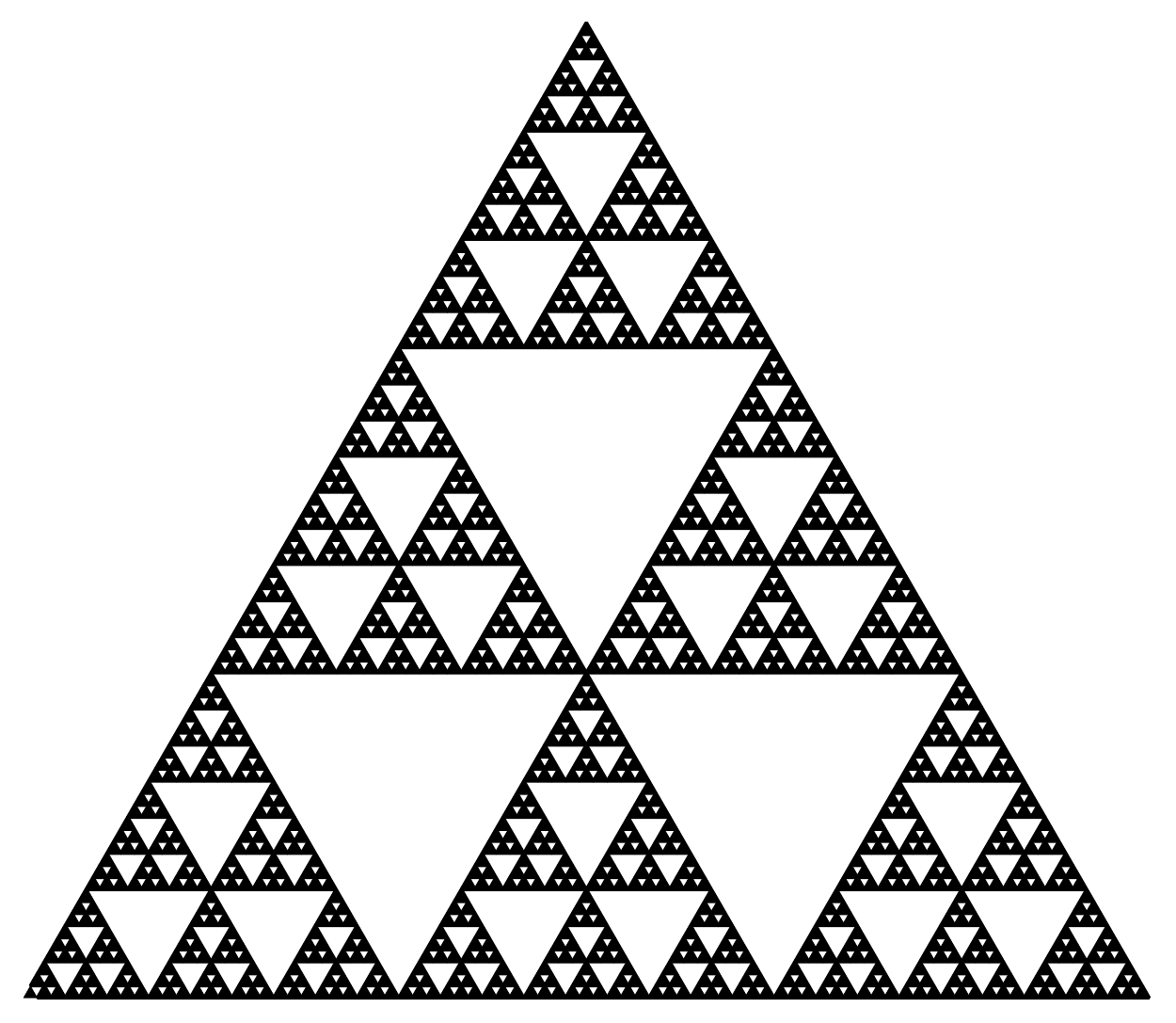}
\begin{center}
 \caption{The level $3$ Sierpinski gasket $\mathcal{SG}_3$.}
\end{center}
\end{figure}

\begin{figure}[h]
\centering
\includegraphics[width=6cm]{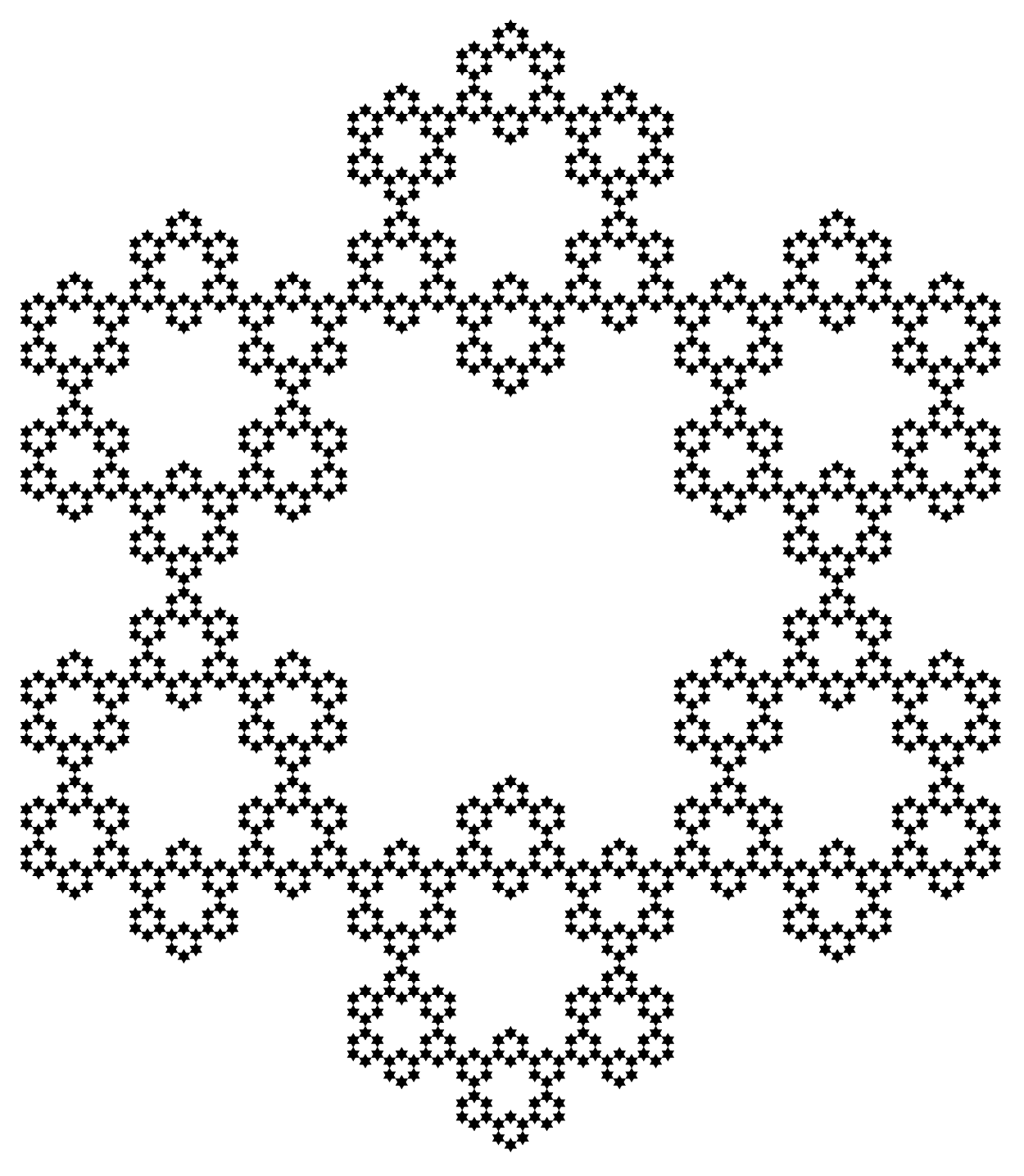}
\begin{center}
 \caption{The hexagasket $\mathcal{HG}$.}
\end{center}
\end{figure}

The \textit{normal derivatives} of a function $f$ at the boundary point $q_l$ is defined as
$$\partial_n f(q_l)=\lim_{m\rightarrow\infty}r^{-m}(2f(q_l)-f(F_{l}^m q_{l+1})-f(F_{l}^m q_{l-1}))$$
(cyclic notation $q_{l+3}=q_l$), while the \textit{transverse derivatives} at $q_l$ is defined as
$$\partial_T f(q_l)=\lim_{m\rightarrow\infty}\lambda^{-m}(f(F_{l}^m q_{l+1})-f(F_{l}^m q_{l-1})),$$ providing the limits exist. For hormonic functions, these derivatives can be evaluated without taking the limit. 

All the above notations and results are from global viewpoint. Now we turn to the localized ones. 

We could localize the definition of derivatives as follows. Let $x$ be a boundary point of cell $F_wK$, that is, there exists a $q_l$ such that $x=F_wq_l$. We define the normal derivative at $x$ with respect to $F_wK$ by
$$\partial_n^w f(x)=\lim_{m\rightarrow\infty}r_w^{-1}r^{-m}(2f(x)-f(F_wF_l^mq_{l+1})-f(F_wF_l^mq_{l-1}))$$
if the limit exists. We denote the superscript $w$ since $x$ may be boundary points for different cells in the same level. We will drop it when no confusion occurs. For $f\in dom\Delta_\mu$, the sum of all normal derivatives of $f$ at $x$ must vanish. This is called the \textit{match condition}. In general, it is necessary and sufficient for gluing together two functions whose Laplacian is defined on neighboring cells.

Also at $x=F_wq_l$, there is a transverse derivative
$$\partial_T^w f(x)=\lim_{m\rightarrow\infty}r_w^{-1}\lambda^{-m}(f(F_wF_l^mq_{l+1})-f(F_wF_l^mq_{l-1}))$$
if the limit exists. For $f\in dom\Delta_\mu$, the transverse derivatives at a point $x$ with respect to different cells may be unrelated.

There are scaling identities for localized derivatives.
$$\partial_n^w f(F_wq_l)=r_w^{-1}\partial_n(f\circ F_w)(q_l),$$
and
$$\partial_T^w f(F_wq_l)=r_w^{-1}\partial_T(f\circ F_w)(q_l).$$

Let $x\in V_*\setminus V_0$. Suppose $m_0$ is the first value for which $x\in V_{m_0}$. We say $x$ is a \textit{junction vertex} if there is at least two different $m_0$-cells containing $x$, i.e., $x$ has at least two different representations $x=F_wq_l$ with $|w|=m_0$. Otherwise, we call $x$  a \textit{nonjunction vertex}, which has exactly one representation $x=F_wq_l$. For both the two different types of vertices, there is a canonical system of \textit{neighborhoods} for each $x$. 
On each certain neighborhood, there is a space of \textit{local multiharmonic functions}. Our definition is slightly different from the definition in {[}S3{]}.

\textbf{Definition 2.1.}

\textit{(a) For $x\in V_m\setminus V_0$, define the m-neighborhood of $x$ as 
\begin{equation*}
U_{m}(x)=\bigcup \{F_wK|x\in F_wK,|w|=m\}.
\end{equation*}
Write $U(x)=U_{m_0}(x)$ for the sake of simplicity, which obviously is the largest one. The boundary of the m-neighborhood $U_m(x)$ is
\begin{equation*}
\partial U_m(x)=\{y\in V_m|y\sim_mx\}. 
\end{equation*}}

\textit{(b) On each $U_m(x)$, define local $(n+1)$-harmonic functions to be those functions $h$ on $U_m(x)$, with $h\circ F_w \in \mathcal{H}_n$ for each $w$, and $\Delta_\mu^i h$ satisfying the matching conditions  at $x$ for all $0\leq i\leq n$. (If $x$ is a nonjunction vertex, we say the matching conditions means $\partial_n \Delta_\mu^i h(x)=0$ for all $0\leq i\leq n$). Write the space of all such functions $\mathcal{H}_n(U_m(x))$.}

Here we remark that our notations differ from that in {[}S3{]} when $x$ is a nonjunction vertex. In our setting, we always view $x$ as an inner point in $U_m(x)$.

Let $W(x)$ denote the set of words of length $m_0$ such that there is a $q_l$ with $x=F_wq_l$. Call $\# W(x)$ the \textit{order} of $x$. Obviously, $\# W(x)\geq 2$ when $x$ is a junction vertex, while $\# W(x)=1$ when $x$ is a nonjunction vertex.  It is easy to verify that the dimension of the space $\mathcal{H}_n(U_m(x))$ is exactly $(n+1)(N_0-1) \# W(x)$ for any $x\in V_m\setminus V_0$. 

For convenience, we always sort the elements in $W(x)$ in lexicographical order. We use $F_x$ to denote the contractive mapping on $U(x)$ with $F_x(y)=F_wF_lF_w^{-1}(y)$ for $y\in F_wK$ and $w\in W(x)$. It is easy to see that $F_x(U_m(x))=U_{m+1}(x)$.

More generally, for a simple connected set $A=\bigcup_{w\in \Omega(A)}F_wK$, where $\Omega(A)$ is a finite set of words, we may define the \textit{boundary} $\partial A$ to be the vertices satisfying

1) $y=F_wq_l,w\in \Omega(A),l=0,1,2$.

2) $y\in V_0$ or $U_m(y)$ is not a subset of $A$ for any $m$.

Then, analogous to the global case, we could define the local energy on $A$ as
\begin{equation*}
\mathcal{E}^A(f,g)=\sum_{w\in \Omega(A)} r_w^{-1}\mathcal{E}(f\circ F_w,g\circ F_w).
\end{equation*}
The domain $dom(\mathcal{E},A)$ is the space of continuous functions on $A$ having finite energy, and $dom_0(\mathcal{E},A)$ is the subspace of such functions which vanish at ${\partial A}$. The Laplacian localized to $A$ could be defined by the weak formulation in an analogous way. We denote  $dom(\Delta_\mu,A)$ the domain of $\Delta_\mu $ on $A$. We need to point out that if $f\in dom\Delta_\mu$,  then $f|_A\in dom(\Delta_\mu, A)$.  Additionally, the local multiharmonic function space is denote as 
$
\mathcal{H}_n(A)$ with dimension  $(n+1)\#\partial A$.

The following Gauss-Green's formula holds.
\begin{equation}
\mathcal{E}^A(f,g)=-\int_A \Delta_\mu f g d\mu+\sum_{z\in \partial A}\partial_n f(z)g(z)
\end{equation}
for any $g\in dom(\mathcal{E},A)$ and $f\in dom(\Delta_\mu,A)$.  The matching condition holds at all vertices except the boundary points. For nonjunction vertices, we still use the term ``matching condition'' like that used in Definition 2.1. In other words, nonjunction vertices always are viewed as inner points. It should be careful that there exists the possibility that some boundary points of $A$ may belong to more than one component cells of $A$ simultaneously. It never occurs for $\mathcal{SG}$ case. See Fig. 2.3 for an example simple set in $\mathcal{SG}_3$. Of course, at these boundary points, the matching condition does not need hold.

\begin{figure}[h]
\centering
\includegraphics[width=4cm]{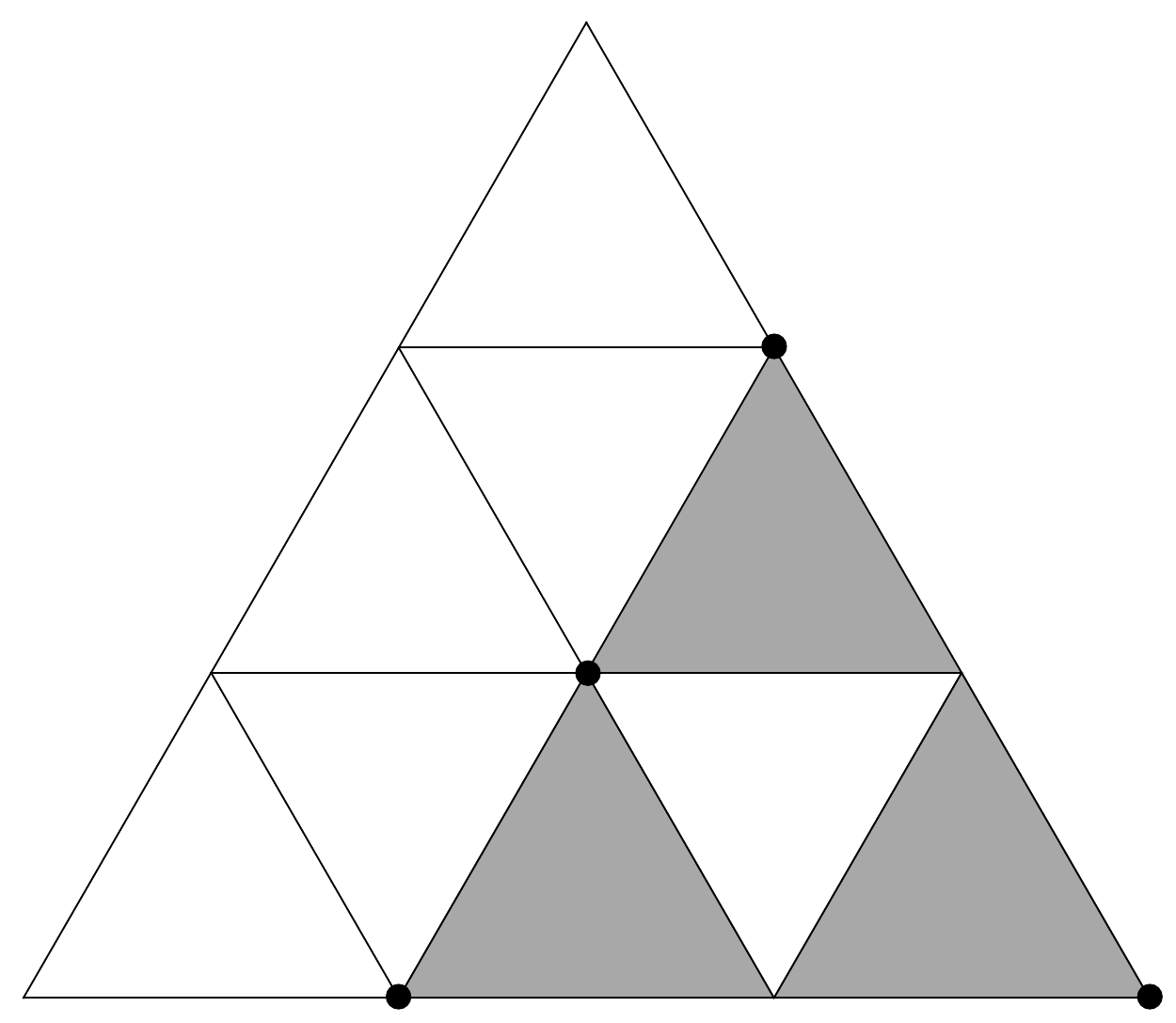}
\begin{center}
 \caption{The shade area is a simple subset in $\mathcal{SG}_3$, with boundary points denoted by dots. Here the center point is a boundary point which belongs to two component cells.}
\end{center}
\end{figure}

\section{Local multiharmonic functions}

There are several bases of $\mathcal{H}_n$ for different purposes. In {[}NSTY{]},  to develop a theory of the local behavior of functions at a single boundary point, a basis of $\mathcal{H}_n$, analogous to the monomials $x^j/j!$ on the unit interval, was described and studied on the Sierpinski gasket $\mathcal{SG}$. This could be easily extended to fractals whose all structures posses full $D3$ symmetry, as follows. Throughout this section, we drop the subscript $\mu$ of $\Delta_\mu$ for simplicity. 

\textbf{Definition 3.1.} \textit{Fix a boundary point $q_l$. The monomials $Q_{jk}^{(l)}$ for $k=1,2,3$ and $0\leq j\leq n$  in $\mathcal{H}_n$ are  the multiharmonic functions satisfying
\begin{eqnarray*}
&&\Delta^i Q_{jk}^{(l)}(q_l)=\delta_{ij}\delta_{k1},\\
&&\partial_n\Delta^i Q_{jk}^{(l)}(q_l)=\delta_{ij}\delta_{k2},\\
&&\partial_T\Delta^i Q_{jk}^{(l)}(q_l)=\delta_{ij}\delta_{k3}.
\end{eqnarray*} }

It is easy to verify that these monomials form a basis of $\mathcal{H}_n$ with the dimension $3(n+1)$. These monomials are related by the following identity, 
$$\Delta Q_{jk}^{(l)}=Q_{(j-1)k}^{(l)}.$$ By the $D3$ symmetry, $Q_{jk}^{(l)}$ for different $l$'s are same under simply rotations. $Q_{j1}^{(l)}$ and $Q_{j2}^{(l)}$ are symmetric while  $Q_{j3}^{(l)}$ is skew-symmetric with respect to the reflection symmetry that fixes $q_l$ and interchanges the other two boundary points. Moreover, the monomials satisfy the following self-similar identities that describe the decay ratios of these functions near $q_l$.

\begin{eqnarray}
&&Q_{j1}^{(l)}(F_l^mx)=\rho^{jm}Q_{j1}^{(l)}(x),\\
&&Q_{j2}^{(l)}(F_l^mx)=r^m\rho^{jm}Q_{j2}^{(l)}(x),\\
&&Q_{j3}^{(l)}(F_l^mx)=\lambda^m\rho^{jm}Q_{j3}^{(l)}(x).
\end{eqnarray}

Denote
\begin{equation}
\alpha_j=Q^{(0)}_{j1}(q_1),\beta_j=Q_{j2}^{(0)}(q_1),\gamma_j=Q_{j3}^{(0)}(q_1),
\end{equation}
for $j\geq 0$. In  {[}NSTY{]}, one can find an elaborate recursion algorithm  of these numbers on the Sierpinski gasket $\mathcal{SG}$. An important observation is that all these numbers are not equal to $0$. The calculation there is quite technical and is hard to be extended  to the general $D3$ case. However, we still could verify that $\alpha_j, \beta_j, \gamma_j$ are never equal to $0$ for some typical fractals with full $D3$ symmetric structures, for example, the level $3$ Sierpinski gasket $\mathcal{SG}_3$ and the hexagasket $\mathcal{HG}$. 
We thus make a following technical assumption.

\textbf{Assumption 3.2.} \textit{All the numbers $\alpha_j, \beta_j$ and $\gamma_j$ are not equal to $0$.}

We will give the calculations of $\alpha_j,\beta_j$ and $\gamma_j$ for $\mathcal{SG}$, $\mathcal{SG}_3$ and $\mathcal{HG}$ in the Appendix, by using a new  algorithm modified from that  in {[}NSTY{]}, which seems to be shorter and more direct. 

We need to extend the above definitions and discussions to all vertices in $V_*$. Naturally, we have the following localized version of monomials, which will play an essential role in answering both the two questions listed in the introduction section.

\textbf{Definition 3.3.} \textit{Fix a vertex $x\in V_*\setminus V_0$. The monomials $P_{jk}^w$ in $\mathcal{H}_n(U(x))$ for $k=1,2,3$, $0\leq j\leq n$ and $w\in W(x)$ are the local multiharmonic functions satisfying}
\begin{eqnarray*}
&&\Delta^i P_{jk}^w(x)=\delta_{ij}\delta_{k1},\\
&&\partial_n^{w''}\Delta^i P_{jk}^{w}(x)=\delta_{ij}\delta_{k2}\delta_{ww''}-\delta_{ij}\delta_{k2}\delta_{w'w''},\\
&&\partial_T^{w''}\Delta^i P_{jk}^{w}(x)=\delta_{ij}\delta_{k3}\delta_{ww''},
\end{eqnarray*}
\textit{ where $w'$ is the next word to $w$ in $W(x)$ in lexicographical order.}

\textbf{Remark 1.}  For $k=1$, the superscript $w$ is unnecessary, and we may not add it when discuss $P_{j1}^{w}$ seperately. For nonjunction vertices, there are no monomials in $k=2$ case. 

\textbf{Remark 2.} 
It is easy to check that $\{P_{jk}^{w}|_{U_m(x)}\}$  forms a basis of $\mathcal{H}_{n}(U_m(x))$.

Similar to (3.1)-(3.3), the following self-similar identities hold
\begin{eqnarray*}
&&P_{j1}(F_x^my)=\rho^{jm}P_{j1}(y),\\
&&P_{j2}^w(F_x^my)=r^m\rho^{jm}P_{j2}^w(y),\\
&&P_{j3}^w(F_x^my)=\lambda^m\rho^{jm}P_{j3}^w(y),
\end{eqnarray*}
which describe the decay behaviors of these monomials near $x$.  It is apparent that these monomials have symmetric properties analogous to the global case.

Denote $R_i$ the rotations in $D3$ symmetric group, with $R_i(q_l)=q_{l+i}$ (cyclic notation).

\textbf{Definition 3.4.} \textit{Fix a vertex $x\in V_*\setminus V_0$. For $k=1,2,3$, let $P_k$ be a linear projection from  $\mathcal{H}_n (U(x))$ into itself satisfying
\begin{eqnarray}
&&\Delta^i P_k(h)(x)=\delta_{k1}\Delta^i h(x), \\
&&\partial_n^w\Delta^i P_k(h)(x)=\delta_{k2}\partial_n^w\Delta^i h(x),\\
&&\partial_T^w\Delta^i P_k(h)(x)=\delta_{k3}\partial_T^w\Delta^i h(x),
\end{eqnarray}
for any $h\in \mathcal{H}_n(U(x))$, $w\in W(x)$, $0\leq i\leq n$. Let $R$ be a linear mapping on $\mathcal{H}_n (U(x))$, defined by
$$
R(h)(y)=(\sum_{w'\in W(x)}r_{w'}^{-1})^{-1}\sum_{w'\in W(x)} r_{w'}^{-1}h\circ F_{w'} \circ R_{l'-l}\circ F_w^{-1}(y),
$$
for $y\in F_wK, w\in W(x)$, for any $h$ in $\mathcal{H}_n (U(x))$.} 

Clearly, $P_k(h)$ is a linear combination of the monomials $P_{jk}^{w}$, and it is easy to check that $$P_1(h)+P_2(h)+P_3(h)=h.$$ As for $R$, roughly speaking, it is an operator on $\mathcal{H}_n (U(x))$ which first rotates variables around $x$, then takes mean values with weights proportional to $r_{w'}^{-1}$.

Let $g_x$ be the local symmetry in $U(x)$, which fixes $F_wq_l$ and permutes the other two boundary points of $F_wK$ for each $w\in W(x)$.

\textbf{Theorem 3.5.} \textit{Assume $x\in V_*\setminus V_0$ and $h\in\mathcal{H}_n(U(x))$, then the following identities hold
\begin{eqnarray}
&&P_1(h)=\frac{1}{2}(R(h)\circ g_x+R(h)), \\
&&P_2(h)=\frac{1}{2}(h+h\circ g_x-R(h)\circ g_x-R(h)), \\
&&P_3(h)=\frac{1}{2}(h-h\circ g_x).
\end{eqnarray}}
\textit{Proof.}  The following equalities are consequences of symmetric definitions of $\Delta$, $\partial_n$ and $\partial_T$.
\begin{eqnarray*}
&&\Delta^i R(h)(x)=\Delta^i h(x), \\
&&\partial_n^w\Delta^i R(h)(x)=0,
\end{eqnarray*}
and
\begin{eqnarray*}
&&\Delta^i h\circ g_x(x)=\Delta^i h(x),\\
&&\partial_n^w\Delta^i h\circ g_x(x)=\partial_n^w\Delta^i h(x),\\
&&\partial_T^w\Delta^i h\circ g_x(x)=-\partial_T^w\Delta^i h(x),
\end{eqnarray*}
which hold for $w\in W(x)$, $0\leq i\leq n$, for any $h\in\mathcal{H}_n(U(x))$. These yield the result of the theorem. $\Box$

\textbf{Corollary 3.6.} \textit{Assume $x\in V_*\setminus V_0$ and $h\in\mathcal{H}_n(U(x))$, then for each $m\geq m_0$,  $h|_{\partial U_m(x)}=0$ if and only if $P_k(h)|_{\partial U_m(x)}=0$ for $k=1,2,3$.}\\

In the rest of this section, we give an application of the local monomials, to show that the higher order weak tangents of smooth functions $f$ at any fixed vertex, could be expressible as limits of local multiharmonic functions that agree with $f$ at the boundary of $U_m(x)$. This is  Question 1 which we want to solve in this paper. To be more precise, we need the following definition of \textit{higher order weak tangents}.

\textbf{Definition 3.7.} \textit{Let $x$ be a vertex in $V_*\setminus V_0$ and $f$ a function defined in a neighborhood of $x$. We say that an $(n+1)$-harmonic function $h$ is a weak tangent of order $n+1$ of $f$ at $x$ if } 
\begin{equation}
(f-h)|_{\partial U_m(x)}=o((\rho^{n}r)^m)
\end{equation}
and
\begin{equation}
(f-h-(f-h)\circ g_x)|_{\partial U_m(x)}=o((\rho^{n}\lambda)^m).
\end{equation}

\textbf{Theorem 3.8.} \textit{Assume Assumption 3.2 holds. Let $x\in V_*\setminus V_0$. Then the following two conclusions hold.}

\textit{(a)An $(n+1)$-harmonic function $h$ on $U(x)$ is uniquely determined by the values $h|_{\partial U_{m+i}},0\leq i\leq n$, and any such values may be freely assigned.}

\textit{(b)Let $f$ be a continuous function defined in a neighborhood of $x$, and assume $f$ has a weak tangent of order n+1 at $x$, denoted by $h$. Let $h_m$ be the (n+1)-harmonic function defined in $U(x)$, assuming the same values as $f$ at the boundary points of $U_{m+i}(x)$ for all $0\leq i\leq n$.
Then $h_m$ converges to $h$  uniformly on $U(x)$.}

\textbf{Remark 1.} This theorem extends the previous result in {[}CQ,S3{]} for the $1$-order tangents and $1$-order harmonic functions. For the nonjunction vertices, there is an implicit restriction that $f$, $h$ and $h_m$ should satisfy the equation $\partial_n \Delta^iu(x)=0$ for all $0\leq i\leq n$, since we always view $x$ as an inner point in $U(x)$.

\textbf{Remark 2.}  There are some sufficient conditions to ensure the existence of the weak tangents. One can find more detailed discussion on the weak tangents (and \textit{tangents}, \textit{strong tangents}) in {[}S3{]}.

\textit{Proof of Theorem 3.8.} (a) The map from $\mathcal{H}_n(U(x))$ to the values $h|_{\partial U_{m+i}(x)},i=0,1,...,n$ is obviously a linear map, and the dimemsion of $\mathcal{H}_n(U(x))$ is $2(n+1)\#W(x)$, which is exactly equal to $\#\bigcup_{0\leq i\leq n}\partial U_{m+i}(x)$. Thus to proof (a), we only need to show that the map is injective.

Fix a word $w\in W(x)$ with $x=F_wq_l$. Let $h\in \mathcal{H}_n(U(x))$. For $k=1,2,3$, notice that $P_k(h)\circ F_w $ is a linear combination of $Q_{jk}^{(l)}$,  denote the combination coefficients of $Q_{jk}^{(l)}$ by $a_{jk}^{w}$. We have the following equalities 
\begin{equation}
\begin{aligned}
P_k(h)(F_x^{m-|w|+i}F_wq_{l+1})&=P_k(h)(F_wF_l^{m-|w|+i}q_{l+1})\\
&=\sum_{j=0}^{n} a_{jk}^wQ_{jk}^{(l)}(F_l^{m-|w|+i}q_{l+1})\\
&=\sum_{j=0}^{n} (A_k)_{ij}a_{jk}^w,
\end{aligned} 
\end{equation}
for
$(A_k)_{ij}=Q_{jk}^{(l)}(F_l^{m-|w|+i}q_{l+1})=Q_{jk}^{(l)}(F_l^{m'+i}q_{l+1}),$ where we denote $m'=m-|w|$ for convenience. Thus the $(n+1)\times (n+1)$ matrix $A_k$ induce a linear map from  $\{a_{jk}^w\}_j$ to the values $\{P_k(h)(F_x^{m-|w|+i}F_wq_{l+1})\}_i$. 

We now show that the matrix $A_k$ is invertible for $k=1,2,3$. Denote by $\Gamma^{(n)}$ an $(n+1)\times(n+1)$ matrix with
\begin{equation*}
  \Gamma^{(n)}=\begin{pmatrix}
    1&1&1&...& 1\\1&\rho&\rho^2&...&\rho^n\\\vdots&\vdots&\vdots& &\vdots\\1&\rho^n&\rho^{2n}&...&\rho^{n^2}
  \end{pmatrix},
\end{equation*}
which is obviously invertible. Then by using the self-similar identities $(3.1)-(3.3)$, we have
\begin{eqnarray*}
&&(A_1)_{ij}=\rho^{m'j+ij}\alpha_j=(\Gamma^{(n)})_{ij}\rho^{m'j}\alpha_j,\\
&&(A_2)_{ij}=r^{m'+i}\rho^{m'j+ij}\beta_j=r^{m'+i}(\Gamma^{(n)})_{ij}\rho^{m'j}\beta_j,\\
&&(A_3)_{ij}=\lambda^{m'+i}\rho^{m'j+ij}\gamma_j=\lambda^{m'+i}(\Gamma^{(n)})_{ij}\rho^{m'j}\gamma_j,
\end{eqnarray*}
which can be rewritten in matrix notation,
\begin{eqnarray*}
&&A_1=\Gamma^{(n)}diag(\alpha_0,\rho^{m'}\alpha_1,...,\rho^{m'n}\alpha_n),\\
&&A_2=diag(r^{m'},...,r^{m'+n}) \Gamma^{(n)}diag(\beta_0,\rho^{m'}\beta_1,...,\rho^{m'n}\beta_n),\\
&&A_3=diag(\lambda^{m'},...,\lambda^{m'+n}) \Gamma^{(n)}diag(\gamma_0,\rho^{m'}\gamma_1,...,\rho^{m'n}\gamma_n),
\end{eqnarray*}
from which  it is obviously that all the matrices $A_1,A_2,A_3$ are invertible. 

The above discussion shows that $P_k(h)$ vanishs at $\partial U_{m+i}(x)$ if and only if $P_k(h)=0$. According to  Corollary 3.6, $h$ vanishs at $\partial U_{m+i}(x)$ if and only if all $P_k(h)$ vanishs at $\partial U_{m+i}(x)$. Thus we have proved (a). 

(b) We need to study the (n+1)-harmonic functions $h-h_m$. Notice that formula (3.13) still holds for $h-h_m$.

For $k=1$, we have $P_1(h_m-h)(F_x^{m'+i}F_wq_{l+1})=\sum_{j=0}^n (A_1)_{ij}a_{j1}$. Thus
\begin{equation*}
\begin{aligned}
a_{j1}&=\sum_{i=0}^{n}(A_1^{-1})_{ji}P_1(h_m-h)(F_x^{m'+i}F_wq_{l+1})\\&=\alpha_j^{-1}\rho^{-jm'}\sum_{i=0}^{n}(\Gamma^{(n)})^{-1}_{ji}P_1(h_m-h)(F_x^{m'+i}F_wq_{l+1}).
\end{aligned}
\end{equation*}
According to $(3.11)$, by using Theorem 3.5, we have
$P_1(h-h_m)|_{\partial U_{m+i}(x)}=o(r^m\rho^{mn})$, which gives that
\begin{equation*}
a_{j1}=o(r^m\rho^{m(n-j)}).
\end{equation*}

For $k=2$,  a similar discussion shows that 
\begin{equation*}
\begin{aligned}
a_{j2}^w&=\sum_{i=0}^n (A_2^{-1})_{ji}P_2(h_m-h)(F_x^{m'+i}F_wq_{l+1})\\
&=\beta_j^{-1}\rho^{-jm'}\sum_{i=0}^{n}(\Gamma^{(n)})^{-1}_{ji}r^{-m'-i}P_2(h_m-h)(F_x^{m'+i}F_wq_{l+1}).
\end{aligned}
\end{equation*}
According to (3.11), still using Theorem 3.5, we have
$P_2(h-h_m)|_{\partial U_{m+i}(x)}=o(r^m\rho^{mn})$, and hence
\begin{equation*}
a_{j2}^w=o(\rho^{m(n-j)}).
\end{equation*}

For $k=3$,  the same argument yields that
\begin{equation*}
\begin{aligned}
a_{j3}^w&=\sum_{i=0}^n (A_3^{-1})_{ji}P_3(h_m-h)(F_x^{m'+i}F_wq_{l+1})\\
&=\gamma_j^{-1}\rho^{-jm'}\sum_{i=0}^{n}(\Gamma^{(n)})^{-1}_{ji}\lambda^{-m'-i}P_3(h_m-h)(F_x^{m'+i}F_wq_{l+1}).
\end{aligned}
\end{equation*}
Using (3.12) and Theorem 3.5, we can get $P_3(h-h_m)|_{\partial U_{m+i}(x)}=o(\lambda^m\rho^{mn})$, so
\begin{equation*}
a_{j3}^w=o(\rho^{m(n-j)}).
\end{equation*}
Thus we have proved that for $k=1,2,3$, $P_k(h-h_m)$ converges  uniformly to zero on each cell $F_wK$, which yields that $h_m$ converges uniformly to $h$ on $U(x)$.  $\Box$

\section{pointwise formula for the higher order Laplacians}

In this section, we will deal with Question 2. We still restrict to consider the $D3$ symmetric fractals. The subscript $\mu$ of the Laplacian $\Delta_\mu$ is still dropped for simplicity. 

\subsection{Definition of pointwise formula}\par

Analogous to the pointwise formula of the Laplacian, we will show that we can approach the $n$-order Laplacian  by the $n$ times iterating of the renormalized discrete Laplacian, which means 
\begin{equation}
\Delta^n f(x)=\lim_{m\to \infty}\tilde{\Delta}_m^n f(x).
\end{equation} 
Notice that  $\tilde{\Delta}_m^n$ may not be defined on all vertices in $V_m\setminus V_0$ for $n\geq 2$. For example, when $n=2$, it is exact those vertices which are not connected with the boundary $V_0$ having the operation  $\tilde{\Delta}_m^2$ well-defined.

It is convenient to define the following notations.

\textbf{Definition 4.1.}\textit{ For two vertices $x,y\in V_m$, the m-distance $d_m(x,y)$ between them is the minimal number of edges which connect $x$ to $y$ in $G_m$.}

It is easy to check that any vertex satisfying $d_m(x, V_0)\geq n$ has a well-defined $\tilde{\Delta}_m^n f(x)$, and we denote
\begin{equation*}
V^n_m=\{x\in V_m:d_m(x,V_0)\geq n\}
\end{equation*}
the domain of the definition of $\tilde{\Delta}^n_m$. See Fig. 4.1 for $V_2^2$, the domain of $\tilde{\Delta}_2^2$ for  $\mathcal{SG}$.

 \begin{figure}[htbp]
     \centering
     \includegraphics[width=0.4\textwidth]{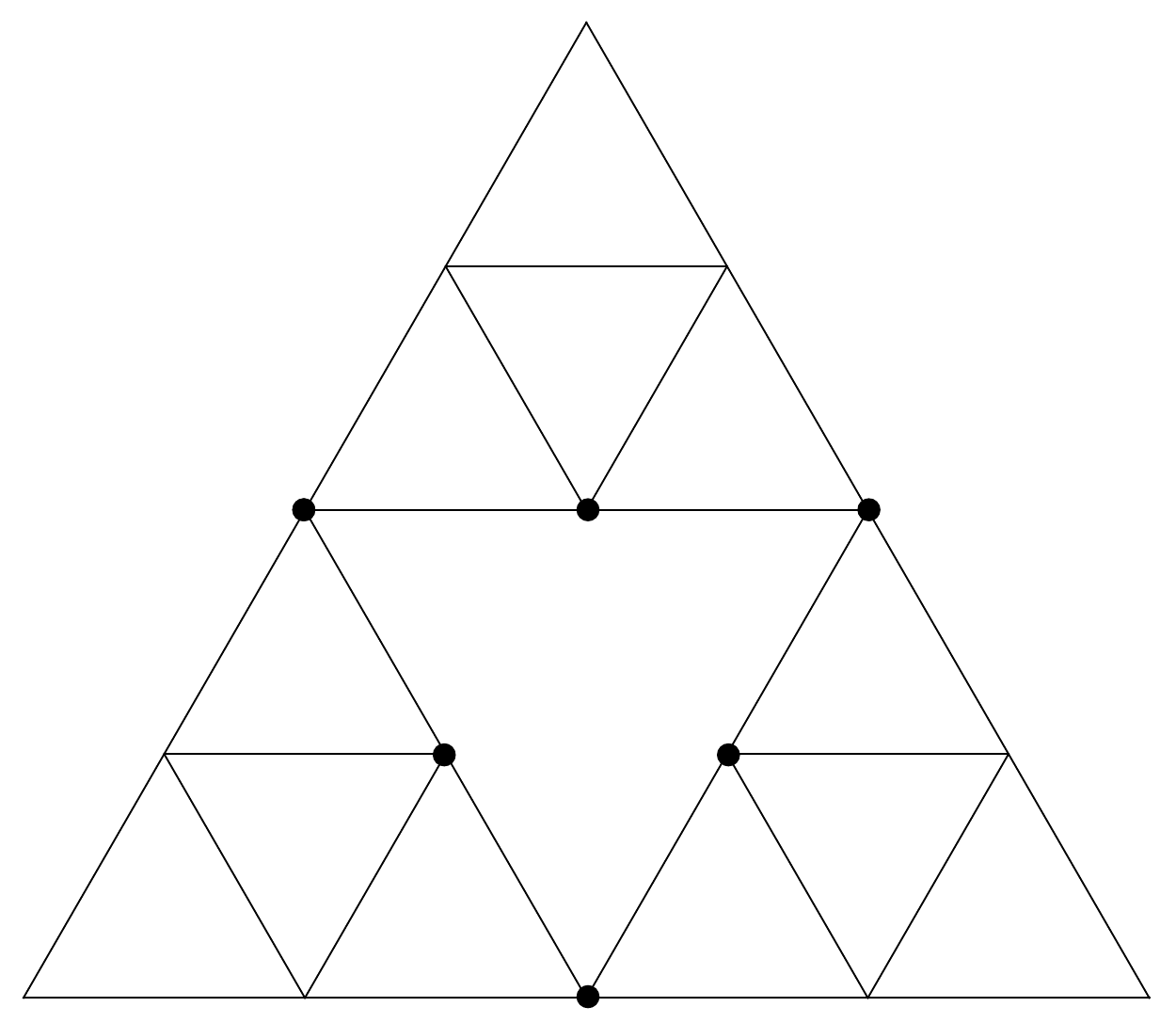}
     \caption{ The domain of $\tilde{\Delta}_2^2$ for $\mathcal{SG}$, denoted by dots. }
\end{figure}

For fixed $x\in V_m^n$, the calculation of $\tilde{\Delta}_m^n f(x)$ involves the values of $f$ at those vertices with $m$-distance to $x$ no more than $n$, which are collected as
\begin{equation}
\begin{aligned}
L_m^n(x)&=\{y\in V_m:d_m(x,y)\leq n\}\\
        &=\bigcup\{L_m^1(y):y\in L_m^{n-1}(x)\}.
\end{aligned}
\end{equation}
The area bounded by these vertices is obviously a neighborhood of $x$, which may be written as $U_m^n(x)$ (see Fig. 4.2), with the following identity holds,
\begin{equation}
U_m^n(x)=\bigcup\{U_m(y):y\in L_{m}^{n-1}(x)\} .
\end{equation}

 \begin{figure}[htbp]   
     \centering
     \includegraphics[width=0.4\textwidth]{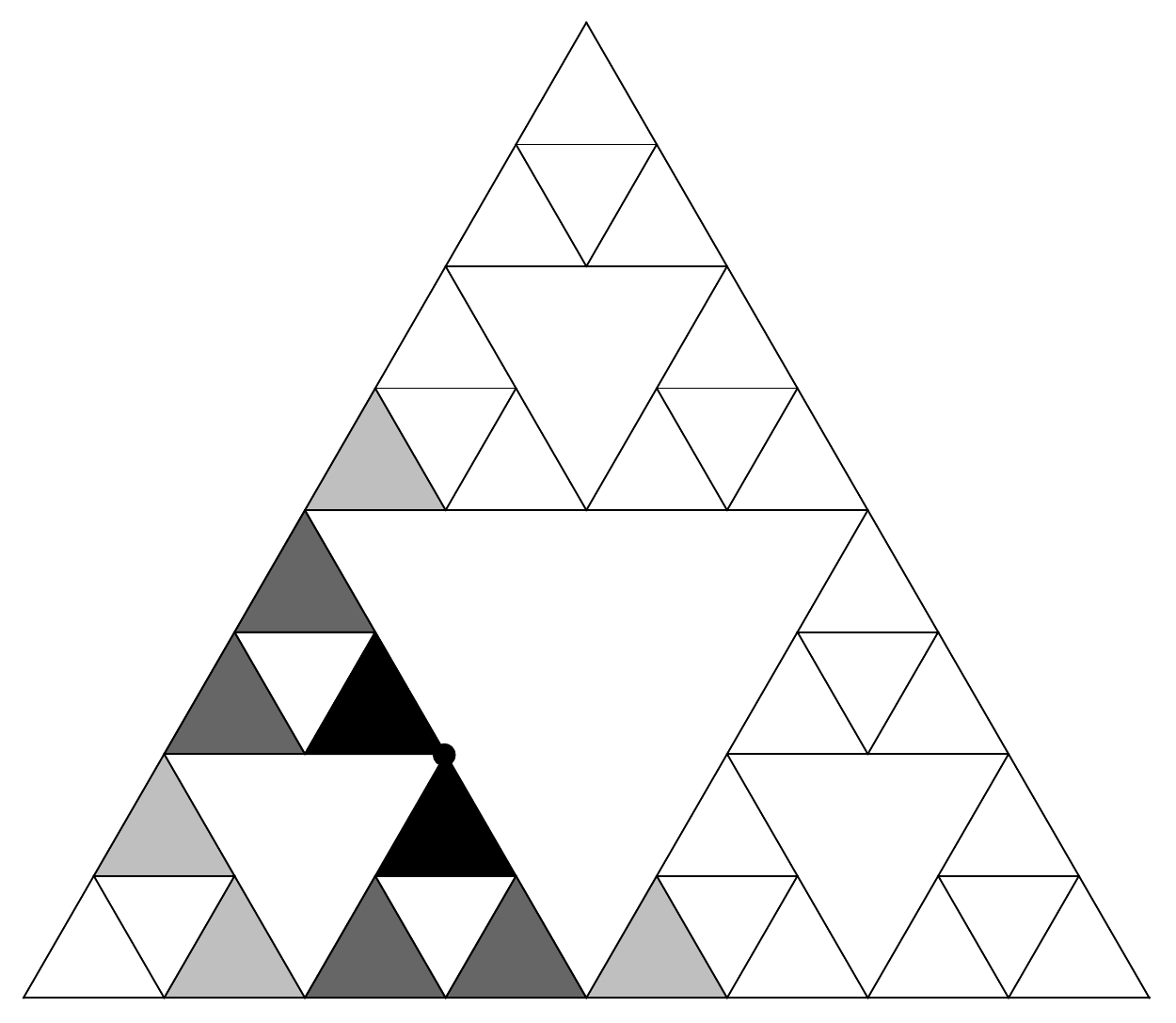}\qquad\qquad
     \includegraphics[width=0.4\textwidth]{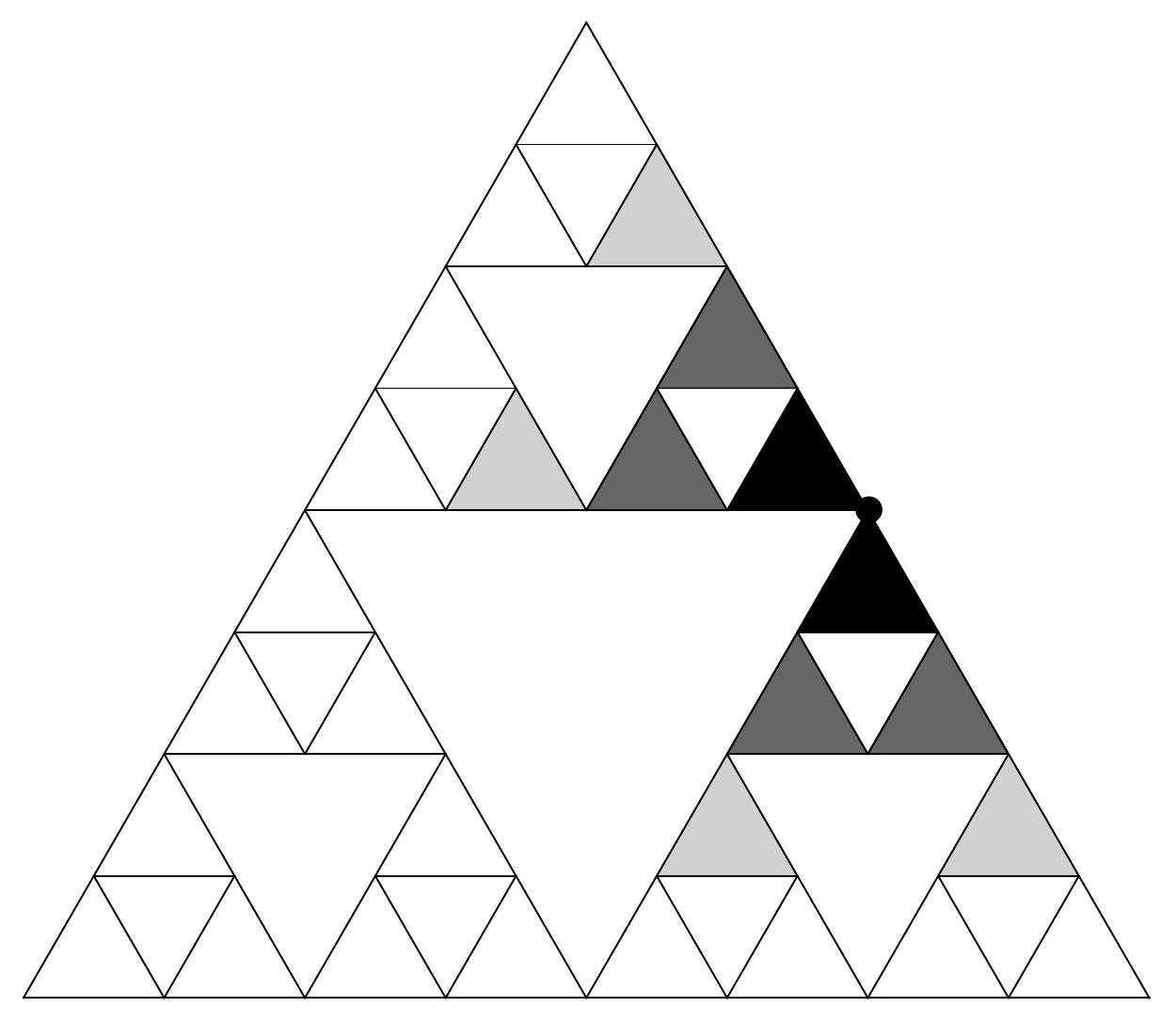}
     \caption{ Some examples of $U_3^n(x)$ with $n\leq 3$ in the Sierpinski gasket.}
\end{figure}

It is natural that the \textit{boundary} of $U_m^n(x)$ is
\begin{equation*}
\partial U_m^n(x)=L_m^n\setminus\{y\in V_m\setminus V_0: U_m(y)\subset U_m^n(x)\}, 
\end{equation*}
which is consistent with the boundary of $U_m(x)$ and the boundary of simple set $A$ as introduced in Section 2. In our setting, nonjunction vertices still always be viewed as inner points. It is easy to check that $\partial U_{m}^n(x)\subset L_m^n(x)\setminus L_m^{n-1}(x)$. It should be careful that there indeed exist vertices that belong to $ L_m^n(x)\setminus L_m^{n-1}(x)
$, which are not  boundary points of $U_m^n(x)$. For example, it is the case when we choose $x$ to be the bottom dotted vertex in Fig 4.1 for $\mathcal{SG}$ for $n=m=2$.

\textbf{Remark.}  The shape of $U_m^n(x)$ varies for $x$ in $V_m^n$ and $m\geq 0$. We could give  a classification of them. Let $x\in V_m^n$ and $y\in V_{m'}^n$. We say $U_m^n(x)$ and $U_{m'}^n(y)$ belong to a same \textit{type}  if there exists some mapping $F$ which is a combination of rotations, reflections and scalings such that $FU_m^n(x)=U_{m'}^n(y)$.

We conclude that there are only \textit{finite types} of $U_m^n(x)$ for any fixed $n$. In fact, the second equality  of (4.2) shows that if there are finite types of $U_m^{n-1}(x)$, then the types of $U_m^n(x)$ is also finite. This observation will be useful in the proof of the uniform convergence of the pointwise formula. See Fig. 4.3 for the total types of $U_m^2(x)$ in $\mathcal{SG}$.

 \begin{figure}[htbp]
     \centering
     \includegraphics[width=0.28\textwidth]{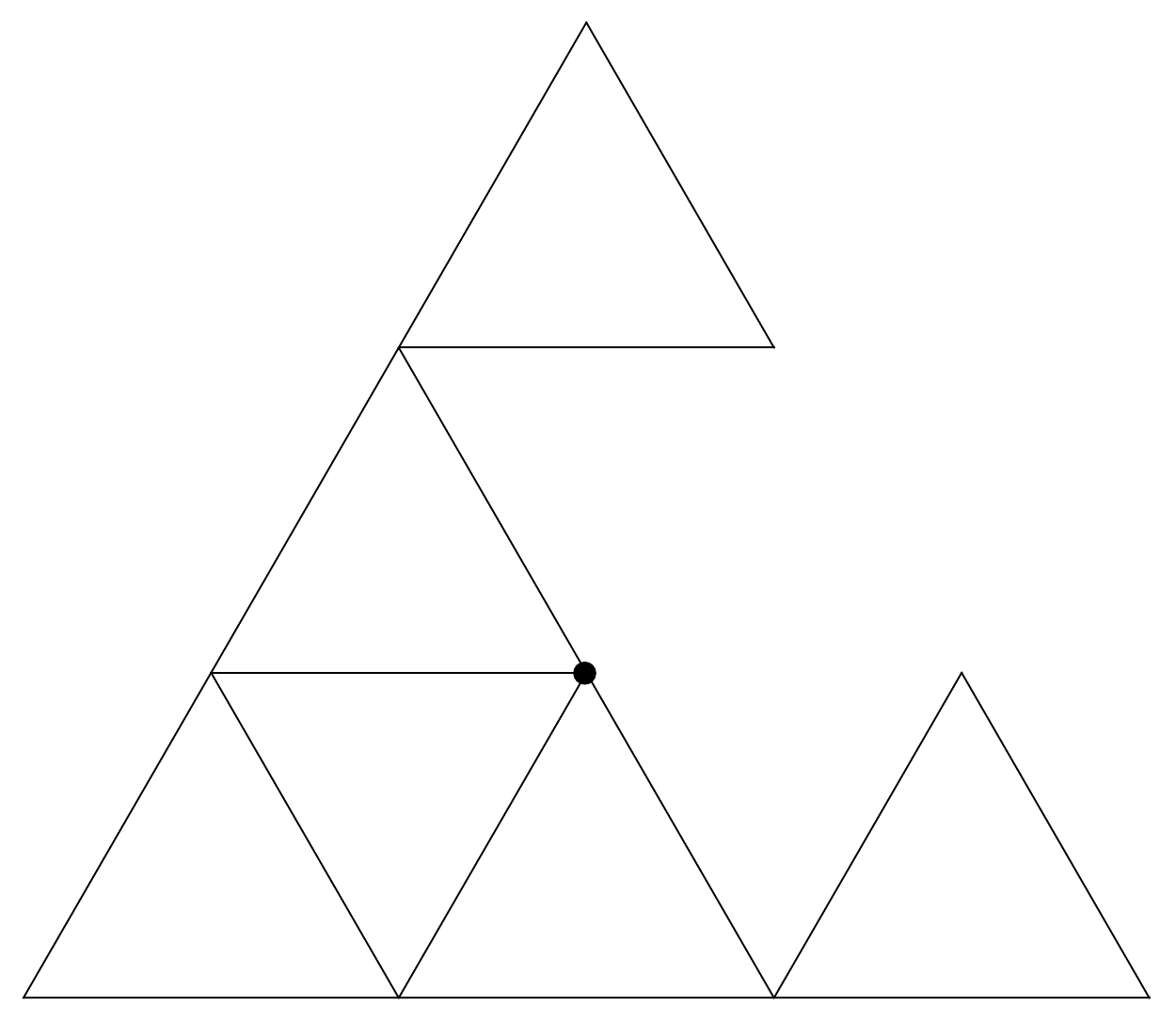}
     \includegraphics[width=0.38\textwidth]{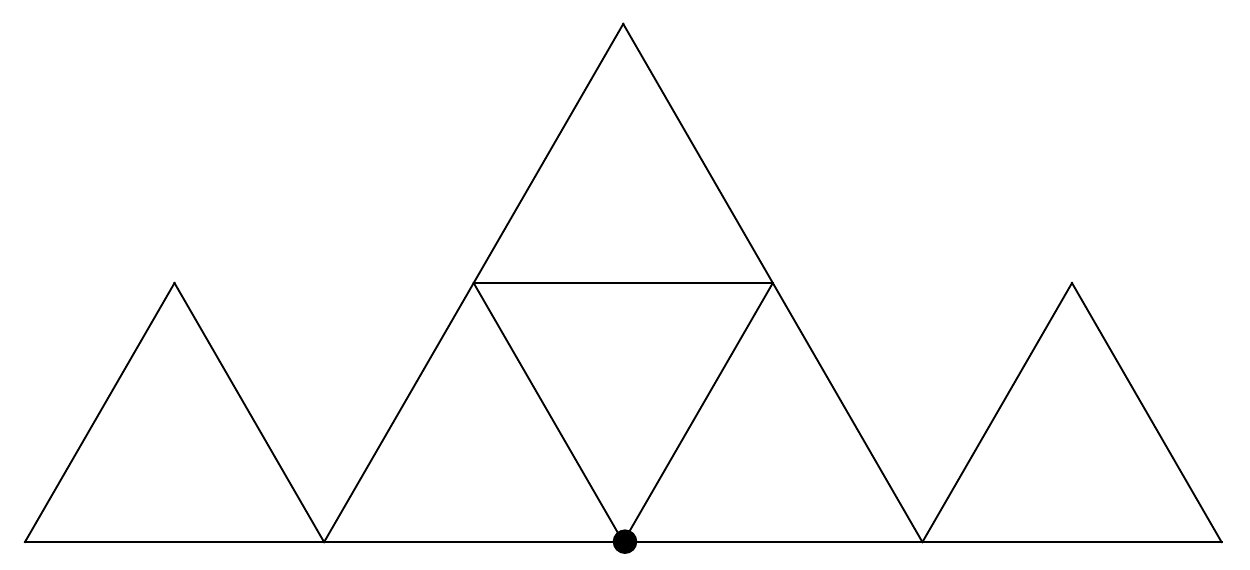}
     \includegraphics[width=0.28\textwidth]{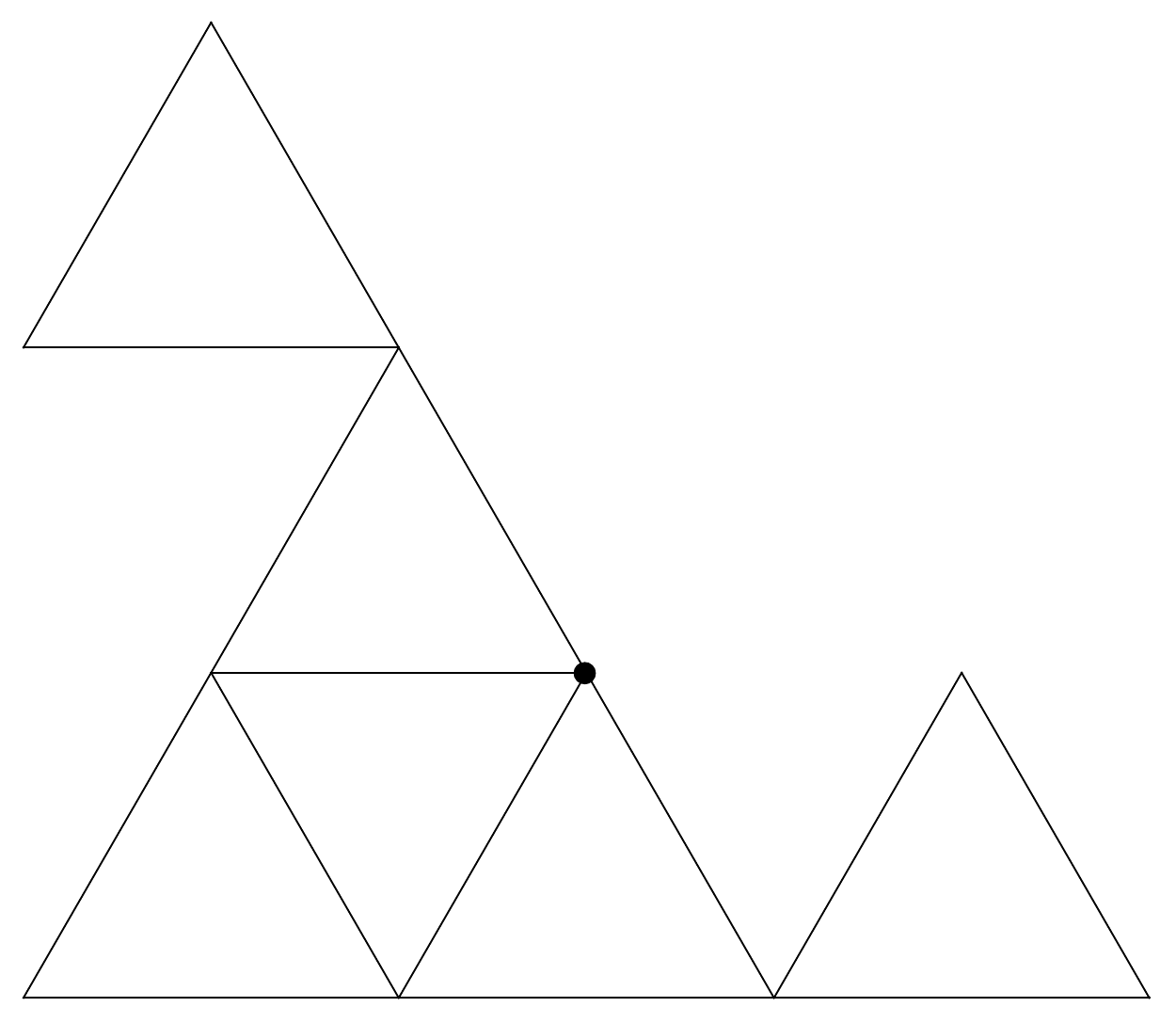}
     \includegraphics[width=0.38\textwidth]{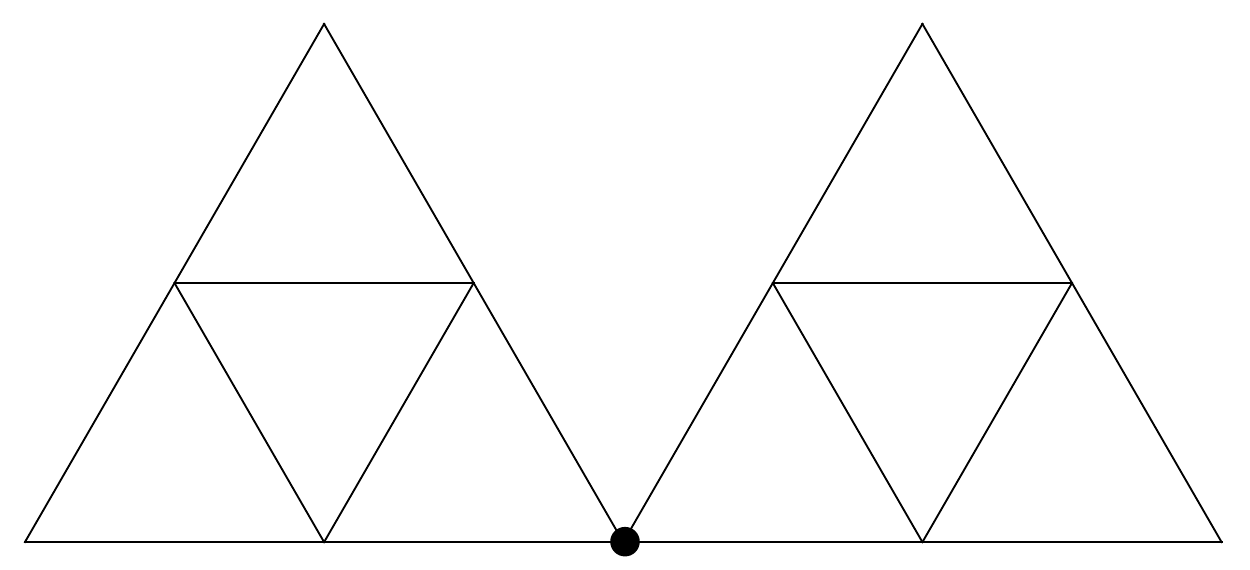}
     \caption{The total types of $U_m^2(x)$ in $\mathcal{SG}$.}
\end{figure}

The following theorem is an answer of Question 2, which will be proved in the subsequent  subsections.

\textbf{Theorem 4.2.} \textit{(a) Assume $f\in dom(\Delta^n)$. Then the pointwise formula (4.1) holds with the limit uniformly on $V_*\setminus V_0$.
(b) Conversely, let $f\in C(K)$ and the right side of (4.1) converges uniformly to a continuous function $u$ on $V_*\setminus V_0$. Then $f\in dom(\Delta^n,K\setminus V_0)$ with $\Delta^n f=u$ on $K\setminus V_0$.}

Before proving, we remark that it looks that the (b) part of this theorem does not involve the whole information of the function $f$, and it has something to do with the existence of harmonic functions with singularities at boundary points (See more explanation on point singularities in [BSSY]). In fact, the conclusion is equivalent to  that for any $g\in dom(\Delta^n)$ with $\Delta^n g=u$, we have $f-g$ may be an $n$-harmonic function with singularity. However, if there is no multiharmonic function with singularity, for example,  the unit interval case, we can say that $f\in dom(\Delta^n)$ already.

\subsection{Proof of Theorem 4.2(a).}
We will take two steps to prove part (a) of Theorem 4.2. First, we deal with those functions which are local $(n+1)$-harmonic near $x$ with $x\in V_m^n$, to get that $(4.1)$ holds without taking the limit. Then, we prove the result for general functions in $dom(\Delta^n)$.

\textbf{Lemma 4.3.} \textit{Let $x$ be a vertex in $V_m\setminus V_0$,  $h$ be an $(n+1)$-harmonic function in $\mathcal{H}_n(U_m(x))$. Then the following equality 
holds
\begin{equation}
\tilde{\Delta}_m h(x)=\sum_{j=1}^{n} \rho^{m(j-1)}\alpha_1^{-1}\alpha_j \Delta^j h(x).
\end{equation}
In particular,  $\alpha_1=1/6$.}

\textit{Proof}. 
Fix a word $w\in W(x)$ with $x=F_wq_l$. Note that $P_1(h\circ F_x^{m-|w|})\circ F_w$ is a linear combination of monomials $Q_{j1}^{(l)}$. In fact
\begin{equation*}
P_1(h\circ F_x^{m-|w|})\circ F_w=\sum_{j=0}^{n} \rho^{mj}\Delta^jh(x)Q_{j1}^{(l)},
\end{equation*}
by comparing the values at $x$ when applying  $\Delta^j$ on both sides. Thus we have
\begin{equation*}
\begin{aligned}
P_1(h\circ F_x^{m-|w|})(F_wq_{l+1})&=\sum_{j=0}^n \rho^{mj}\Delta^jh(x)Q_{j1}^{(l)}(q_{l+1})\\
&=\sum_{j=0}^n \rho^{mj}\alpha_j \Delta^j h(x).
\end{aligned}
\end{equation*}
On the other hand, according to (3.8), 
$$\begin{aligned}
P_1(h\circ F_x^{m-|w|})(F_wq_{l+1})&=(2\sum_{w'\in W(x)}r_{w'}^{-1})^{-1}\sum_{y\in F_{w'}K,y\sim_m x}r_{w'}^{-1}h(y)\\
&=(\sum_{y\sim_m x}c_{xy})^{-1}\sum_{y\sim_m x}c_{xy}h(y).
\end{aligned}
$$
Thus,
\begin{equation*}
\begin{aligned}
\tilde{\Delta}_m h(x)&=\frac{\sum_{y\sim_m x}c_{xy}}{\int \psi^m_{x}d\mu}\left(P_1(h\circ F_x^{m-|w|})(F_wq_{l+1})-h(x)\right)\\&=\frac{2\sum_{w\in W(x)}r_w^{-1}r^{-(m-|w|)}}{\int \psi^m_{x}d\mu}\sum_{j=1}^{n} \rho^{mj}\alpha_j\Delta^j h(x)\\
&= 6\rho^{-m} \sum_{j=1}^{n}\rho^{mj}\alpha_j\Delta^j h(x).
\end{aligned}
\end{equation*}
The second equality above comes from the fact that $\alpha_0$ always equal to $1$. From the arbitrariness of $h$, if we choose $h$ to satisfy  $\Delta h=1$ in the above equality, this gives 
$$\tilde{\Delta}_m h(x)=6\alpha_1$$ for all $m$. By passing $m$ to infinity, we get that $\alpha_1=1/6$. Thus we have proved the lemma.  $\Box$

\textbf{Remark.} We use the decomposition of $h$ based on the monomials in the above proof, which requires the harmonic extension matrices $M_0,M_1,M_2$ to be non-degenerate. However, this requirement is not necessary essentially. In fact, we could alternatively start from an ``easy'' basis, which extends the conclusion to the degenerate cases. We will give this method in Section 6 when discussing the $D4$ symmetric fractals. 

It is interesting that the constant $\alpha_1=1/6$ is universal for all $D3$ symmetric fractals, which is an initial value for the calculations in Appendix.

\textbf{Lemma 4.4.} \textit{let x be a vertex in $V_m^n$ and $h\in \mathcal{H}_n(U_m^n(x))$, then 
\begin{equation*}
\tilde{\Delta}_m^n h(x)=\Delta^n h(x).
\end{equation*}                 }
\textit{Proof}. It is obvious that $h|_{U_m(y)}\in \mathcal{H}_n(U_m(y))$ for any $y$ with $d_m(x,y)\leq n-1$. So we can apply Lemma 4.3 to all points in $L_m^{n-1}(x)$. Thus we have
\begin{equation*}
\begin{aligned}
\tilde{\Delta}^n_m h(x)&=\tilde{\Delta}^{n-1}_m(\tilde{\Delta}_m h)(x)\\
                       &=\tilde{\Delta}^{n-1}_m(\sum_{j=1}^{n} \rho^{m(j-1)}\alpha_1^{-1}\alpha_j \Delta^j h)(x)\\
                       &=\sum_{j=1}^{n}\rho^{m(j-1)}\alpha_1^{-1}\alpha_j\tilde{\Delta}_m^{n-1}(\Delta^jh)(x).
\end{aligned}
\end{equation*}
Since for each $j\geq 1$, $\Delta^jh$ belongs to $\mathcal{H}_{n-1}(U_m^{n-1}(x))$, by a  standard inductive argument, we then have
$$\tilde{\Delta}_m^nh(x)=\sum_{j=1}^{n}\rho^{m(j-1)}\alpha_1^{-1}\alpha_j\Delta^{n+j-1}h(x)=\Delta^{n}h(x).$$
Hence we have proved the lemma. $\Box$

Now for each $x\in V_m^n$, we will define a global function  $\phi_{m,x}^{(n)}$ supported in $U_m^n(x)$, belonging to $dom(\Delta ^{n-1}, U_m^n(x))$, called \textit{$n$-tent function},  which is piecewise $n$-harmonic and  for any $f\in dom(\Delta^n)$, it holds that
\begin{equation*}
\tilde{\Delta}_m^n f(x)=-\mathcal{E}(\phi_{m,x}^{(n)},\Delta^{n-1} f).
\end{equation*}

We will define this function by an inductive argument. 

In fact, when $n=1$, we just need to choose $\phi_{m,x}^{(1)}=(\int \psi_x^m d\mu)^{-1}\psi_x^m$, and we have
\begin{equation*}
\tilde{\Delta}_m f(x)=-\mathcal{E}(\phi_{m,x}^{(1)},f).
\end{equation*}

When $n=2$, first we define
\begin{equation*}
\tilde{\phi}_{m,x}^{(2)}=\sum_{y\sim_m x}(\int \psi_x^m d\mu)^{-1}c_{xy}(\phi_{m,y}^{(1)}-\phi_{m,x}^{(1)}).
\end{equation*}
It is easy to check that
\begin{equation*}
\begin{aligned}
\tilde{\Delta}^2_m f(x)&=\sum_{y\sim_m x} (\int \psi_x^m d\mu)^{-1}c_{xy}\left(\tilde{\Delta}_m f(y)-\tilde{\Delta}_m f(x)\right)\\
&=-\sum_{y\sim_m x}(\int \psi_x^m d\mu)^{-1}c_{xy}\left(\mathcal{E}(\phi_{m,y}^{(1)},f)-\mathcal{E}(\phi_{m,x}^{(1)},f)\right)\\
&=-\mathcal{E}(\tilde{\phi}_{m,x}^{(2)},f).
\end{aligned}
\end{equation*}
Now let 
\begin{equation*}
\phi_{m,x}^{(2)}(\cdot)=-\int_K G_{m,2}(\cdot,z)\tilde{\phi}_{m,x}^{(2)}(z)d\mu(z),
\end{equation*}
where $G_{m,2}(\cdot,\cdot)$ is the local Green's function (See [KSS, S4]) on $U_m^2(x)$.

We will show that it satisfies both the Dirichlet and Neumann boundary conditions at the boundary of $U_m^2(x)$, and thus could be extended to the whole $K$ by zero extension. Then by Gauss-Green's formula, we have
\begin{equation*}
\tilde{\Delta}_m^2f(x)=-\mathcal{E}(\phi_{m,x}^{(2)},\Delta f).
\end{equation*}

More generally, assume we already have constructed the $(n-1)$-level tent function  $\phi_{m,x}^{(n-1)}$ for vertices $x$ in $V_m^{n-1}$, with the Dirichlet boundary condition at $\partial U_m^{n-1}(x)$, satisfying $\tilde{\Delta}^{n-1}_m f(x)=-\mathcal{E}(\phi_{m,x}^{(n-1)},\Delta^{n-2}f)$. 
We will first define
\begin{equation}
\tilde{\phi}_{m,x}^{(n)}=\sum_{y\sim_m x}(\int \psi_x^m d\mu)^{-1}c_{xy}(\phi_{m,y}^{(n-1)}-\phi_{m,x}^{(n-1)}), 
\end{equation}
which obviously satisfies that
\begin{equation}
\tilde{\Delta}^n_m f(x)=-\mathcal{E}(\tilde{\phi}_{m,x}^{(n)},\Delta^{n-2}f),
\end{equation}
then define
\begin{equation}
\phi_{m,x}^{(n)}(\cdot)=-\int_K G_{m,n}(\cdot,z)\tilde{\phi}_{m,x}^{(n)}(z)d\mu(z),
\end{equation}
where $G_{m,n}(\cdot,\cdot)$ is the local Green's function on $U_m^n(x)$.

We will show that $\phi_{m,x}^{(n)}$ satisfies both the Dirichlet and Neumann boundary conditions at $\partial U_m^{n}(x)$, and then extend it to the whole $K$ by zero extension. Then
\begin{equation}
\tilde{\Delta}^{n}_m f(x)=-\mathcal{E}(\phi_{m,x}^{(n)},\Delta^{n-1}f).
\end{equation}
Thus for $n\geq 2$, $\phi_{m,x}^{(n)}$ defined using the above recipe  will satisfy both the Dirichlet and the Neumann boundary conditions at the boundary $\partial U_m^n(x)$.  The following lemma remains to be proved.

\textbf{Lemma 4.5.} \textit{Let $n\geq 2$. Suppose we have defined $\phi_{m,\cdot}^{(n-1)}$ for vertices in $V_m^{n-1}$, satisfying $\tilde{\Delta}_m^{n-1}f(\cdot)=-\mathcal{E}(\phi_{m,\cdot}^{(n-1)},\Delta^{n-2}f)$, with the Dirichlet boundary condition holding at $\partial U_m^{n-1}(\cdot)$. Then  the function $\phi_{m,x}^{(n)}$ defined by (4.5) and (4.7) satisfies both the Dirichlet and Neumann boundary conditions at the boundary of $U_m^n(x)$. Moreover, the equality (4.8) holds.}

\textit{Proof.}
Let $h$ be a multiharmonic harmonic function in $\mathcal{H}_{n-1}(U_m^n(x))$. By the definition of $\tilde{\phi}_{m,x}^{(n)}$, it is easy to check that
\begin{equation*}
\begin{aligned}
\mathcal{E}(\tilde{\phi}_{m,x}^{(n)}, \Delta^{n-2}h)&=\sum_{y\sim_m x}(\int\psi_x^md\mu)^{-1}c_{xy}\left(\mathcal{E}(\phi_{m,y}^{(n-1)},\Delta^{n-2}h)-\mathcal{E}(\phi_{m,x}^{(n-1)},\Delta^{n-2}h)\right)\\
&=-\sum_{y\sim_mx}(\int\psi_x^md\mu)^{-1}c_{xy}\left(\tilde{\Delta}_m^{n-1}h(y)-\tilde{\Delta}_m^{n-1}h(x)\right)\\
&=-\sum_{y\sim_mx}(\int\psi_x^md\mu)^{-1}c_{xy}\left(\Delta^{n-1}h(y)-\Delta^{n-1}h(x)\right)\\&=0,
\end{aligned}
\end{equation*}
where the second equality comes from the assumption of $\phi_{m,\cdot}^{(n-1)}$,  the third equality is result of Lemma 4.4, and the forth equality follows from the fact that $\Delta^{n-1}h\in \mathcal{H}_0(U_m(x))$. 

On the other hand, by using the Dirichlet boundary condition of $\phi_{m,\cdot}^{(n-1)}$ and $\phi_{m,x}^{(n)}$ at the boundary of $U_m^n(x)$, and repeatedly using the Gauss-Green's formula, we have
\begin{equation*}
\begin{aligned}
\mathcal{E}(\tilde{\phi}_{m,x}^{(n)}, \Delta^{n-2}h)&=-\int_{U_m^n(x)}\tilde{\phi}_{m,x}^{(n)}\Delta^{n-1}hd\mu+\sum_{z\in\partial U_m^n(x)}\tilde{\phi}_{m,x}^{(n)}(z)\partial_n\Delta^{n-1}h(z)\\
&=\int_{U_m^n(x)}\phi_{m,x}^{(n)}\Delta^nhd\mu-\int_{U_m^n(x)}\Delta\phi_{m,x}^{(n)}\Delta^{n-1}hd\mu\\
&=\sum_{z\in\partial U_m^n(x)}\phi_{m,x}^{(n)}(z)\partial_n\Delta^{n-1}h(z)-\sum_{z\in\partial U_m^n(x)}\partial_n\phi_{m,x}^{(n)}(z)\Delta^{n-1}h(z)\\
&=-\sum_{z\in\partial U_m^n(x)}\partial_n\phi_{m,x}^{(n)}(z)\Delta^{n-1}h(z).
\end{aligned}
\end{equation*}

Thus we have proved that 
$\sum_{z\in\partial U_m^n(x)}\partial_n\phi_{m,x}^{(n)}(z)\Delta^{n-1}h(z)=0$ holds for any $h\in\mathcal{H}_{n-1}(U_m^n(x))$, which yields that
$$\sum_{z\in\partial U_m^n(x)}\partial_n\phi_{m,x}^{(n)}(z)h(z)=0$$ holds for any $h\in\mathcal{H}_{0}(U_m^n(x)).$ By the arbitrariness of $h$, we have proved that $\phi_{m,x}^{(n)}$ satisfies both the Dirichlet and Neumann boundary conditions at $\partial U_m^n(x)$.

Now for general function $f\in dom (\Delta^n)$, by (4.6), using the boundary conditions of $\phi_{m,x}^{(n)}$ at $\partial U_m^n(x)$, and by using the Gauss-Green's formula, we finally prove that
$$\tilde{\Delta}_m^n f(x)=-\mathcal{E}(\tilde{\phi}_{m,x}^{(n)},\Delta^{n-2}f)=-\mathcal{E}(\Delta\phi_{m,x}^{(n)},\Delta^{n-2}f)=-\mathcal{E}(\phi_{m,x}^{(n)},\Delta^{n-1}f).\quad \Box$$

\textbf{Lemma 4.6.} \textit{ For any $m$, any $x\in V_m^n$, it holds that $\int \phi_{m,x}^{(n)}d\mu=1.$ Furthermore, for any same type sets $U_m^n(x)$ and $U_{m'}^n(y)$, we have
$
\|\phi_{m,x}^{(n)}\|_1=\|\phi_{{m}',y}^{(n)}\|_1
$ if the $m$-level conductances on $U_m^n(x)$ are proportional to those on $U_{m'}^n(y)$.}

\textit{Proof}. Applying the Gauss-Green's formula to (4.8), using the Dirichlet boundary condition of $\phi_{m,x}^{(n)}$ at $\partial U_m^n(x)$, we have
\begin{equation}
\tilde{\Delta}_m^n f(x)=-\mathcal{E}(\phi_{m,x}^{(n)},\Delta^{n-1}f)=\int \phi_{m,x}^{(n)}\Delta^nf d\mu
\end{equation}
 for any $f\in dom(\Delta^n)$.
Choosing a multiharmonic  function $h\in \mathcal{H}_n (U_m^n(x))$ with $\Delta^nh=1$ and taking it into the above equality, using Lemma 4.4,  we get
\begin{equation*}
\int \phi_{m,x}^{(n)}d\mu=1.
\end{equation*}

Let $U_{m}^n(x)$ and $U_{m'}^n(y)$ be in same type. It means there is a mapping $F$ which is a combination of rotations, reflections and scalings, satisfying $FU_m^n(x)=U_{m'}^n(y)$. It is easy to find that
$$\phi_{m',y}^{(n)}=\frac{\mu(U_m^n(x))}{\mu(U_{m'}^n(y))}\phi_{m,x}^{(n)}\circ F^{-1},$$
by scaling. Hence $
\|\phi_{m,x}^{(n)}\|_1=\|\phi_{{m}',y}^{(n)}\|_1.
$ $\Box$

Since there are only finite types of $U_m^n(x)$, and for each type, there are only finite subtypes with  proportional conductances, we have

\textbf{Corollary 4.7.} \textit{Let $n\geq 2$ be fixed. for any $m$, any $x\in V_m^n$,  $\|\phi_{m,x}^{(n)}\|_1$ is uniformly bounded.}

\textit{Proof of Theorem 4.2(a)}. Applying  (4.9), Lemma 4.6 and Corollary 4.7, we have
\begin{equation*}
\begin{aligned}
|\tilde{\Delta}_m^nf(x)-\Delta^n f(x)|&=|\int \phi_{m,x}^{(n)}(z)(\Delta^nf(z)-\Delta^nf(x)) d\mu(z)|\\
&\leq \|\phi_{m,x}^{(n)}\|_1\omega_{\Delta^nf}(U_m^n(x))\\
&\leq C\omega_{\Delta^nf}(U_m^n(x))
\end{aligned}
\end{equation*}
for some constant $C>0$, where $\omega_{\Delta^nf}(U_m^n(x))$ is the oscillation of $\Delta^nf$ in $U_m^n(x)$. Since $\Delta^nf$ is continuous on $K$, $\omega_{\Delta^nf}(U_m^n(x))$ will go to zero uniformly as $m$ goes to infinity. Thus we have (4.1) holds uniformly.  $\Box$

\textbf{Remark.} The proof provides that the ratio of the convergence in $(4.1)$ depends only on the \textit{modulus of continuity} of $\Delta^nf$.

\subsection{Proof of Theorem 4.2(b)} In this subsection we will give the proof of the second part of Theorem 4.2.

For a simple set $A$ in $K$, we use $S(\mathcal{H}_0,V_m,A)$ to denote the space of \textit{harmonic splines}, which are harmonic on each $m$-level cell in $A$. For those harmonic splines vanishing at  the boundary of $A$, we denote the collection of them by $S_0(\mathcal{H}_0,V_m,A)$. For any function $u\in l(V_m\cap (A\setminus\partial A))$, there is a unique solution $\psi\in S_0(\mathcal{H}_0, V_m, A)$ satisfying 
\begin{equation*}
\tilde{\Delta}_m \psi(x)=u(x),   \quad \forall x\in V_m\cap (A\setminus\partial A). 
\end{equation*}
In fact, for $\psi_1\neq \psi_2\in S_0(\mathcal{H}_0,V_m,A)$, we have $\tilde{\Delta}_m\psi_1\neq \tilde{\Delta}_m\psi_2$, so that $\tilde{\Delta}_m$ is an injection. Comparing the dimension, one can find that $\tilde{\Delta}_m$ is reversible. For convenience, we define  $G_{m,A}$ the inverse operator of $-\tilde{\Delta}_m$. It means that for any  $u\in l(V_m\cap (A\setminus\partial A))$,  $G_{m,A}u\in S_0(\mathcal{H}_0,V_m,A)$ and $$-\tilde{\Delta}_m G_{m,A}u=u.$$

We denote $G_A$ the local Green's operator on $A$, i.e., for any continuous function $u$ on $A$,  $G_Au\in dom_0(\Delta,A)$, and satisfies 
\begin{equation*}
-\Delta G_Au=u.
\end{equation*}
The following lemma shows that $G_{m,A}$ will go to $G_A$ as $m$ goes to infinity, in some sense.

\textbf{Lemma 4.8.} \textit{For any simple set $A$ in $K$, let $f\in dom_0(\Delta,A)$ and $\psi_m\in S_0(\mathcal{H}_0,V_m,A)$. If
$\tilde{\Delta}_m \psi_m$ converges to $\Delta f$ uniformly as $m$ goes to infinity,
then $\psi_m$ converges to $f$ uniformly as $m$ goes to infinity.}

\textit{Proof}. First we will show $\{\psi_m\}$ are equicontinuous and uniformly bounded. Observe that for any function $g\in dom_0(\mathcal{E},A)$,
\begin{equation*}
\begin{aligned}
\mathcal{E}^A(\psi_m,g)&=-\sum_{x\in V_m\cap (A\setminus \partial A)}g(x)\Delta_m \psi_m(x)\\
&=-\sum_{x\in V_m\cap (A\setminus \partial A)} (\int_K\psi_x^m d\mu) g(x)\tilde{\Delta}_m \psi_m(x)\\
&=-\sum_{\substack{F_wK\subset A\\|w|=m}}\sum_{x\in F_wV_0} \frac{1}{3}\mu_wg(x)\tilde{\Delta}_m \psi_m(x).
\end{aligned}
\end{equation*}
Noticing that $\psi_m$ satisfies the Dirichlet boundary condition at $\partial A$, combining with the estimate $|\psi_m(x)-\psi_m(y)|^2\leq R(x,y)\mathcal{E}^A(\psi_m)$, with $R(x,y)$  the effective resistance metric (See [Ki4, S4]) between $x$ and $y$ on A, we have
\begin{equation*}
\|\psi_m\|_\infty^2\leq c_1\mathcal{E}^A(\psi_m)\leq c_2\|\psi_m\|_\infty \|\tilde{\Delta}_m\psi_m\|_\infty,
\end{equation*}
for some constants $c_1,c_2>0$, which results that
\begin{equation*}
\|\psi_m\|_\infty^2\leq c_1\mathcal{E}^A(\psi_m)\leq c_2^2 \|\tilde{\Delta}_m\psi_m\|^2_\infty.
\end{equation*}

Since we have $\tilde{\Delta}_m \psi_m\rightarrow \Delta f$ uniformly, $\|\tilde{\Delta}_m\psi_m\|_\infty$ is uniformly bounded. So we have $\{\psi_m\}$ are uniformly bounded and equicontinuous. Denote $f-\psi_m=g_m$, then $\{g_m\}$ are also uniformly bounded and equicontinuous. Moreover,
\begin{equation*}
\begin{aligned}
\mathcal{E}^A(g_m)&=-\int_A g_m\Delta fd\mu+\sum_{\substack{F_wK\subset A\\|w|=m}}\sum_{x\in F_wV_0} \frac{1}{3}\mu_wg_m(x)\tilde{\Delta}_m \psi_m(x)\\
&=\sum_{\substack{F_wK\subset A\\|w|=m}}\left(-\int_{F_wK} g_m\Delta fd\mu+\sum_{x\in F_wV_0}\frac{1}{3}\mu_wg_m(x)\tilde{\Delta}_m \psi_m(x)\right)\\
&= \sum_{\substack{F_wK\subset A\\|w|=m}}\mu_w\left(-\int_K g_m\circ F_w(\Delta f)\circ F_wd\mu+\sum_{x\in F_wV_0}\frac{1}{3}g_m(x)\tilde{\Delta}_m \psi_m(x)\right).
\end{aligned}
\end{equation*}
All terms in the sum would converge to $0$ uniformly as $m$ goes to $\infty$. In fact,
\begin{equation*}
\begin{aligned}
&\left|-\int_K g_m\circ F_w(\Delta f)\circ F_wd\mu+\sum_{x\in F_wV_0}\frac{1}{3}g_m(x)\tilde{\Delta}_m \psi_m(x)\right|\\
\leq&\sum_{x\in F_wV_0}\frac{1}{3}\left|-\int_K g_m\circ F_w(\Delta f)\circ F_wd\mu+g_m(x)\tilde{\Delta}_m \psi_m(x)\right|\\
\leq&\sum_{x\in F_wV_0} \frac{1}{3}\sup_{y\in F_wK}|g_m(y)\Delta f(y)-g_m(x)\tilde{\Delta}_m \psi_m(x)|.
\end{aligned}
\end{equation*}
Using the equicontinuous and uniformly boundedness of $\{g_m\}$, we then have $\lim_{m\to\infty}\mathcal{E}^A(f-\psi_m)=0$. Together with the fact that $f-\psi_m$ satisfies the Dirichlet boundary condition at $\partial A$, it yields the result of Lemma 4.8. $\Box$

\textbf{Remark.} \textit{We may restate Lemma 4.8 as follows.  Suppose $\phi_m\in S(\mathcal{H}_0,V_m,A)$ converges uniformly to a continuous function $u$, then 
\begin{equation}
\lim_{m\to\infty}G_{m,A}\phi_m=G_Au
\end{equation}
holds uniformly.}

\textit{Proof of Theorem 4.2(b)}. Assume we have $\lim_{m\to \infty}\tilde{\Delta}_m^n f(x)=u(x)$ uniformly on $V_*\setminus V_0$. Then by repeatedly using Lemma 4.8, on any $A$ not intersecting the boundary $V_0$, we have
\begin{equation*}
\lim_{m\to \infty}(-G_{m,A})^n\tilde{\Delta}_m^nf=(-G_A)^n u
\end{equation*}
converges uniformly. 
So we have $f-(-G_{m,A})^n\tilde{\Delta}_m^nf$ converges uniformly to the function $f-(-G_A)^n u$. 

Now we prove $f-(-G_A)^n u\in\mathcal{H}_{n-1}(A)$.

Recall that in Lemma 4.3, we have shown that for any $(n+1)$-harmonic function $h\in \mathcal{H}_n(A)$, $\tilde{\Delta}_m h$ must be some $n$-harmonic function, see (4.4). It is not hard to verify that the right side of (4.4) could go through all $n$-harmonic functions. Thus we have an inverse conclusion that for any $n$-harmonic function $h'$ on $A$, there is a $(n+1)$-harmonic function $h\in\mathcal{H}_n(A)$ such that $\tilde{\Delta}_mh=h'$ on $V_m\cap(A\setminus\partial A)$. 

Now we apply the above discussion in our proof. First we have the following equality,
\begin{equation*}
\tilde{\Delta}_m (\tilde{\Delta}_m^{n-1}f+G_{m,A}\tilde{\Delta}_m^nf)=0.
\end{equation*}
Thus $\tilde{\Delta}_m^{n-1}f+G_{m,A}\tilde{\Delta}_m^nf$ equals to some   harmonic function on $A$. Since for each $1< i\leq n$,
\begin{equation*}
\tilde{\Delta}_m (\tilde{\Delta}_m^{n-i}f-(-G_{m,A})^{i}\tilde{\Delta}_m^nf)=\tilde{\Delta}_m^{n-i+1}f-(-G_{m,A})^{i-1}\tilde{\Delta}_m^nf,
\end{equation*}
by repeatlly using the above discussion, we have that  $f-(-G_{m,A})^n\tilde{\Delta}_m^nf$ equals to some $n$-harmonic function on $A$. Noticing that the space of $n$-harmonic functions on $A$ is of finite dimension, the uniform limit $f-(-G_A)^n u$ of  $f-(-G_{m,A})^n\tilde{\Delta}_m^nf$ is of course a $n$-harmonic function. 

Thus we have $f=(-G_A)^nu+(f-(-G_A)^nu)\in dom(\Delta^n, A)$, and obviously $\Delta^nf=u$ on $A$. By the arbitrariness of $A$, we finally have proved $\Delta^nf=u$ on $K\setminus V_0.  \quad \Box$

\section{Pointwise formula of $\Delta_\mu^n$ on general p.c.f fractals}
We have no idea on how to extend the previous results to other p.c.f. fractals. However, we still have some pointwise calculations of the higher order Laplacian in general.\\

The following is an extension of the \textit{mean value property} of harmonic functions.

\textbf{Lemma 5.1.} \textit{Let $l\in \mathbb{N}$,  $\{y_j\}_{j=1}^l\subset V_*$, and  $\{a_j\}_{j=1}^{l}\subset\mathbb{R}$ with
$\sum_{j=1}^{l} a_j=0.$
Then there exists a function $\phi\in dom\mathcal{E}$ such that  
\begin{equation}
\sum_{j=1}^{l}a_jf(y_j)=-\mathcal{E}(f,\phi)
\end{equation}
holds for any $f\in dom\mathcal{E}$.
Furthermore, if we additionally assume $\sum_{j=1}^{l}a_jh(y_j)=0$ holds for any $h\in\mathcal{H}_0$, then there is a unique such $\phi$ satisfying the  $0$ boundary condition $\phi|_{V_0}=0$.}

\textit{Proof}. Assume $\{y_j\}_{j=1}^{l}\subset V_m$ for some $m$, we now find the function $\phi$ in $S(\mathcal{H}_0,V_m)$, the space of continuous functions which are harmonic on each $m$-level cell of $K$.

In fact, For any $\psi\in S(\mathcal{H}_0,V_m)$, we have
\begin{equation*}
-\mathcal{E}(f,\psi)=\sum_{y\in V_m\setminus V_0}\Delta_m \psi(y)f(y)-\sum_{y\in V_0}\partial_n \psi(y)f(y), \forall f\in dom\mathcal{E}.
\end{equation*}
It is easy to check that the map 
$\psi\rightarrow \{{\partial_n \psi|_{V_0},\Delta_m \psi|_{V_m\setminus V_0}}\}$ is injective from the space $S(\mathcal{H}_0,V_m)$ modulo constants to $\mathbb{R}^{\#V_m}$.
Additionally, by the Gauss-Green's formula, it always holds 
\begin{equation*}
-\sum_{y\in V_m\setminus V_0}\Delta_m \psi(y)+\sum_{y\in V_0}\partial_n \psi(y)=0.
\end{equation*}
 Thus the map is a bijection from $S(\mathcal{H}_0,V_m)$ modulo constants to a $(\#V_m-1)$-dimensional subspace of $\mathbb{R}^{\#V_m}$ by a counting dimension argument. Thus, we have proved the existence of function $\phi$ satisfying (5.1). 

If we additionally assume $\sum_{j=1}^la_jh(y_j)=0$, then we have 
$\mathcal{E}(h,\phi)=0$
for any $h\in \mathcal{H}_0$. Thus
\begin{equation*}
\mathcal{E}_0(h,\phi)=0, \forall h\in \mathcal{H}_0,
\end{equation*}
which yields that $\phi|_{V_0}=C$, since $\mathcal{E}_0(\cdot,\cdot)$ is an inner product on $\mathbb{R}^{\# V_0}$ modulo constants. So there exists a unique $\phi$ satisfying (5.1)
with the Dirichlet boundary condition. $\Box$

\textbf{Theorem 5.2. (Calculation of $\Delta_\mu$)} \textit{Let $l\in\mathbb{N}$,  $\{y_j\}_{j=1}^l\subset V_*$ and $\{{a_j}\}_{j=1}^l\subset \mathbb{R}$. Assume the following condition holds}
\begin{equation*}
\begin{cases}
\sum_{j=1}^{l} a_jh(y_j)=0, \forall  h\in \mathcal{H}_0,\\
\sum_{j=1}^{l} a_jh'(y_j)=A, \forall  h'\emph{ with }\Delta_\mu h'=1,
\end{cases}
\end{equation*}
\textit{for some constant $A\neq 0$. Then 
\begin{equation*}
\Delta_\mu f(x)=\lim_{m\to\infty} A^{-1}(r_{[w]_m}\mu_{[w]_m})^{-1}\sum_{j=1}^{l}a_jf(F_{[w]_m}y_j)
\end{equation*}
uniformly on $K$ for any function $f\in dom\Delta_\mu$. Here $w$ is an infinite word corresponding to $x$ and $ F_{[w]_m}K$ denotes the according $m$-cell containing $x$ for each $m\geq 0$. }

\textit{Proof}. According to Lemma 5.1, there exists a piecewise harmonic spline $\phi$ satisfying (5.1) with the Dirichlet boundary condition. Using the Gauss-Green's formula, for each $f\in dom\Delta_\mu$,
\begin{equation}
\sum_{j=1}^{l}a_jf(y_j)=-\mathcal{E}(f,\phi)=\int_K \phi\Delta_\mu fd\mu.
\end{equation}
Consider a $h'$ with $\Delta_\mu h'=1$, then $\sum_{j=1}^{l}a_jh'(y_j)=A$, and so $\int_K \phi d\mu=A$.
By using a scaling of the identity (5.2), for each $m\geq 0$, we get
\begin{equation*}
\sum_{j=1}^{l}a_jf\circ F_{[w]_m}(y_j)=\int_K \phi\Delta_\mu (f\circ F_{[w]_m})d\mu=r_{[w]_m}\mu_{[w]_m}\int_K \phi(\Delta_\mu f)\circ F_{[w]_m}d\mu.
\end{equation*}
Taking the limit as $m\rightarrow\infty$, we have proved the theorem. $\Box$
\\

Now we turn to the higher order case.

\textbf{Lemma 5.3.}\textit{ Let $l\in \mathbb{N}$,  $\{{y_j}\}_{j=1}^l\subset V_*$ and $\{a_j\}_{j=1}^l\in\mathbb{R}$ with
$
\sum_{j=0}^{l} a_jh(y_j)=0
$
holding for any  $h\in \mathcal{H}_{n-1}$, then there exists a function $\phi_n$ with the Dirichlet boundary condition, such that
\begin{equation}
\sum_{j=0}^{l} a_jf(y_j)=-\mathcal{E}(\Delta_\mu^{n-1}f,\phi_n)
\end{equation}
holds for any  $f\in dom(\Delta_\mu^{n-1})$.}

\textit{Proof}. For $n=1$, It is just what Lemma 5.1 says, so we get the initial function $\phi_1$. 

Now, we assert that we could choose
\begin{equation*}
\phi_n=(-1)^{n-1}\int G(x,z_1)\cdots G(z_{n-2},z_{n-1})\phi_1(z_{n-1})d\mu(z_{n-1})\cdots d\mu(z_1).
\end{equation*} where $G(\cdot,\cdot)$ is the Green's function solving the Dirichlet problem of the Poisson equation on $K$.

In fact, assume Lemma 5.3 holds for $n-1$ case, then $\forall f\in dom(\Delta_\mu^{n-1})$,
\begin{equation*}
\sum_{j=0}^{l} a_jf(y_j)=-\mathcal{E}(\Delta_\mu^{n-2}f,\phi_{n-1}).
\end{equation*}
Using the Gauss-Green's formula, by the assumption of $\{a_j\}_{j=1}^l$, and using the Dirichlet boundary condition of $\phi_{n-1}$, $\phi_n$, we have
\begin{equation*}
\begin{aligned}
0&=\sum_{j=0}^{l} a_jh(y_j)=-\mathcal{E}(\Delta_\mu^{n-2}h,\phi_{n-1})\\
 &=\int_K \phi_{n-1}\Delta_\mu^{n-1}hd\mu-\sum_{z\in V_0}\phi_{n-1}(z)\partial_n\Delta_\mu^{n-2}h(z)\\
 &=\sum_{z\in V_0}\left(\Delta_\mu^{n-1}h(z)\partial_n \phi_n(z)-\phi_n(z)\partial_n \Delta_\mu^{n-1}h(z)\right)\\
 &=\sum_{z\in V_0}\Delta_\mu^{n-1}h(z)\partial_n \phi_n(z)
\end{aligned}
\end{equation*} holds for any $h\in\mathcal{H}_{n-1}$.
Noticing that $\Delta_\mu^{n-1}h$ goes through the whole space $\mathcal{H}_0$, we could get
\begin{equation*}
\partial_n\phi_n|_{V_0}=0.
\end{equation*}
Thus by using the Gauss-Green's formula again, we have
\begin{equation*}
\sum_{j=0}^{l} a_jf(y_j)=-\mathcal{E}(\Delta_\mu^{n-2}f,\phi_{n-1})
=\int_K \phi_{n-1}\Delta_\mu^{n-1}fd\mu
=-\mathcal{E}(\Delta_\mu^{n-1}f,\phi_n). \quad\Box
\end{equation*}
\\
\textbf{Theorem 5.4. (Calculation of $\Delta_\mu^n$)} \textit{Let $l\in\mathbb{N}$, $\{y_j\}_{j=1}^l\subset V_*$, and $\{a_j\}_{j=1}^{l}\subset\mathbb{R}$. Assume the following condition holds
\begin{equation*}
\begin{cases}
\sum_{j=1}^{l} a_jh(y_j)=0, \forall  h\in \mathcal{H}_{n-1},\\
\sum_{j=1}^{l} a_jh'(y_j)=A, \forall h' \emph{ with }  \Delta^n h'=1,
\end{cases}
\end{equation*}
for some constant $A\neq 0$. Then
\begin{equation*}
\Delta_\mu^n f(x)=\lim_{m\to\infty} A^{-1}(r_{[w]_m}\mu_{[w]_m})^{-n}\sum_{i=1}^{l}a_jf(F_{[w]_m}y_j),
\end{equation*}
uniformly on $K$ for any function $f\in dom(\Delta_\mu^n)$. }

\textit{Proof.} By Lemma 5.3, there exists a function $\phi_n$ satisfying (5.3) with the Dirichlet boundary condition. Combining it with the Gauss-Green's formula, we get
\begin{equation}
\sum_{j=1}^{l}a_jf(y_j)=-\mathcal{E}(\Delta_\mu^{n-1}f,\phi_n)=\int_K \phi_n\Delta_\mu^n fd\mu.
\end{equation}
Let $h'$ be any function with $\Delta_\mu^n h'=1$. Then $\sum_{j=1}^{l}a_jh'(y_j)=A$ and thus $\int_K \phi_n d\mu=A$.
Scaling the identity  (5.4), for each $m\geq 0$, we have
\begin{equation*}
\sum_{j=1}^{l}a_jf\circ F_{[w]_m}(y_j)=\int_K \phi_n\Delta_\mu^n (f\circ F_{[w]_m})d\mu=(r_{[w]_m}\mu_{[w]_m})^n\int_K \phi_n(\Delta_\mu^n f)\circ F_{[w]_m}d\mu.
\end{equation*}
Taking  the limit as $m\rightarrow\infty$, we have proved the theorem. $\Box$

\textbf{Remark.} From the proof of the theorem, it is easy to find that the ratio of the uniform convergence depends only on the modulus of continuity of $\Delta_\mu^n f$ as stated in Section 4.

\section{Extension to the $D4$ symmetric p.c.f. fractals}
In this section, we manage to extend the previous results to the $D4$ symmetric p.c.f. fractals. (Actually it could be extended to those fully symmetric  p.c.f. fractals with regular harmonic structure.)  In this case, $N_0=4$, $r=r_l$, $\mu=\mu_l$ for $l=0,1,2,3$, and we denote $\rho=r\mu$. Similar to the previous discussion, throughout this section, we assume
\begin{equation*}
r_l\mu_l=\rho, 0\leq l< N.
\end{equation*}

It is well known that the \textit{Vicsek set} $\mathcal{VS}$ and the \textit{tetrahedral Sierpinski gasket} $\mathcal{SG}^4$ are two typical $D4$ symmetric examples. Consider a square with corners $\{q_0,q_1,q_2,q_3\}$ and center $q_4$. Let $F_l$ be contractive mappings with ratio $1/3$ and fixed points $q_l$. The generated invariant set is called the Vicsek set. The fractal and the second step graph is as shown in Fig. 6.1. It could be naturally extended to the $n$-branch Vicsek set $\mathcal{VS}_n$. We omit it here. Please see [Z] and [CSW] for the spectral analysis of the Laplacians on this family of sets. $\mathcal{SG}^4$ is the fractal satisfying $\mathcal{SG}^4=\bigcup_{l=0}^{3}F_l\mathcal{SG}^4$ where $F_lx=\frac{1}{2}(x+q_l)$ and $\{q_0,q_1,q_2,q_3\}$ are the four vertices of a tetrahedron.  See Fig. 6.2.
\begin{figure}[h]
\begin{center}
\includegraphics[width=4.2cm]{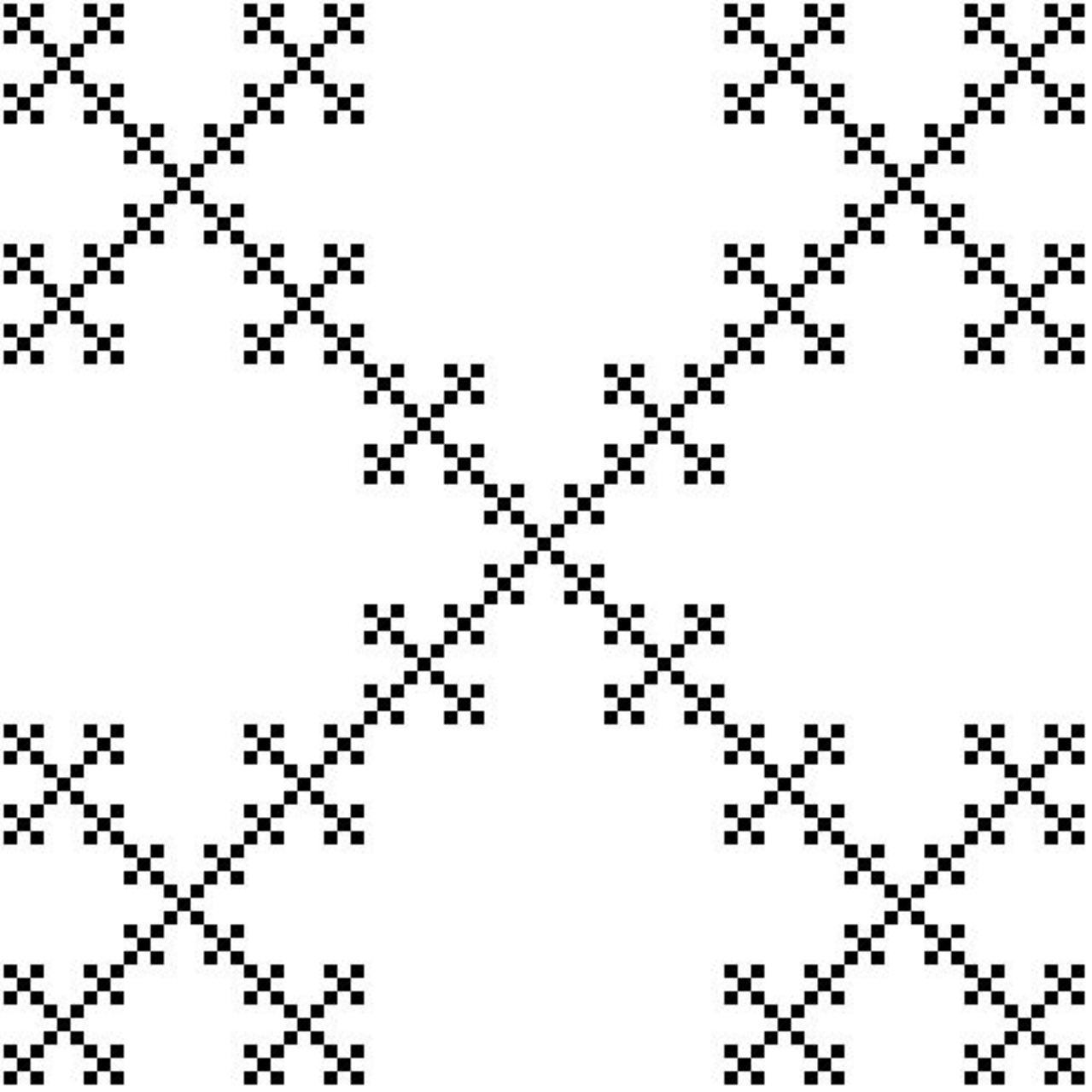}\hspace{2cm}
\includegraphics[width=4.2cm]{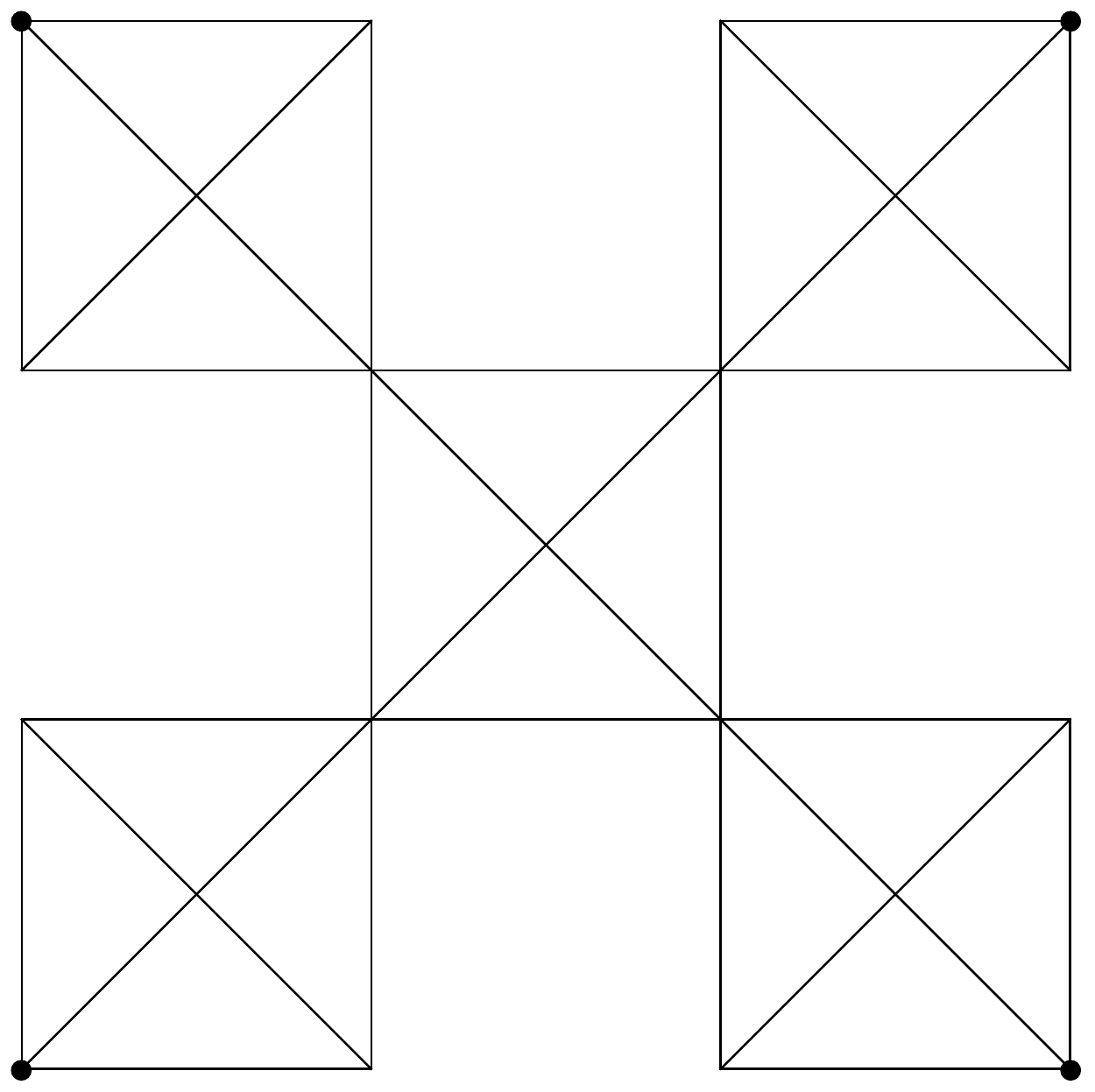}
\caption{The Vicsek set $\mathcal{VS}$.}
\end{center}
\end{figure}

 \begin{figure}[htbp]   
     \centering
     \includegraphics[width=0.4\textwidth]{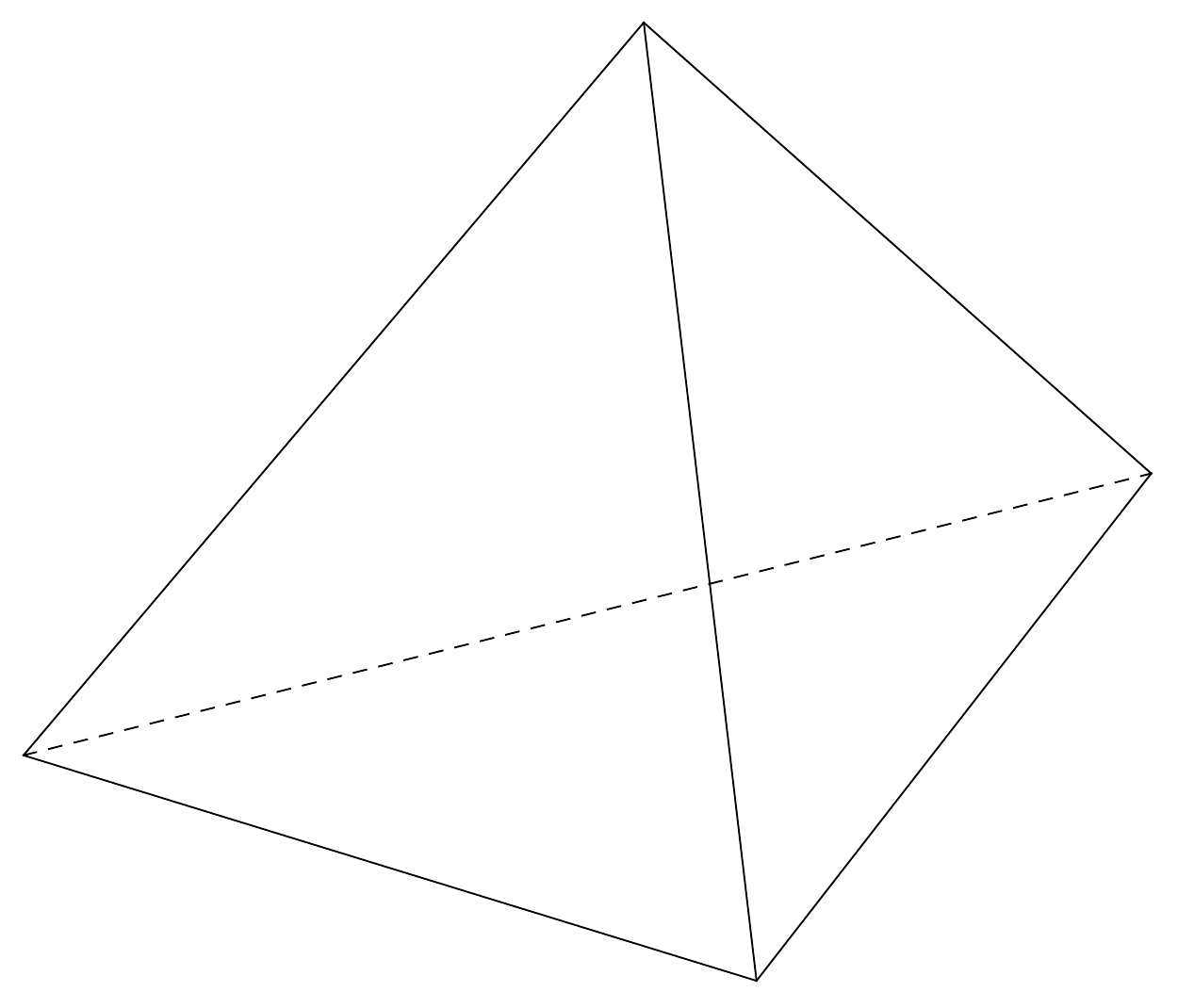}
     \includegraphics[width=0.4\textwidth]{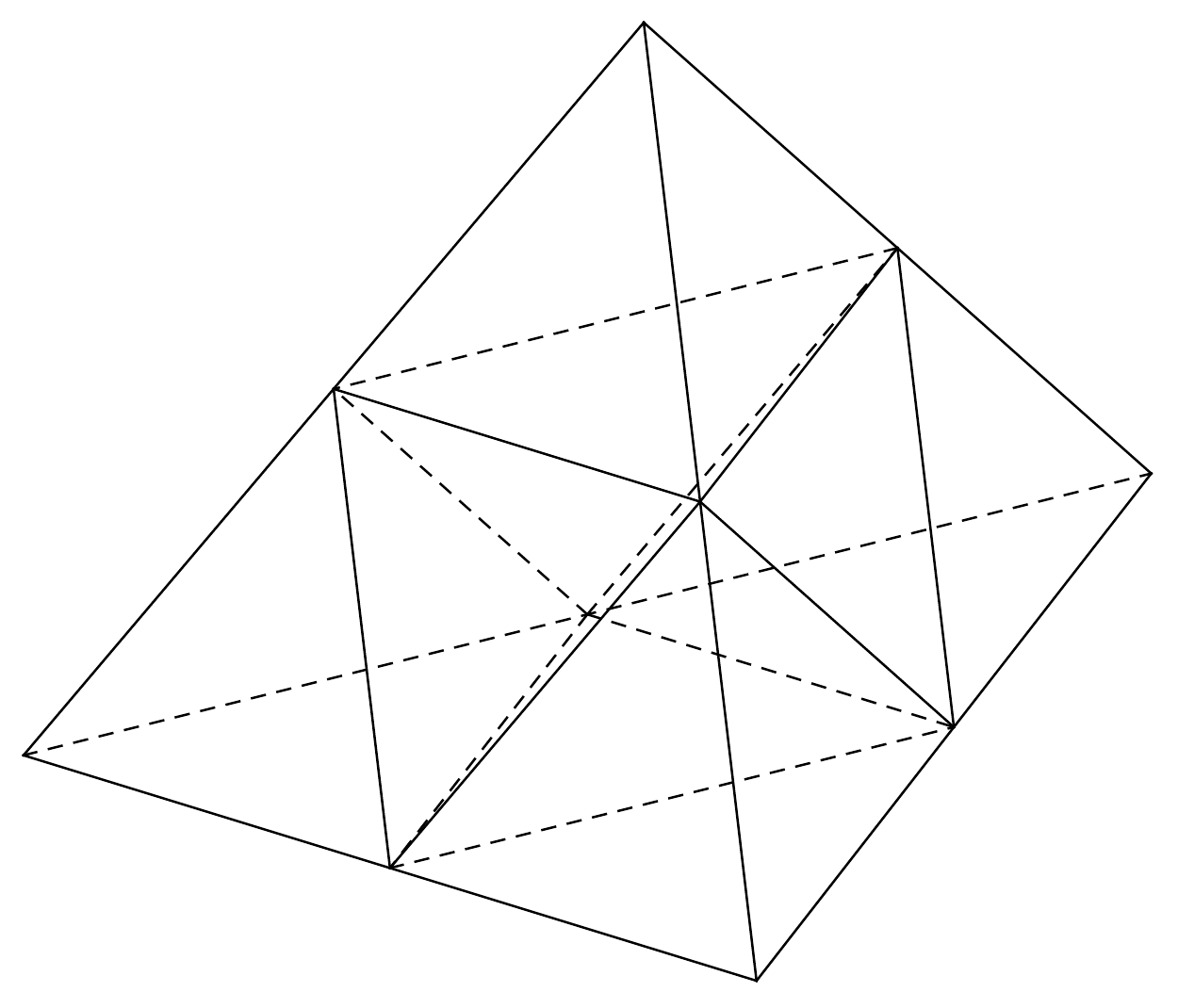}
     \caption{First and second step graphs of $\mathcal{SG}^4$.}
\end{figure}

Firstly, the results in Section 4 could be extended to the $D4$ symmetric case. Noticing that the harmonic extension matrices may be degenerate, for example, in the case of $\mathcal{VS}$, to avoid the inconvenience of the definition of monomials, we will use the ``easy'' basis $\{f_{jl}\}$ instead, which was introduced in [SU]. The argument below provides another proof of Lemma 4.3, in the $D4$ version,  with suitable modifications. Here, the notations $U_m(x)$, $U_m^n(x)$, the boundary of simple sets $\partial A$ and the renormalized graph Laplacian $\tilde{\Delta}_m$ are same as the previous ones.\\

\textbf{Definition 6.1.} \textit{For $j\geq 0$, $0\leq l< 4$, define $f_{jl}\in \mathcal{H}_j$} satisfying
\[
\Delta^i f_{jl}(q_{l'})=\delta_{ij}\delta_{ll'}\quad\text{for $0\leq l<4$ and $i\geq 0$.}
\]

It is easy to check that $\{f_{jl}\}_{j\leq n, l<4}$ form an  ``easy''  basis of $\mathcal{H}_n$ with the dimension $4(n+1)$.
By the symmetry, we denote
\begin{equation*}
\partial_n f_{jl}(q_{l})=a_j, \mbox{ and }\partial_n f_{jl}(q_{l'})=b_j\mbox{ for } l'\neq l.
\end{equation*}
Since $\Delta^i f_{jl}=f_{(j-i)l}$, for $i\leq j$, it is obvious that 
\begin{equation*}
\partial_n\Delta^i f_{jl}(q_{l})=\begin{cases}
a_{j-i},&\text{ if }i\leq j,\\0,&\text{ if }i>j,\end{cases}
\mbox{ and }\quad
\partial_n\Delta^i f_{jl}(q_{l'})=\begin{cases}
b_{j-i},&\text{ if }i\leq j,\\0,&\text{ if }i>j.\end{cases}
\end{equation*}

More generally, for any $h\in \mathcal{H}_n$ and any $i\leq n$, we have
\begin{equation*}
\Delta^i h=\sum_{j=0}^{n-i}\sum_{l'=0}^{3} \Delta^{i+j}h(q_{l'})f_{jl'},
\end{equation*} 
which results that
\begin{equation*}
\partial_n\Delta^i h(q_l)=\sum_{j=0}^{n-i} \Delta^{i+j} h(q_l)a_j+\sum_{j=0}^{n-i}\sum_{l'\neq l} \Delta^{i+j} h(q_{l'})b_j.
\end{equation*}

\textbf{Lemma 6.2.(matching condition)} \textit{Let $x$ be a vertex in $V_*\setminus V_0$, $h$ be any $(n+1)$-harmonic function in $\mathcal{H}_n(U_m(x))$. Then for any $i\leq n$, the following identity holds,
\begin{equation}
\sum_{j=0}^{n-i}\sum_{y\sim_m x}c_{xy}\rho^{mj}(\frac{1}{3}\Delta^{i+j}h(x)a_j+\Delta^{i+j}h(y)b_j)=0.
\end{equation}
} 

\textit{Proof.} The matching condition at $x$ is 
\begin{equation*}
\sum_{w\in W(x)}\partial^w_n\Delta^i h(x)=0.
\end{equation*}
For each summand above, by scaling, we have
\[\begin{aligned}
\partial^w_n\Delta^i h(x)=&r_w^{-1}r^{-(m-|w|)}\rho^{-mi}\partial_n \Delta^i(h\circ F_x^{m-|w|}\circ F_w)(q_l)\\
=&r_w^{-1}r^{-(m-|w|)}\rho^{-mi}\left(\sum_{j=0}^{n-i} \Delta^{i+j} (h\circ F_x^{m-|w|}\circ F_w)(q_l)a_j\right.\\
&+\left.\sum_{j=0}^{n-i}\sum_{l'\neq l} \Delta^{i+j} (h\circ F_x^{m-|w|}\circ F_w)(q_{l'})b_j\right)\\
=&r_w^{-1}r^{-(m-|w|)}\left(\sum_{j=0}^{n-i}\Delta^{i+j}h(x)\rho^{mj}a_j+\sum_{j=0}^{n-i}\sum_{y\sim_m x,y\in F_wK}\Delta^{i+j} h(y)\rho^{mj}b_j\right),
\end{aligned}\]
where $w\in W(x)$ with $x=F_wq_l$. Thus we can rewrite the matching condition into
(6.1),
 using the fact that $r_w^{-1}r^{-(m-|w|)}=c_{xy}$ and $N_0=4$. $\Box$

\textbf{Lemma 6.3.} \textit{Let $x$ be a vertex in $V_m\setminus V_0$, $h$ be any $(n+1)$-harmonic function in $\mathcal{H}_n(U_m(x))$. Then the following equality holds, 
\begin{equation}
\tilde{\Delta}_m h(x)=12\sum_{j=1}^{n} \rho^{m(j-1)}c_j\Delta^j h(x),
\end{equation}
where $\{c_j\}_{j=0}^{\infty}$ is a sequence of constants satisfying
\begin{equation*}
c_j=\frac{1}{3}a_j+\sum_{s=0}^{j-1}b_{j-s}c_s,
\end{equation*}
with initial value $a_0=3$, $b_0=-1$, $c_0=1$. 
}

\textit{Proof.} In fact, (6.2) follows from a transformation of the matching condition (6.1). 

When $i=n$, (6.1) is just
$\sum_{y\sim_m x}c_{xy}(\Delta^nh(x)-\Delta^nh(y))=0$, and for convenience, we write it as
\[\sum_{y\sim_m x}c_{xy}\Delta^nh(y)=\Delta^nh(x)\sum_{y\sim_m x}c_{xy}.\]
When $i=n-1$, (6.1) becomes
\begin{equation*}
\sum_{y\sim_m x} c_{xy}(\Delta^{n-1}h(x)-\Delta^{n-1}h(y))+\rho^m\sum_{y\sim_m x} c_{xy}(\frac{1}{3}\Delta^nh(x)a_1+\Delta^nh(y)b_1)=0,
\end{equation*} 
which yields that
\begin{equation*}
\begin{aligned}
\sum_{y\sim_m x}c_{xy}\Delta^{n-1}h(y)&=\sum_{y\sim_m x}c_{xy}\Delta^{n-1}h(x)+\sum_{y\sim_m x}c_{xy}\rho^m(\frac{1}{3}a_1+b_1)\Delta^n h(x)\\
&=\sum_{y\sim_m x}c_{xy}(c_0\Delta^{n-1}h(x)+\rho^mc_1\Delta^n h(x)).
\end{aligned}
\end{equation*}
Recursively, for $i<n-1$, (6.1) becomes
\begin{equation*}
\sum_{y\sim_m x} c_{xy}(\Delta^{i}h(x)-\Delta^{i}h(y))+\sum_{j=1}^{n-i}\sum_{y\sim_m x}\rho^{mj}c_{xy}(\frac{1}{3}\Delta^{i+j}h(x)a_j+\Delta^{i+j}h(y)b_j)=0,
\end{equation*}
so that
\[\begin{aligned}
\sum_{y\sim_m x} c_{xy}\Delta^i h(y)=&\sum_{y\sim_m x} c_{xy}\Delta^ih(x)+(\sum_{y\sim_m x}c_{xy})\cdot\\&\sum_{j=1}^{n-i}\rho^{mj}(\frac{1}{3}a_j\Delta^{i+j}h(x)+b_j\sum_{s=0}^{n-i-j}\rho^{ms}c_{s}\Delta^{i+j+s}h(x))\\
=&\sum_{y\sim_m x} c_{xy}\Delta^ih(x)+(\sum_{y\sim_m x}c_{xy})\sum_{j=1}^{n-i}\rho^{mj}(\frac{1}{3}a_j+\sum_{s=0}^{j-1}c_sb_{j-s})\Delta^{i+j}h(x)\\
=&(\sum_{y\sim_m x}c_{xy})\sum_{j=0}^{n-i}\rho^{mj}c_j\Delta^{i+j} h(x).
\end{aligned}\]
Thus, we have
\begin{equation*}
\sum_{y\sim_m x} c_{xy}h(y)=(\sum_{y\sim_m x}c_{xy})\sum_{j=0}^{n} \rho^{mj}c_j\Delta^j h(x).
\end{equation*}
Multiplying the both sides of the above equality with $(\int_K \psi_x^md\mu)^{-1}$, we finally get
\begin{equation*}
\tilde{\Delta}_m h(x)=12\sum_{j=1}^{n} \rho^{m(j-1)}c_j\Delta^j h(x).\quad  \Box
\end{equation*}

Having this lemma, the $D4$ version of Lemma 4.3, we could follow the same steps to prove the pointwise formula of $\Delta^n$ as  shown in Section 4 for $D4$ symmetric fractals. The argument is very similar, so we omit it. Thus Question $2$ still has a positive answer in $D4$ symmetric case.

Secondly, we turn to Question $1$. For simplicity, we use the simplest $D4$ fractal $\mathcal{SG}^4$ as an example to show how to deal with the $D4$ symmetric fractals. Here we should require that all the harmonic extension matrices (differ only by permutations) to be nondegenerate.

It is easy to check that for $D4$ symmetric fractals, the third and fourth eignvalues coincide, which we still denote by $\lambda$. Thus, we would view the related two derivatives as two components of the ``transverse derivative''. We then define all the four types of derivatives as introduced in [S3].

The \textit{normal derivative} of a function $f$ at the boundary point $q_l$ is defined as
$$\partial_n f(q_l)=\lim_{m\rightarrow\infty}r^{-m}(3f(q_l)-f(F_{l}^m q_{l+1})-f(F_{l}^m q_{l+2})-f(F_{l}^m q_{l+3}))$$ (cyclic notation $q_{l+4}=q_l$),
providing the limit exists, while the \textit{transverse derivatives} at $q_l$ are defined as
$$
\begin{cases}
\partial_{T,1} f(q_l)=\lim_{m\rightarrow\infty}\lambda^{-m}(2f(F_{l}^m q_{l+1})-f(F_{l}^m q_{l+2})-f(F_{l}^m q_{l+3})),\\
\partial_{T,2} f(q_l)=\lim_{m\rightarrow\infty}\lambda^{-m}(2f(F_{l}^m q_{l+2})-f(F_{l}^m q_{l+1})-f(F_{l}^m q_{l+3})),
\end{cases}
$$
providing the limits exist. Of course, the definitions of these derivatives could also be localized to all vertices. We omit the details. 

The following is the definition of  \textit{monomials} on $\mathcal{SG}^4$.

\textbf{Definition 6.4.} \textit{Fix a boundary point $q_l$. The monomials $Q_{jk}^{(l)}$ for $1\leq k\leq 4$ and $0\leq j\leq n$  in $\mathcal{H}_n$ are  the multiharmonic functions satisfying
\[
\begin{aligned}
\Delta^i Q_{jk}^{(l)}(q_l)=\delta_{ij}\delta_{k1},\\
\partial_n\Delta^i Q_{jk}^{(l)}(q_l)=\delta_{ij}\delta_{k2},\\
\partial_{T,1}\Delta^i Q_{jk}^{(l)}(q_l)=\delta_{ij}\delta_{k3},\\
\partial_{T,2}\Delta^i Q_{jk}^{(l)}(q_l)=\delta_{ij}\delta_{k4}.
\end{aligned}
\] }
The $\{Q_{j1}^{(l)}\}$ and $\{Q_{j2}^{(l)}\}$ are symmetric under rotations and reflections which fix $q_l$, while $\{Q_{j3}^{(l)}\}$  are skew-symmetric under the reflection $g_{l2}$ and $\{Q_{j4}^{(l)}\}$ are skew-symmetric under the reflection $g_{l1}$, where $g_{li}$ denote the reflection which preserves $q_l$ and $q_{l+i}$ and permutes the other two boundary points. The self-similar identities (3.1)-(3.3) for $D3$ symmetric cases still hold for $\mathcal{SG}^4$ under suitable modification. We keep using the notations
\[
\alpha_j=Q^{(0)}_{j1}(q_1),\beta_j=Q_{j2}^{(0)}(q_1),\gamma_j=Q_{j3}^{(0)}(q_1)=Q_{j4}^{(0)}(q_2).
\]
Analogous to $D3$ symmetric cases, we can define the local monomials. We omit the details. 

Still we hope to decompose the multiharmonic functions as in Theorem 3.5. Denote $R_i$ the rotations in the $D4$ group with $R_i(q_l)=q_{l+i}$ (cyclic notation). 

Let $x\in V_*\setminus V_0$ and  $h\in\mathcal{H}_n(U(x))$.  Define $R(h)$ to be a function satisfying 
$$
R(h)(y)=\#W(x)^{-1}\sum_{w'\in W(x)} h\circ F_{w'} \circ R_{l'-l}\circ F_w^{-1}(y),
$$
for $y\in F_wK, w\in W(x)$.

For $i=1,2,3$, let $g_{x,i}$ be the local symmetry in $U(x)$, which is $F_w\circ g_{li}\circ F_w^{-1}$ on each component $F_wK$ with $x=F_w q_l$. 

Now,  for any $h\in \mathcal{H}_n(U(x))$, we could write $h=P_1(h)+P_2(h)+P_3(h)$ with 
\[
\begin{aligned}
P_1(h)&=\frac{1}{3}R(\sum_{i=1}^3 h\circ g_{x,i}),\\
P_2(h)&=\frac{1}{3}\sum_{i=1}^3 h\circ g_{x,i}-P_1(h),\\
P_3(h)&=h-P_1(h)-P_2(h).\\
\end{aligned}
\]

It is easy to check that for $i\leq n, w\in W(x)$,
\[
\begin{aligned}
\Delta^i P_k(h)(x)&=\delta_{k1}\Delta^i h(x),\\
\partial_n^w\Delta^i P_k(h)(x)&=\delta_{k2}\partial_n^w\Delta^i h(x),\\
\partial_{T,1}^w\Delta^i P_k(h)(x)=\delta_{k3} \partial_{T,1}^w\Delta^i h(x)&,\quad\partial_{T,2}^w\Delta^i P_k(h)(x)=\delta_{k3} \partial_{T,2}^w\Delta^i h(x).
\end{aligned}
\]

We need the definition of \textit{weak tangent} analogous to Definition 3.7.

\textbf{Definition 6.5.}($\mathcal{SG}^4$) \textit {Let $f$ be a  function which is continuous in a neighborhood of a vertex $x\in V_*\setminus V_0$. An $(n+1)$-harmonic function $h$ on $U(x)$ is said to be a weak tangent of order $n+1$ to $f$ at $x$ if
\begin{equation}
(f-h)|_{\partial U_m(x)}=o((\rho^nr)^m),
\end{equation}
and
\begin{equation}
\begin{cases}
((f-h)-(f-h)\circ g_{x,1})|_{\partial U_m(x)}=o((\rho^n\lambda)^m),\\
((f-h)-(f-h)\circ g_{x,2})|_{\partial U_m(x)}=o((\rho^n\lambda)^m).
\end{cases}
\end{equation}}
 
Now we could prove analogous results as stated in Theorem 3.8, for $\mathcal{SG}^4$, providing that all the numbers $\alpha_j,\beta_j,\gamma_j$ are not equal to $0$. The method is similar, except that when discussing the transverse derivatives, we need to look at the values $2P_3(h)(F_x^{m-|w|+i}F_wq_{l+1})-P_3(h)(F_x^{m-|w|+i}F_wq_{l+2})-P_3(h)(F_x^{m-|w|+i}F_wq_{l+3})$ and $2P_3(h)(F_x^{m-|w|+i}F_wq_{l+2})-P_3(h)(F_x^{m-|w|+i}F_wq_{l+1})-P_3(h)(F_x^{m-|w|+i}F_wq_{l+3})$ instead of the left side of $(3.13)$. For the calculations of $\alpha_j,\beta_j,\gamma_j$, see Appendix.

\section{Appendix}

As an appendix of this paper, we focus on the calculation of $\alpha_j$, $\beta_j$ and $\gamma_j$, the boundary values of the monomials $\{Q_{jk}^{(l)}\}$. We will mainly discuss the $D3$ symmetric fractals. The most typical example $\mathcal{SG}$ has been well studied in [NSTY], where an iterated calculation of the values as well as the derivatives of  $\{Q_{jk}^{(l)}\}$ at the boundary were given. However, their method is indirect, since it involves the boundary values and inner products of functions in the ``easy'' basis, and need to transform data from the ``easy'' basis to our ``monomial'' basis. Here we provide a new algorithm, which is more direct and shorter, using which, we could calculate $\alpha_j, \beta_j$ and $\gamma_j$ on some other examples, including $\mathcal{SG}_3$,$\mathcal{HG}$. Moreover, with some suitable modification, our algorithm will still be valid on $D4$ symmetric fractals. We will explain it on $\mathcal{SG}^4$.

Our approach is based on the relation of the Laplacian and the graph Laplacian of multihamonic functions, established in Lemma 4.3. 
Taking $m=1$ in (4.4), we get a recursion relation,
\begin{equation}
\begin{aligned}
\tilde{\Delta}_1 Q_{jk}^{(l)}(x)&=\sum_{i=1}^{j} \rho^{i-1}\alpha_1^{-1}\alpha_i\Delta^i Q_{jk}^{(l)}(x)\\
&=\sum_{i=1}^{j} \rho^{i-1}\alpha_1^{-1}\alpha_iQ_{(j-i)k}^{(l)}(x)
\end{aligned}
\end{equation}
holding at all vertices $x\in V_1\setminus V_0$ for all $j\geq 1$. Here in the second line of (7.1), we use the identity $\Delta^i Q_{jk}^{(l)}=Q_{(j-i)k}^{(l)}$. 

Thus, assuming we already have the values $\alpha_j,j\geq 0$, (7.1) as well as the self-similar identities (3.1)-(3.3) form a system of equations to calculate $Q_{jk}^{(l)}|_{V_1}$ from the values $Q_{ik}^{(l)}|_{V_1},0\leq i<j$. We do use this idea to solve the  $k\geq 2$ cases.

For $k=1$ case, it is a bit complicated, since we need to calculate all $\alpha_j$ simultaneously. We will give a theorem to show that $\alpha_{j}$ can be determined recursively  by using (7.1). 

For convenience of the readers, we first introduce the new calculation on $\mathcal{SG}$ as an example,  then give the proof for general $D3$ symmetric cases.

First we introduce some observations as well as some notations, some of which are same as those in [NSTY]. 

Simplify (7.1), we get
\begin{equation*}
\Delta_1 Q_{jk}^{(l)}(x)=(\sum_{y\sim_1 x} c_{xy})\sum_{i=1}^{j} \rho^i\alpha_i Q_{(j-i)k}^{(l)}(x).
\end{equation*}
Noticing that $\alpha_0=1$, we could rewrite the above identity into 
\begin{equation}
\sum_{y\sim_1 x} c_{xy}Q_{jk}^{(l)}(y)=(\sum_{y\sim_1 x} c_{xy})\sum_{i=0}^{j} \rho^i\alpha_i Q_{(j-i)k}^{(l)}(x)
\end{equation}
for $j\geq 1$. 

Also, we need a notation of \textit{infinite dimensional semi-circulant matrices}  $\alpha,\beta,\gamma$. For example, $\alpha=\{\alpha_{ij}\}_{i,j=0,1,2\cdots}$, has $\alpha_{ij}=\alpha_{i-j}$ for $i\geq j$ and $\alpha_{ij}=0$ for $i<j$.
It is easy to check  $(\alpha\beta)_{ij}=\sum_{l=i}^{j}\alpha_{il}\beta_{lj}=\sum_{l=0}^{i-j}\alpha_{l}\beta_{i-j-l}$ for $i\geq j$, and the multiplications among these matices are commutable. We will  need a linear operator $\tau$ on such matrices defined by
\begin{equation}
\tau\begin{pmatrix}    
     d_0&0& \\
     d_1&d_0&0&\\ \vspace{-0.2cm}
     d_2&d_1&d_0&0 &\\ \vspace{-0.1cm}
     d_3& d_2&d_1 &d_0 &\ddots\\
     \vdots & & & \ddots &\ddots
     \end{pmatrix}
=\begin{pmatrix}    
     d_0&0& \\
     \rho^{-1}d_1&d_0&0&\\ \vspace{-0.2cm}
     \rho^{-2}d_2&\rho^{-1}d_1&d_0&0 &\\ \vspace{-0.1cm}
     \rho^{-3}d_3& \rho^{-2}d_2&\rho^{-1}d_1 &d_0 &\ddots\\
     \vdots & & & \ddots &\ddots
     \end{pmatrix},
\end{equation}
where $\rho$ is the scaling constant of the Laplacian defined before.

\textbf{Example 7.1.} The monomials have been well studied in [NSTY], with $\alpha_j,\beta_j,\gamma_j$ exactly calculated. The resursion relations are
\[\begin{aligned}
\alpha_j&=\frac{4}{5^j-5}\sum_{i=1}^{j-1} \alpha_{j-i}\alpha_i, \forall j\geq 2,\\
\gamma_j&=\frac{4}{5^{j+1}-5}\sum_{i=0}^{j-1} \alpha_{j-i}\gamma_i, \forall  j\geq 1,\\
\beta_j=\frac{1}{5^j-1}\sum_{i=0}^{j-1}&(\frac{2}{5}5^{j-i}\alpha_{j-i}\beta_i-\frac{2}{3}\alpha_{j-i}5^i\beta_i+\frac{4}{5} \alpha_{j-i}\beta_i), \forall  j\geq 1,
\end{aligned}
\]
with initial data $\alpha_0=1, \alpha_1=1/6, \beta_0=-1/2, \gamma_0=1/2$. 

Now, we give a different calculation.

First, for $k=1$, by considering the symmetry, (7.2) becomes 
\[
\begin{cases}
\frac{a_{j}}{5^j}+\alpha_j+\frac{\alpha_j}{5^j}=4 \sum_{i=0}^{j}\frac{\alpha_i}{5^i}\frac{\alpha_{j-i}}{5^{j-i}},\\
2\alpha_j+\frac{2}{5^j}\alpha_j=4\sum_{i=0}^{j} \frac{\alpha_i}{5^i}\frac{a_{j-i}}{5^{j-i}},
\end{cases}
\]
for $j\geq 1$, where we denote $a_j=5^j Q_{j1}^{(0)}(F_1q_2)$. In addition, for $j=0$, we have
\[
\begin{cases}
a_{0}+2\alpha_0+1=4\alpha_0,\\
4\alpha_0=4 a_0.
\end{cases}
\] 
We could rewrite the above identities in matrix notation,
\begin{equation*}
\begin{cases}
a+\alpha+\tau(\alpha)+I=4\alpha^2,\\
2\alpha+2\tau(\alpha)=4\alpha a,
\end{cases}
\end{equation*}
by multiplying them with $5^j$ on both sides, where $a$ is the infinite matrix defined with $a_{ij}=a_{i-j}$ for $i\geq j$ and $a_{ij}=0$ for $i<j$.  Eliminating $a$, we get
\begin{equation*}
8\alpha^3-2\alpha^2-3\alpha=2\tau(\alpha)\alpha+\tau(\alpha),
\end{equation*}
which results that $\tau(\alpha)=4\alpha^2-3\alpha$. So we get the recursion relation for $\alpha_j$.

For $k=2$, we can write (7.2) into
\begin{equation*}
\begin{cases}
\frac{b_{j}}{5^j}+\frac{3\beta_j}{5^{j+1}}+\beta_j=4\sum_{i=0}^{j}\frac{\alpha_i}{5^i}\frac{3\beta_{j-i}}{5^{j-i+1}},\\
2\frac{3\beta_j}{5^{j+1}}+2\beta_j=4\sum_{i=0}^{j}\frac{\alpha_i}{5^i}\frac{b_{j-i}}{5^{j-i}},
\end{cases}
\end{equation*}
for $j\geq 0$, where we denote $b_j=5^jQ_{j2}^{(0)}(F_1q_2)$. Thus by multiplying both sides with $5^j$, we have
\begin{equation*}
\begin{cases}
b+\frac{3}{5}\beta+\tau(\beta)=\frac{12}{5}\alpha\beta,\\
\frac{6}{5}\beta+2\tau(\beta)=4\alpha b.
\end{cases}
\end{equation*}
With some calculation, we get
\[
\frac{3}{5}\beta(2\alpha-I)(4\alpha+I)=\tau(\beta)(2\alpha+I),
\]
which gives the recursion relation of $\beta_j$.

For $k=3$, we have $Q_{j3}^{(0)}(F_1q_2)=0$ by symmetry, so only one system of  equations need to be considered, which  immediately yields the recursion relation of $\gamma_j$.\\

Now, we turn to the general cases.

\textbf{Theorem 7.2.} \textit{Let $j\geq 2$.  Then $Q^{(l)}_{j1}|_{V_1}$ is uniquely determined by  the values of $Q^{(l)}_{i1}|_{V_1}$, $0\leq i<j$, by the following relations
\begin{equation}
\begin{cases}
\tilde{\Delta}_1 Q^{(l)}_{j1}(x)=\sum_{i=1}^{j}\rho^{i-1}\alpha_1^{-1}\alpha_iQ_{(j-i)1}^{(l)}(x),\forall x\in V_1\setminus V_0,\\
Q^{(l)}_{j1}|_{F_lV_0}=\rho^{j}Q^{(l)}_{j1}|_{V_0}.
\end{cases}
\end{equation}}

\textit{Proof.}  Obviously, $Q^{(l)}_{j1}|_{V_1}$ indeed satisfies the equations (7.4), which could be rewritten into an explicit form
\begin{equation*}
\begin{cases}
\tilde{\Delta}_1 Q^{(l)}_{j1}(x)-\rho^{j-1}\alpha_1^{-1}\alpha_j=\sum_{i=1}^{j-1}\rho^{i-1}\alpha_1^{-1}\alpha_iQ_{(j-i)1}^{(l)}(x),\forall x\in V_1\setminus V_0,\\
Q^{(l)}_{j1}|_{F_lV_0}=\rho^{j}Q^{(l)}_{j1}|_{V_0}.
\end{cases}
\end{equation*}

Thus, to prove that $Q^{(l)}_{j1}|_{V_1}$ is determined by (7.4) uniquely, we only need to prove the equations
\begin{equation}
\begin{cases}
\tilde{\Delta}_1 f(x)-\rho^{j-1}\alpha_1^{-1}f(q_{l+1})=0, \forall x\in V_1\setminus V_0,\\
f|_{F_lV_0}=\rho^{j}f|_{V_0},\\
f\circ g_l=f \mbox{ on } V_1
\end{cases}
\end{equation}
has a unique solution $f|_{V_1}=0$, where $g_l$ is the symmetry that fixes $q_l$ and interchanges the other two  vertices of $V_0$.

First we need to look at the equation
\begin{equation}
\tilde{\Delta}_1 h=1, h|_{V_0}=0.
\end{equation}
It is not hard to check that $h=(Q_{11}^{(l)}-h')|_{V_1}$ is the unique solution of (7.6), where $h'$ is the harmonic function with the same boundary values as those of $Q_{11}^{(l)}$, from which, one can find that 
\begin{equation}
h(F_lq_{l+1})=(\rho-r)\alpha_1.
\end{equation}

Now suppose $f$ is a solution of (7.5). Write $f=\rho^{j-1}\alpha_1^{-1}f(q_{l+1})h+\tilde{f}$. It is easy to check $\tilde{\Delta}_1 \tilde {f}=0$, and $f|_{V_0}=\tilde{f}|_{V_0}$. Moreover, by using (7.7) and $h|_{V_0}=0$, the relation $f|_{F_lV_0}=\rho^j f|_{V_0}$ implies
$$\tilde{f}\circ F_l(q_{l+1})+\rho^{j-1}f(q_{l+1})(\rho-r)=\rho^{j}f(q_{l+1}).$$
This could be simplified into 
$$rf(q_{l+1})+\rho^{j-1}f(q_{l+1})(\rho-r)=\rho^{j}f(q_{l+1}),$$
since $f|_{V_0}=\tilde{f}|_{V_0}$ and  $f\circ g_l=f$.
Thus we have 
$(r-\rho^{j-1}r)f(q_{l+1})=0$, which implies that $f(q_{l+1})=0$, and thus $f=\tilde{f}=0$ on $V_1$.

Hence we have proved the equations (7.5) only has a zero solution on $V_1$, which yields the result of the theorem. $\Box$

In the remaining section, we give the calculations of $\alpha_j$, $\beta_j$ and $\gamma_j$, as well as the numerical data, case by case, for $\mathcal{SG}_3$,$\mathcal{HG}$, and $\mathcal{SG}^4$.

\subsection{The level $3$ Sierpinski gasket $\mathcal{SG}_3$.} 

The first level graph of $\mathcal{SG}_3$ contains 7 inner vertices, and we take the following notations
\begin{equation*}
\rho^ja^{(k)}_j=Q^{(0)}_{jk}(F_3q_0), \rho^jb^{(k)}_j=Q^{(0)}_{jk}(F_1q_0), \rho^j c^{(k)}_j=Q^{(0)}_{jk}(F_1q_2).
\end{equation*}
The values of the functions $Q_{jk}^{(0)}$ on $V_1$ for $k=1,3$ are shown in Figure 7.1. Analogous to $\alpha, \beta, \gamma$, let $a^{(k)}, b^{(k)}, c^{(k)}$ be the infinite dimensional semi-circulant matrices generated by $\{a_j^{(k)}\}, \{b_j^{(k)}\}, \{c_j^{(k)}\}$.
\begin{figure}[h]
\begin{center}
\includegraphics[width=4.5cm,totalheight=4cm]{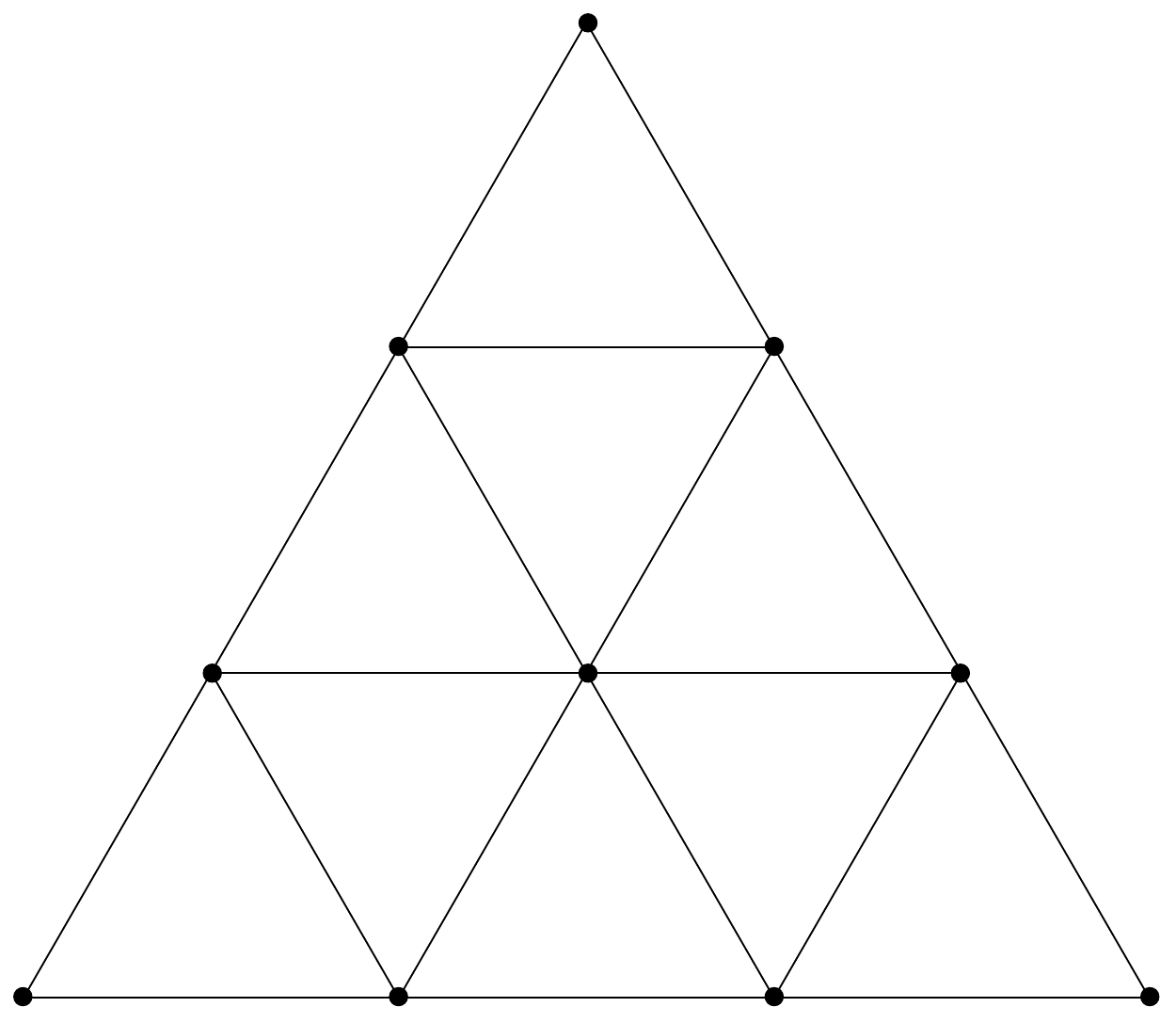}\hspace{2cm}
\includegraphics[width=4.5cm,totalheight=4cm]{sg3gamma1.pdf}
\setlength{\unitlength}{1cm}
\begin{picture}(0,0) \thicklines
\put(-11.7,-0){$\alpha_j$}
\put(-6.85,0){$\alpha_j$}
\put(-8.7,-0.3){$\rho^jc_j^{(1)}$}
\put(-10.15,-0.3){$\rho^jc_j^{(1)}$}
\put(-9.4,1.5){$\rho^ja_j^{(1)}$}
\put(-11.3,1.3){$\rho^jb_j^{(1)}$}
\put(-7.6,1.3){$\rho^jb_j^{(1)}$}
\put(-8.3,2.6){$\rho^j \alpha_j$}
\put(-10.55,2.6){$\rho^j \alpha_j$}
\put(-5.05,-0){$\gamma_j$}
\put(-0.25,0){$-\gamma_j$}
\put(-2.37,-0.3){$-\rho^jc_j^{(3)}$}
\put(-3.53,-0.3){$\rho^jc_j^{(3)}$}
\put(-2.58,1.5){$0$}
\put(-4.65,1.3){$\rho^jb_j^{(3)}$}
\put(-0.95,1.3){$-\rho^jb_j^{(3)}$}
\put(-1.7,2.6){$-\lambda\rho^j \gamma_j$}
\put(-4.1,2.6){$\lambda\rho^j \gamma_j$}

\put(-9.3,3.0){$F_0$}
\put(-2.7,3.0){$F_0$}
\put(-10.7,0.4){$F_1$}
\put(-4.1,0.4){$F_1$}
\put(-7.8,0.4){$F_2$}
\put(-1.2,0.4){$F_2$}
\put(-9.3,0.4){$F_3$}
\put(-2.7,0.4){$F_3$}
\put(-8.5,1.7){$F_4$}
\put(-1.9,1.7){$F_4$}
\put(-10,1.7){$F_5$}
\put(-3.4,1.7){$F_5$}
\end{picture}
\caption{The values of $Q_{j1}^{(0)}$ and $Q_{j3}^{(0)}$ on $V_1$. (The values of $Q_{j2}^{(0)}$ are similar to those of $Q_{j1}^{(0)}$.)}
\end{center}
\end{figure}

For $k=1$, by considering the symmetry, (7.2) becomes 
\[\begin{cases}
\rho^j\alpha_j+\rho^ja^{(1)}_j+\rho^jb^{(1)}_j=4\sum_{i=0}^{j} \rho^i\alpha_i\rho^{j-i}\alpha_{j-i},\\
\rho^j\alpha_j+\rho^jb^{(1)}_j+\rho^jc^{(1)}_j=3\sum_{i=0}^{j} \rho^i\alpha_i\rho^{j-i}a^{(1)}_{j-i},\\
\rho^j\alpha_j+\alpha_j+\rho^ja^{(1)}_j+\rho^jc^{(1)}_j=4\sum_{i=0}^{j} \rho^i\alpha_i\rho^{j-i}b^{(1)}_{j-i},\\
\alpha_j+\rho^ja^{(1)}_j+\rho^jb_j^{(1)}+\rho^jc^{(1)}_j=4\sum_{i=0}^{j} \rho^i\alpha_i\rho^{j-i}c^{(1)}_{j-i},
\end{cases}\]
 for $j\geq 1$. Combining with the $j=0$ case, we can rewrite them in matrix notation,
\[\begin{cases}
I+\alpha+a^{(1)}+b^{(1)}=4\alpha^2,\\
\alpha+b^{(1)}+c^{(1)}=3\alpha a^{(1)},\\
\alpha+\tau(\alpha)+a^{(1)}+c^{(1)}=4\alpha b^{(1)},\\
\tau(\alpha)+a^{(1)}+b^{(1)}+c^{(1)}=4\alpha c^{(1)}.
\end{cases}\]

For $k=2$, we write (7.2) into
\[\begin{cases}
r\rho^j\beta_j+\rho^ja^{(2)}_j+\rho^jb^{(2)}_j=4\sum_{i=0}^{j} \rho^i\alpha_ir\rho^{j-i}\beta_{j-i},\\
r\rho^j\beta_j+\rho^jb^{(2)}_j+\rho^jc^{(2)}_j=3\sum_{i=0}^{j} \rho^i\alpha_i\rho^{j-i}a^{(2)}_{j-i},\\
r\rho^j\beta_j+\beta_j+\rho^ja^{(2)}_j+\rho^jc^{(2)}_j=4\sum_{i=0}^{j} \rho^i\alpha_i\rho^{j-i}b^{(2)}_{j-i},\\
\beta_j+\rho^ja^{(2)}_j+\rho^jb_j^{(2)}+\rho^jc^{(2)}_j=4\sum_{i=0}^{j} \rho^i\alpha_i\rho^{j-i}c^{(2)}_{j-i},
\end{cases}\]
for all $j\geq 0$. Rewrite them in matrix notation, we get
\[\begin{cases}
r\beta+a^{(2)}+b^{(2)}=4r\alpha\beta,\\
r\beta+b^{(2)}+c^{(2)}=3\alpha a^{(2)},\\
r\beta+\tau(\beta)+a^{(2)}+c^{(2)}=4\alpha b^{(2)},\\
\tau(\beta)+a^{(2)}+b^{(2)}+c^{(2)}=4\alpha c^{(2)}.
\end{cases}\]

For $k=3$, we write (7.2) into
\[\begin{cases}
-\lambda\rho^j\gamma_j+\rho^jb^{(3)}_j=4\sum_{i=0}^{j} \rho^i\alpha_i\lambda\rho^{j-i}\gamma_{j-i},\\
\lambda\rho^j\gamma_j+\gamma_j+\rho^jc^{(3)}_j=4\sum_{i=0}^{j} \rho^i\alpha_i\rho^{j-i}b^{(3)}_{j-i},\\
\gamma_j+\rho^jb_j^{(3)}-\rho^jc^{(3)}_j=4\sum_{i=0}^{j} \rho^i\alpha_i\rho^{j-i}c^{(3)}_{j-i},
\end{cases}\]
for all $j\geq 0$. Then rewrite them in matrix notation,
\[\begin{cases}
-\lambda\gamma+b^{(3)}=4\lambda\alpha\gamma,\\
\lambda\gamma+\tau(\gamma)+c^{(3)}=4\alpha b^{(3)},\\
\tau(\gamma)+b^{(3)}-c^{(3)}=4\alpha c^{(3)}.
\end{cases}\]

The recursion relations of $\alpha,\beta,\gamma$ can be solved from the above equations by eliminating $a^{(k)}, b^{(k)}, c^{(k)}$. In matrix notation, the solutions are
\[\begin{aligned}
(1+6\alpha)\tau(\alpha)&=1+12\alpha-6\alpha^2-96\alpha^3+96\alpha^4,\\
(1+8\alpha+12\alpha^2)\tau(\beta)&=
(3+6\alpha-60\alpha^2-96\alpha^3+192\alpha^4)r\beta,\\
\tau(\gamma)&=(-1+16\alpha^2)\lambda\gamma.
\end{aligned}\] 
We can thus give the recursion relations of $\alpha_j,\beta_j,\gamma_j$ by using $r=7/15, \lambda=1/15$ and $\rho=7/90$.
\[
\begin{aligned}
\alpha_j=&\frac{1}{7(\frac{90}{7})^j-90}\left(96\sum_{i_1=1}^{j-1}\sum_{i_2=0}^{j-i_1}\sum_{i_3=0}^{j-i_1-i_2}\alpha_{i_1}\alpha_{i_2}\alpha_{i_3}\alpha_{j-i_1-i_2-i_3}\right.\\&\left.-6\sum_{i=1}^{j-1}(1+(\frac{90}{7})^{j-i})\alpha_i\alpha_{j-i}\right)\mbox{ for } j\geq 2,\\
\beta_j=&\frac{1}{15}\frac{1}{(\frac{90}{7})^j-1}\left(64\sum_{i_1=1}^{j}\sum_{i_2=0}^{j-i_1}\sum_{i_3=0}^{j-i_1-i_2}\sum_{i_4=0}^{j-i_1-i_2-i_3}\alpha_{i_1}\alpha_{i_2}\alpha_{i_3}\alpha_{i_4}\beta_{j-i_1-i_2-i_3-i_4}\right.\\
& +32\sum_{i_1=1}^{j}\sum_{i_2=0}^{j-i_1}\sum_{i_3=0}^{j-i_1-i_2}\alpha_{i_1}\alpha_{i_2}\alpha_{i_3}\beta_{j-i_1-i_2-i_3}\\&
\quad+\sum_{i_1=1}^{j}\sum_{i_2=0}^{j-i_1}(12-\frac{60}{7}(\frac{90}{7})^{j-i_1-i_2})\alpha_{i_1}\alpha_{i_2}\beta_{j-i_1-i_2}\\
&\left.\qquad +\sum_{i=1}^{j}(14-\frac{100}{7}(\frac{90}{7})^{j-i})\alpha_i\beta_{j-i}\right) \mbox{ for } j\geq 1,\\
\gamma_j=&\frac{16}{15}\frac{1}{(\frac{90}{7})^j-1}\left(\sum_{i_1=1}^{j}\sum_{i_2=0}^{j-i_1}\alpha_{i_1}\alpha_{i_2}\gamma_{j-i_1-i_2}+\sum_{i=1}^{j}\alpha_i\gamma_{j-i}\right) \mbox{ for } j\geq 1,
\end{aligned}
\] 
with initial data $\alpha_0=1, \alpha_1=1/6, \beta_0=-1/2, \gamma_0=1/2$.

\subsection{The Hexagasket $\mathcal{HG}$.}
The first level graph of $\mathcal{HG}$ contains 9 inner vertices. We take the following notations
\[
\rho^j a_j^{(k)}=Q_{jk}^{(0)}(F_4q_0), \rho^j b_j^{(k)}=Q_{jk}^{(0)}(F_1q_0), \rho^j c_j^{(k)}=Q_{jk}^{(0)}(F_1q_2), \rho^j d_j^{(k)}=Q_{jk}^{(0)}(F_3q_0).
\]

The values of the functions $Q_{jk}^{(0)}$ on $V_1$ for $k=1,3$ are shown in Figure 7.2. As before, let $a^{(k)}, b^{(k)}, c^{(k)}, d^{(k)}$ be the infinite dimensional semi-circulant matrices generated by $\{a_j^{(k)}\}, \{b_j^{(k)}\}, \{c_j^{(k)}\}, \{d_j^{(k)}\}$.

\begin{figure}[h]
\begin{center}
\includegraphics[width=4.5cm]{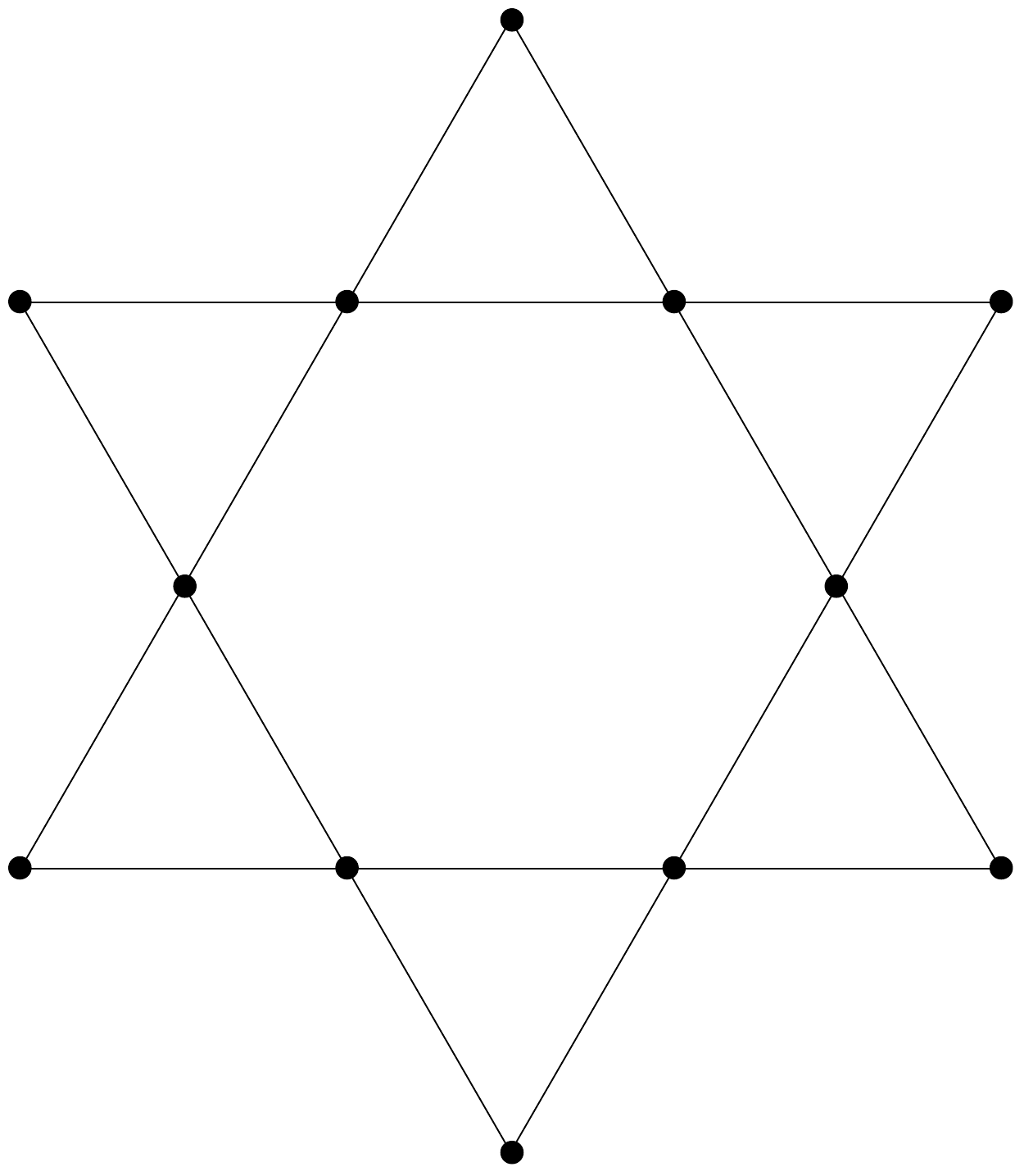}\hspace{2cm}
\includegraphics[width=4.5cm]{hexa1.pdf}
\setlength{\unitlength}{1cm}
\begin{picture}(0,0) \thicklines
\put(-11.7,1.2){$\alpha_j$}
\put(-6.85,1.2){$\alpha_j$}
\put(-12,3.7){$\rho^ja^{(1)}_j$}
\put(-6.85,3.7){$\rho^ja^{(1)}_j$}
\put(-10.5,3.95){$\rho^j\alpha_j$}
\put(-8.4,3.95){$\rho^j\alpha_j$}
\put(-11.3,2.55){$\rho^jb^{(1)}_j$}
\put(-7.6,2.55){$\rho^jb^{(1)}_j$}
\put(-10.6,1){$\rho^jc^{(1)}_j$}
\put(-8.4,1){$\rho^jc^{(1)}_j$}
\put(-9.45,-0.2){$\rho^jd^{(1)}_j$}

\put(-5.1,1.2){$\gamma_j$}
\put(-0.25,1.2){$-\gamma_j$}
\put(-5.4,3.7){$\rho^ja^{(3)}_j$}
\put(-0.25,3.7){$-\rho^ja^{(3)}_j$}
\put(-4.0,3.95){$\lambda\rho^j\gamma_j$}
\put(-1.8,3.95){$-\lambda\rho^j\gamma_j$}
\put(-4.7,2.55){$\rho^jb^{(3)}_j$}
\put(-1,2.55){$-\rho^jb^{(3)}_j$}
\put(-4,1){$\rho^jc^{(3)}_j$}
\put(-1.8,1){$-\rho^jc^{(3)}_j$}
\put(-2.55,-0.2){$0$}

\put(-9.25,4.1){$F_0$}
\put(-2.65,4.1){$F_0$}
\put(-10.65,1.7){$F_1$}
\put(-4.05,1.7){$F_1$}
\put(-7.8,1.7){$F_2$}
\put(-1.2,1.7){$F_2$}
\put(-9.3,0.7){$F_3$}
\put(-2.7,0.7){$F_3$}
\put(-10.75,3.2){$F_4$}
\put(-4.15,3.2){$F_4$}
\put(-7.85,3.2){$F_5$}
\put(-1.25,3.2){$F_5$}
\end{picture}
\caption{The values of $Q_{j1}^{(0)}$ and $Q_{j3}^{(0)}$ on $V_1$. (The values of $Q_{j2}^{(0)}$ are similar to those of $Q_{j1}^{(0)}$.)}
\end{center}
\end{figure}

Now we rewrite the equation (7.2) for all $j$ using the above notations.

For $k=1$, by considering the symmetry, (7.2) becomes 
\[\begin{cases}
\rho^j\alpha_j+\rho^ja^{(1)}_j+\rho^jb^{(1)}_j=4\sum_{i=0}^{j} \rho^i\alpha_i\rho^{j-i}\alpha_{j-i},\\
\rho^j\alpha_j+\rho^jb^{(1)}_j=2\sum_{i=0}^{j} \rho^i\alpha_i\rho^{j-i}a_{j-i}^{(1)},\\
\alpha_j+\rho^j\alpha_j+\rho^ja^{(1)}_j+\rho^jc^{(1)}_j=4\sum_{i=0}^{j} \rho^i\alpha_i\rho^{j-i}b_{j-i}^{(1)},\\
\alpha_j+\rho^jb^{(1)}_j+\rho^jc^{(1)}_j+\rho^jd^{(1)}_j=4\sum_{i=0}^{j} \rho^i\alpha_i\rho^{j-i}c_{j-i}^{(1)},\\
2\rho^jc^{(1)}_j=2\sum_{i=0}^{j} \rho^i\alpha_i\rho^{j-i}d_{j-i}^{(1)},
\end{cases}\]
for $j\geq 1$. Combining with the $j=0$ case, we can rewrite them in matrix notation,
\[\begin{cases}
I+\alpha+a^{(1)}+b^{(1)}=4\alpha^2,\\
\alpha+b^{(1)}=2\alpha a^{(1)},\\
\alpha+\tau(\alpha)+a^{(1)}+c^{(1)}=4\alpha b^{(1)},\\
\tau(\alpha)+b^{(1)}+c^{(1)}+d^{(1)}=4\alpha c^{(1)},\\
c=\alpha d^{(1)}.
\end{cases}\]

For $k=2$, we can write (7.2) into
\[\begin{cases}
r\rho^j\beta_j+\rho^ja^{(2)}_j+\rho^jb^{(2)}_j=4\sum_{i=0}^{j} \rho^i\alpha_ir\rho^{j-i}\beta_{j-i},\\
r\rho^j\beta_j+\rho^jb^{(2)}_j=2\sum_{i=0}^{j} \rho^i\alpha_i\rho^{j-i}a_{j-i}^{(2)},\\
\beta_j+r\rho^j\beta_j+\rho^ja^{(2)}_j+\rho^jc^{(2)}_j=4\sum_{i=0}^{j} \rho^i\alpha_i\rho^{j-i}b_{j-i}^{(2)},\\
\beta_j+\rho^jb^{(2)}_j+\rho^jc^{(2)}_j+\rho^jd^{(2)}_j=4\sum_{i=0}^{j} \rho^i\alpha_i\rho^{j-i}c_{j-i}^{(2)},\\
2\rho^jc^{(2)}_j=2\sum_{i=0}^{j} \rho^i\alpha_i\rho^{j-i}d_{j-i}^{(2)},
\end{cases}\]
for all $j\geq 0$. Rewrite them in matrix notation, we get
\[\begin{cases}
r\beta+a^{(2)}+b^{(2)}=4r\alpha\beta,\\
r\beta+b^{(2)}=2\alpha a^{(2)},\\
r\beta+\tau(\beta)+a^{(2)}+c^{(2)}=4\alpha b^{(2)},\\
\tau(\beta)+b^{(2)}+c^{(2)}+d^{(2)}=4\alpha c^{(2)},\\
c^{(2)}=\alpha d^{(2)}.
\end{cases}\]

For $k=3$,  (7.2) becomes
\[\begin{cases}
-\lambda\rho^j\gamma_j+\rho^ja^{(3)}_j+\rho^jb^{(3)}_j=4\sum_{i=0}^{j} \rho^i\alpha_i\lambda\rho^{j-i}\gamma_{j-i},\\
\lambda\rho^j\gamma_j+\rho^jb^{(3)}_j=2\sum_{i=0}^{j} \rho^i\alpha_i\rho^{j-i}a^{(3)}_{j-i},\\
\gamma_j+\lambda\rho^j\gamma_j+\rho^ja^{(3)}_j+\rho^jc^{(3)}_j=4\sum_{i=0}^{j} \rho^i\alpha_i\rho^{j-i}b_{j-i}^{(3)},\\
\gamma_j+\rho^jb^{(3)}_j-\rho^jc^{(3)}_j=4\sum_{i=0}^{j} \rho^i\alpha_i\rho^{j-i}c_{j-i}^{(3)},
\end{cases}\]
for all $j\geq 0$. Then rewrite them in matrix notation,
\[\begin{cases}
-\lambda\gamma+a^{(3)}+b^{(3)}=4\lambda\alpha\gamma,\\
\lambda\gamma+b^{(3)}=2\alpha a^{(3)},\\
\tau(\gamma)+\lambda\gamma+a^{(3)}+c^{(3)}=4\alpha b^{(3)},\\
\tau(\gamma)+b^{(3)}-c^{(3)}=4\alpha c^{(3)}.
\end{cases}\]

By eliminating $a^{(k)}, b^{(k)}, c^{(k)}, d^{(k)}$ from the above equations, we get the following recursion relations of $\alpha, \beta, \gamma$ in matrix notation, 
\[\begin{aligned}
(-1+2\alpha)\tau(\alpha)&=-1+4\alpha+14\alpha^2-48\alpha^3+32\alpha^4,\\
(-1+4\alpha^2)\tau(\beta)&=r(1+10\alpha-4\alpha^2-64\alpha^3+64\alpha^4)\beta,\\
\tau(\gamma)&=\lambda(-1-8\alpha+16\alpha^2)\gamma.
\end{aligned}\]
Similarly, substituting $r=3/7, \lambda=1/7, \rho=1/14$, we have the recursion relations in the explicit form
\[\begin{aligned}
\alpha_j=&\frac{1}{14^j-14}\left(32\sum_{i_1=1}^{j-1}\sum_{i_2=0}^{j-i_1}\sum_{i_3=0}^{j-i_1-i_2}\alpha_{i_1}\alpha_{i_2}\alpha_{i_3}\alpha_{j-i_1-i_2-i_3}\right.\\&-16\sum_{i_1=1}^{j-1}\sum_{i_2=0}^{j-i_1}\alpha_{i_1}\alpha_{i_2}\alpha_{j-i_1-i_2}
\left.-2\sum_{i=1}^{j-1}(1+14^{j-i})\alpha_i\alpha_{j-i}\right)\mbox{ for } j\geq 2,\\
\beta_j=&\frac{1}{14^j-1}\left(\frac{64}{7}\sum_{i_1=1}^{j}\sum_{i_2=0}^{j-i_1}\sum_{i_3=0}^{j-i_1-i_2}\sum_{i_4=0}^{j-i_1-i_2-i_3}\alpha_{i_1}\alpha_{i_2}\alpha_{i_3}\alpha_{i_4}\beta_{j-i_1-i_2-i_3-i_4}\right.\\
&-\sum_{i_1=1}^{j}\sum_{i_2=0}^{j-i_1}(\frac{4}{7}+\frac{4}{3} \cdot14^{j-i_1-i_2})\alpha_{i_1}\alpha_{i_2}\beta_{j-i_1-i_2}\\&\quad\left.+\sum_{i=1}^j(\frac{6}{7}-\frac{4}{3}\cdot 14^{j-i})\alpha_i\beta_{j-i} \right)\mbox{ for } j\geq 1,\\
\gamma_j=&\frac{1}{14^j-1}\left(\frac{16}{7}\sum_{i_1=1}^{j}\sum_{i_2=0}^{j-i_1}\alpha_{i_1}\alpha_{i_2}\gamma_{j-i_1-i_2}+\frac{8}{7}\sum_{i=1}^{j}\alpha_i\gamma_{j-i}\right)\mbox{ for } j\geq 1,
\end{aligned}
\]

It is easy to check that the initial data are $\alpha_0=1, \alpha_1=1/6, \beta_0=-1/2, \gamma_0=1/2$ as before. 

\subsection{The tetrahedral Sierpinski gasket $\mathcal{SG}^4$.}

For convenience, denote $\rho^j a_j^{(k)}=Q^{(0)}_{jk}(F_1q_2)$ and let $a^{(k)}$ denote the infinite dimensional semi-circulant matrices generated by $\{a_{j}^{(k)}\}$. See the values  of $Q_{j1}^{(0)}$ and $Q_{j3}^{(0)}$ on $V_1$ in Fig. 7.3. We should point out that the identity (7.2) still holds for $\mathcal{SG}^4$, since the harmonic extension matrices are nondegenerate  for $\mathcal{SG}^4$ and hence the argument for proving Lemma 4.3 remains to be valid.

\begin{figure}[h]
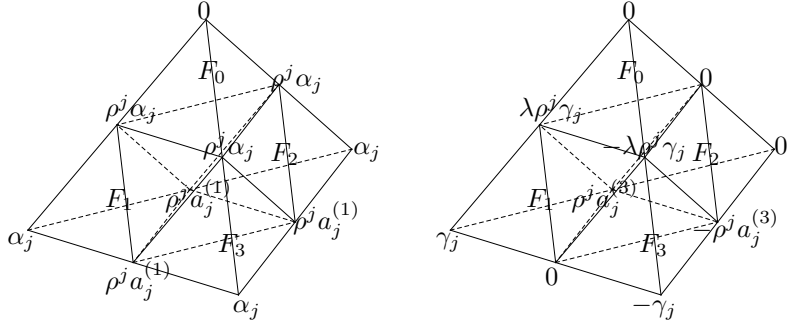

\begin{center}
\includegraphics[width=4.5cm]{sg42.pdf}\hspace{1cm}
\includegraphics[width=4.5cm]{sg42.pdf}
\setlength{\unitlength}{1cm}

\begin{picture}(0,0) \thicklines
\put(-5.3,1.2){$\alpha_j$}
\put(-2.3,0.3){$\alpha_j$}
\put(-0.7,2.4){$\alpha_j$}
\put(-4,2.9){$\rho^j\alpha_j$}
\put(-2.7,2.4){$\rho^j\alpha_j$}
\put(-1.8,3.3){$\rho^j\alpha_j$}
\put(-2.77,4.2){$0$}
\put(-4,0.6){$\rho^j a^{(1)}_j$}
\put(-3.2,1.7){$\rho^j a^{(1)}_j$}
\put(-1.5,1.4){$\rho^j a^{(1)}_j$}

\put(0.4,1.2){$\gamma_j$}
\put(3.0,0.3){$-\gamma_j$}
\put(1.85,0.65){$0$}
\put(1.5,2.9){$\lambda\rho^j\gamma_j$}
\put(2.6,2.4){$-\lambda\rho^j\gamma_j$}
\put(2.2,1.7){$\rho^ja^{(3)}_j$}
\put(3.8,1.3){$-\rho^ja^{(3)}_j$}
\put(2.88,4.2){$0$}
\put(3.9,3.3){$0$}
\put(4.9,2.4){$0$}

\put(-2.77,3.4){$F_0$}
\put(-4.0,1.7){$F_1$}
\put(-1.8,2.3){$F_2$}
\put(-2.5,1.1){$F_3$}
\put(2.8,3.4){$F_0$}
\put(1.6,1.7){$F_1$}
\put(3.8,2.3){$F_2$}
\put(3.1,1.1){$F_3$}
\end{picture}

\caption{The values of $Q_{j1}^{(0)}$ and $Q_{j3}^{(0)}$ on $V_1$. }
\end{center}
\end{figure}

For $k=1$, by considering the symmetry, (7.2) becomes 
\[\begin{cases}
\alpha_j+2\rho^j\alpha_j+2\rho^ja^{(1)}_j=6\sum_{i=0}^{j} \rho^i\alpha_i\rho^{j-i}\alpha_{j-i},\\
2\alpha_j+2\rho^j\alpha_j+2\rho^ja^{(1)}_j=6\sum_{i=0}^{j} \rho^i\alpha_i\rho^{j-i}a_{j-i}^{(1)},
\end{cases}\]
for $j\geq 1$. Combining with the $j=0$ case, we can rewrite them in matrix notation,
\[\begin{cases}
I+\tau(\alpha)+2\alpha+2a^{(1)}=6\alpha^2,\\
\tau(\alpha)+\alpha+a^{(1)}=3\alpha a^{(1)}.
\end{cases}\]

For $k=2$, we can write (7.2) into
\[\begin{cases}
\beta_j+2r\rho^j\beta_j+2\rho^ja^{(2)}_j=6\sum_{i=0}^{j} \rho^i\alpha_ir\rho^{j-i}\beta_{j-i},\\
2\beta_j+2r\rho^j\beta_j+2\rho^ja^{(2)}_j=6\sum_{i=0}^{j} \rho^i\alpha_i\rho^{j-i}a_{j-i}^{(2)},
\end{cases}\]
for all $j\geq 0$. Rewrite them in matrix notation, we get
\[\begin{cases}
\tau(\beta)+2r\beta+2a^{(2)}=6r\alpha\beta,\\
\tau(\beta)+r\beta+a^{(2)}=3\alpha a^{(2)}.
\end{cases}\]

For $k=3$, we can write (7.2) into
\[\begin{cases}
\gamma_j-\lambda\rho^j\gamma_j+\rho^ja^{(3)}_j=6\sum_{i=0}^{j} \rho^i\alpha_i\lambda\rho^{j-i}\gamma_{j-i},\\
\gamma_j+\lambda\rho^j\gamma_j-\rho^ja^{(3)}_j=6\sum_{i=0}^{j} \rho^i\alpha_i\rho^{j-i}a^{(3)}_{j-i},
\end{cases}\]
for all $j\geq 0$. Then rewrite them in matrix notation,
\[\begin{cases}
\tau(\gamma)-\lambda\gamma+a^{(3)}=6\lambda\alpha\gamma,\\
\tau(\gamma)+\lambda\gamma-a^{(3)}=6\alpha a^{(3)}.
\end{cases}\]

By eliminating $a^{(k)}$ in the above matrix equations, we could solve the recursion relations of $\alpha,\beta,\gamma$ as follows, in the matrix notation,
\[
\begin{aligned}
\tau(\alpha)&=1 + 6 \alpha^2-6\alpha,\\
(1+3\alpha)\tau(\beta)&=r(-12\alpha+18\alpha^2)\beta,\\
\tau(\gamma)&=6\lambda\alpha\gamma.
\end{aligned}
\]
Using the fact that $r=2/3, \lambda=1/6,\rho=1/6$, the recursion relations of $\alpha_j,\beta_j,\gamma_j$ in the explicit form are
\[
\begin{aligned}
\alpha_j=&\frac{1}{6^{j-1}-1}\sum_{i=1}^{j-1}\alpha_i\alpha_{j-i}\mbox{ for } j\geq 2,\\
\beta_j=&\frac{1}{6^j-1}\left(3\sum_{i_1=1}^{j}\sum_{i_2=0}^{j-i_1}\alpha_{i_1}\alpha_{i2}\beta_{j-i_1-i_2}+\sum_{i=1}^{j}(1-\frac{3}{4}\cdot 6^{j-i})\alpha_i\beta_{j-i}\right)\mbox{ for } j\geq 1\\
\gamma_j=&\frac{1}{6^j-1}\sum_{i=1}^j \alpha_i\gamma_{j-i} \mbox{ for } j\geq 1.
\end{aligned},
\]

To get the initial data, we need to look at Lemma 4.3 again. In fact, we have
\[
1=\tilde{\Delta}_1 h(x)=\frac{3r^{-1}\#W(x) }{\int \psi_x^1d\mu}\rho\alpha_1,
\]
for $x\in V_1\setminus V_0$, for any $h\in \mathcal{H}_1$ with $\Delta h(x)=1$. This gives that $\alpha_1=1/12$. Thus we have initial data $\alpha_0=1, \alpha_1=1/12, \beta_0=-1/3,\gamma_0=1/3.$ In addition, it is easy to check that $\gamma_j=4\alpha_{j+1}$.

Table 7.1-7.4 present numerical computations of $\alpha_j, \beta_j$ and $\gamma_j$ for $\mathcal{SG}_3$, $\mathcal{HG}$ and $\mathcal{SG}^4$, respectively. (We are grateful to Mr. Wei Wei for providing an effective program.) For $\mathcal{SG}_3$ and $\mathcal{HG}$, it is easy to find that $\alpha_j,\beta_j$ behave similar to geometric progressions when $j$ is large enough, with the reciprocal of common ratio $-124.68442107\cdots$ for $\mathcal{SG}_3$  and $-46.728917838\cdots$ for $\mathcal{HG}$.  As for $\mathcal{SG}^4$, it is quite similar to $\mathcal{SG}$ case, calculated in [NSTY]. All $\alpha_j$ and $\gamma_j$ take positive values, and $\beta_j$ behaves like a geometric progression when $j\geq 11$, with the reciprocal of common ratio being $-338.81012588\cdots$. In fact, for each case, the above mentioned reciprocal of ratio is the opposite of the least  eigenvalue of $-\Delta$,  corresponding to  an eigenfunction, which is symmetric with respect to the reflection fixing $q_0$ with both the value and the normal derivative at $q_0$ vanishing. An explanation of this phenomenon comes from a slight generalization of Theorem 2.9 in [NSTY] (for $\mathcal{SG}$) involving a rather detailed knowledge of the description of eigenfunctions of $-\Delta$ by the spectral decimation. In fact, the natural of eigenvalues and eigenfunctions could be known explicitly via the method of spectral decimation for some fully symmetric p.c.f. fractals(See [FS, MT, Sh1-Sh2, ST,T1]). We refer to the reader to find the spectral decimation recipes for $\mathcal{SG}_3$, $\mathcal{HG}$ and $\mathcal{SG}^4$ in [DS], [BCDEHKMST] and [FS] respectively, 
using which, we could verify that $124.68442107\cdots$ is the eigenvalue of $-\Delta$ on $\mathcal{SG}_3$ with eigenfunction shown in Fig. 7.4(a), $46.728917838\cdots$ is the eigenvalue on $\mathcal{HG}$ with eigenfunction shown in Fig. 7.4(b),  $338.81012588\cdots$ is the second Neumann eigenvalue on $\mathcal{SG}^4$ with eigenfunction shown in Fig. 7.4(c) as the $\mathcal{SG}$ case.

\begin{table}[htbp]
\caption{The data of $\alpha_j,\beta_j,\gamma_j$ for $\mathcal{SG}_3$.}
\centering
\begin{tabular}{l  l  l  l}
\hline
j&$\alpha_j$&$\beta_{j}$&$\gamma_{j}$\\
\hline
0 & 1                             & -0.5000000000                  & 0.5000000000                  \\
1 & 0.1666666667                  & -0.04645247657                 & 0.01499330656                 \\
2 & 0.5332440874$\times 10^{-2}$  & -0.1029307014$\times 10^{-2}$  & 0.1572291087$\times 10^{-2}$  \\
3 & 0.5981637501$\times 10^{-4}$  & -0.8969605760$\times 10^{-5}$  & 0.7920041236$\times 10^{-6}$  \\
4 & 0.3116779311$\times 10^{-6}$  & -0.3931817387$\times 10^{-7}$  & 0.2312923424$\times 10^{-8}$  \\
5 & 0.9411222994$\times 10^{-9}$  & -0.1015006336$\times 10^{-9}$  & 0.4359699965$\times 10^{-11}$ \\
6 & 0.1768209338$\times 10^{-11}$ & -0.1801026965$\times 10^{-12}$ & 0.5680708143$\times 10^{-14}$ \\
7 & 0.2703155181$\times 10^{-14}$ & -0.1686026983$\times 10^{-15}$ & 0.5372013243$\times 10^{-17}$ \\
8 & -0.9794948608$\times 10^{-19}$& -0.5193342268$\times 10^{-18}$ & 0.3827057814$\times 10^{-20}$ \\
9 & 0.2103310426$\times 10^{-19}$ & 0.2555737442$\times 10^{-20}$  & 0.2116632464$\times 10^{-23}$ \\
10& -0.1542910507$\times 10^{-21}$& -0.2156048818$\times 10^{-22}$ & 0.9317275066$\times 10^{-27}$ \\
11& 0.1245430189$\times 10^{-23}$ & 0.1723779132$\times 10^{-24}$  & 0.3333118088$\times 10^{-30}$ \\
12& -0.9985135547$\times 10^{-26}$& -0.1382754075$\times 10^{-26}$ & 0.9861897696$\times 10^{-34}$ \\
13& 0.8008453196$\times 10^{-28}$ & 0.1108998433$\times 10^{-28}$  & 0.2449517096$\times 10^{-37}$ \\
14& -0.6422974455$\times 10^{-30}$& -0.8894451991$\times 10^{-31}$ & 0.5172684514$\times 10^{-41}$ \\
15& 0.5151385008$\times 10^{-32}$ & 0.7133572856$\times 10^{-33}$  & 0.9388504271$\times 10^{-45}$ \\
16& -0.4131538617$\times 10^{-34}$& -0.5721302858$\times 10^{-35}$ & 0.1478428756$\times 10^{-48}$ \\
17& 0.3313596503$\times 10^{-36}$ & 0.4588626969$\times 10^{-37}$  & 0.2036817659$\times 10^{-52}$ \\
18& -0.2657586629$\times 10^{-38}$& -0.3680192727$\times 10^{-39}$ & 0.2471289873$\times 10^{-56}$ \\
19& 0.2131450430$\times 10^{-40}$ & 0.2951605903$\times 10^{-41}$  & 0.2667519036$\times 10^{-60}$ \\
\hline
\end{tabular}\end{table}

\begin{table}[htbp]
\caption{The data of $\alpha_j,\beta_j,\gamma_j$ for $\mathcal{HG}$.}
\centering\begin{tabular}{l l l l}
\hline
j&$\alpha_j$&$\beta_{j}$&$\gamma_{j}$\\
\hline
0 & 1                             & -0.5000000000                  & 0.5000000000                   \\
1 & 0.1666666667                  & -0.04334554334                 & 0.02197802198                  \\
2 & 0.00518925518                 & -0.0008741066739               & 0.0002728244486                \\
3 & 0.4189271589$\times 10^{-4}$  & -0.5294515600$\times 10^{-5}$  & 0.14549879455$\times 10^{-5}$  \\
4 & 0.3320775837$\times 10^{-6}$  & -0.4109349118$\times 10^{-7}$  & 0.41744883710$\times 10^{-8}$  \\
5 & -0.1983647549$\times 10^{-8}$ & 0.3457846603$\times 10^{-9}$   & 0.73217281196$\times 10^{-11}$ \\
6 & 0.5477498983$\times 10^{-10}$ & -0.8500346754$\times 10^{-11}$ & 0.85864350196$\times 10^{-14}$ \\
7 & -0.1153381369$\times 10^{-11}$& 0.1804328531$\times 10^{-12}$  & 0.71435664780$\times 10^{-17}$ \\
8 & 0.2470213442$\times 10^{-13}$ & -0.3862665725$\times 10^{-14}$ & 0.44009750989$\times 10^{-20}$ \\
9 & -0.5286113716$\times 10^{-15}$& 0.8266021437$\times 10^{-16}$  & 0.20790960455$\times 10^{-23}$ \\
10& 0.1131230500$\times 10^{-16}$ & -0.1768931223$\times 10^{-17}$ & 0.75893835050$\times 10^{-27}$ \\
11& -0.2420836036$\times 10^{-18}$& 0.3785517208$\times 10^{-19}$  & 0.24618203914$\times 10^{-30}$ \\
12& 0.5180595119$\times 10^{-20}$ & -0.8101016218$\times 10^{-21}$ & 0.19858084061$\times 10^{-34}$ \\
13& -0.1108648640$\times 10^{-21}$& 0.1733619479$\times 10^{-22}$  & 0.76041321900$\times 10^{-37}$ \\
14& 0.2372510837$\times 10^{-23}$ & -0.3709949982$\times 10^{-24}$ & -0.11425940714$\times 10^{-39}$\\
15& -0.5077179071$\times 10^{-25}$& 0.7939302157$\times 10^{-26}$  & 0.20348945233$\times 10^{-42}$ \\
16& 0.1086517580$\times 10^{-26}$ & -0.1699012630$\times 10^{-27}$ & -0.34990746951$\times 10^{-45}$\\
17& -0.2325150314$\times 10^{-28}$& 0.3635891238$\times 10^{-29}$  & 0.59498904116$\times 10^{-48}$ \\
18& 0.4975827434$\times 10^{-30}$ & -0.7780816262$\times 10^{-31}$ & -0.10013087692$\times 10^{-50}$\\
19& -0.1064828304$\times 10^{-31}$& 0.1665096609$\times 10^{-32}$  & 0.16709434786$\times 10^{-53}$\\
\hline
\end{tabular}
\end{table}

\begin{table}[htbp]
\caption{The data of $\alpha_j,\beta_j$ for $\mathcal{SG}^4$.}
\begin{tabular}{l l l}
\hline
j&$\alpha_j$&$\beta_{j}$\\
\hline
0 & 1                             & -0.3333333333                         \\
1 & 0.08333333333                 & -0.01805555555                        \\
2 & 0.001388888889                & -0.0002199074074                      \\
3 & 0.6613756614$\times 10^{-5}$  & -0.8408186703$\times 10^{-6}$         \\
4 & 0.1409909356$\times 10^{-7}$  & -0.1553094253$\times 10^{-8}$         \\
5 & 0.1600107728$\times 10^{-10}$ & -0.1516427881$\times 10^{-11}$        \\
6 & 0.1100614425$\times 10^{-13}$ & -0.1060002208$\times 10^{-14}$        \\
7 & 0.4989341781$\times 10^{-17}$ & -0.1156439717$\times 10^{-18}$        \\
8 & 0.1578377410$\times 10^{-20}$ & -0.9588249504$\times 10^{-21}$        \\
9 & 0.3637188495$\times 10^{-24}$ & 0.2457583671$\times 10^{-23}$         \\
10& 0.6319196749$\times 10^{-28}$ & -0.7332537520$\times 10^{-26}$        \\
11& 0.8514009570$\times 10^{-32}$ & 0.2162924611$\times 10^{-29}$         \\
12& 0.9104272572$\times 10^{-36}$ & -0.6384045750$\times 10^{-31}$        \\
13& 0.7876496987$\times 10^{-40}$ & 0.1884253065$\times 10^{-33}$         \\
14& 0.5602617844$\times 10^{-44}$ & -0.5561383694$\times 10^{-36}$        \\
15& 0.3321803330$\times 10^{-48}$ & 0.1641445537$\times 10^{-38}$         \\
16& 0.1661241248$\times 10^{-52}$ & -0.4844735773$\times 10^{-41}$        \\
17& 0.7080859932$\times 10^{-57}$ & 0.1429926499$\times 10^{-43}$         \\
18& 0.2596231219$\times 10^{-61}$ & -0.4220436139$\times 10^{-46}$        \\
19& 0.8256443790$\times 10^{-66}$ & 0.1245664110$\times 10^{-48}$      \\
\hline   
\end{tabular}
\end{table}

\begin{table}[htbp]
\caption{The data of ratios of $\alpha_j,\beta_j$.}
\centering
\begin{tabular}{l | l l | l l| l}
\hline
j&$\alpha_{j-1}/\alpha_j(\mathcal{SG}_3)$&$\beta_{j-1}/\beta_{j}(\mathcal{SG}_3)$&$\alpha_{j-1}/\alpha_j(\mathcal{HG})$&$\beta_{j-1}/\beta_{j}(\mathcal{HG})$&$\beta_{j-1}/\beta_{j}(\mathcal{SG}^4)$\\
\hline
0 & /            & /            &/             & /            & / \\
1 & 6            & 10.76368876  & 6            & 11.53521127  & 18.46153846\\
2 & 31.25523013  & 45.12985529  & 32.11764706  & 49.58839080  & 82.10526316\\
3 & 89.14684104  & 114.7549894  & 123.8701068  & 165.0966282  & 261.5396341\\
4 & 191.9172615  & 228.1287475  & 126.1533989  & 128.8407348  & 541.3828996\\
5 & 331.1768633  & 387.3687528  & -167.4075537 & -118.8412787 & 1024.179437\\
6 & 532.2459729  & 563.5708715  & -36.21447589 & -40.67888879 & 1430.589360\\
7 & 654.1279429  & 1068.207675  & -47.49078777 & -47.11085928 & 9166.082698\\
8 & -2759.744119 & 324.6516205  & -46.69156721 & -46.71200305 & 120.6100984\\
9 & -4.656920103 & -203.2032783 & -46.73023651 & -46.72944239 & -390.1494634\\
10& -136.3209607 & -118.5380136 & -46.72888255 & -46.72890235 & -335.1614178\\
11& -123.8857481 & -125.0768604 & -46.72891855 & -46.72891775 & -339.0103141\\
12& -124.7284209 & -124.6627410 & -46.72891783 & -46.72891778 & -338.8015525\\
13& -124.6824487 & -124.6849440 & -46.72891784 & -46.72891783 & -338.8104214\\
14& -124.6844940 & -124.6842902 & -46.72891784 &-46.72891784  & -338.8101179\\
15& -124.6844188 & -124.6843927 & -46.72891784 &-46.72891784  & -338.8101261\\
16& -124.6844211 & -124.6844125 & -46.72891784 &-46.72891784  & -338.8101259\\
17& -124.6844211 & -124.6844186 & -46.72891784 &-46.72891784  & -338.8101259\\
18& -124.6844211 & -124.6844204 & -46.72891784 &-46.72891784  & -338.8101259\\
19& -124.6844211 & -124.6844209 & -46.72891784 &-46.72891784  & -338.8101259\\
\hline
\end{tabular}
\end{table}

\begin{figure}[h]
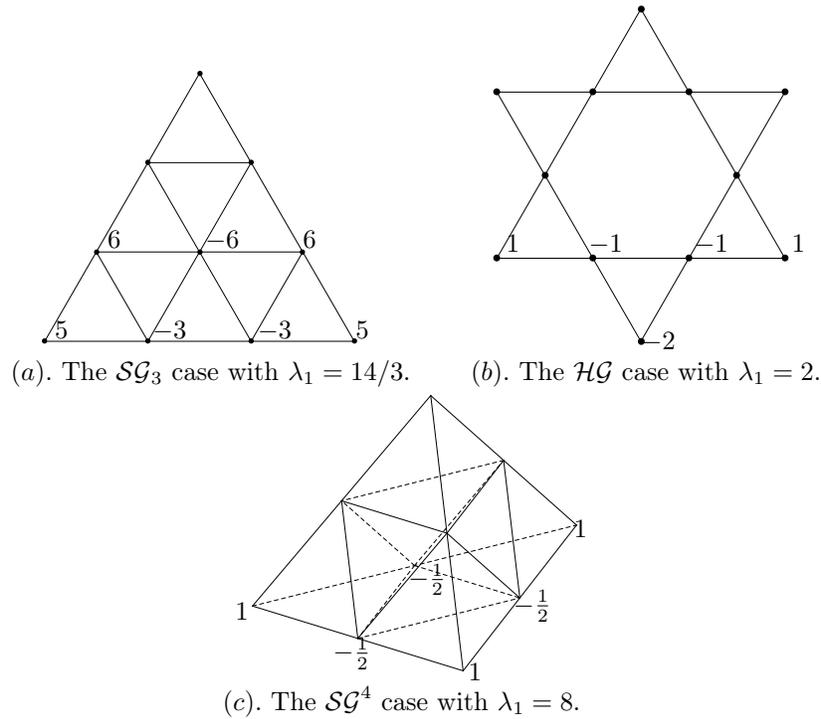

\begin{center}
\includegraphics[width=4.3cm]{sg3gamma1.pdf}\hspace{1.6cm}
\includegraphics[width=4.0cm]{hexa1.pdf}

$(a).$ The $\mathcal{SG}_3$ case with $\lambda_1=14/3$.\hspace{0.8cm}$(b).$ The $\mathcal{HG}$ case with $\lambda_1=2$.

\includegraphics[width=4.5cm]{sg42.pdf}

$(c).$ The $\mathcal{SG}^4$ case with $\lambda_1=8$.
\setlength{\unitlength}{1cm}
\begin{picture}(0,0) \thicklines
\put(-7.3,4.95){$5$}
\put(-3.3,4.95){$5$}
\put(-4.6,4.95){$-3$}
\put(-6.0,4.95){$-3$}
\put(-5.3,6.15){$-6$}
\put(-6.6,6.15){$6$}
\put(-4.0,6.15){$6$}

\put(0.5,4.8){$-2$}
\put(1.2,6.1){$-1$}
\put(-0.2,6.1){$-1$}
\put(-1.3,6.1){$1$}
\put(2.5,6.1){$1$}

\put(-1.8,0.4){$1$}
\put(-4.9,1.2){$1$}
\put(-0.4,2.3){$1$}
\put(-3.6,0.67){$-\frac{1}{2}$}
\put(-1.2,1.3){$-\frac{1}{2}$}
\put(-2.6,1.65){$-\frac{1}{2}$}
\end{picture}
\caption{The values of the ultimate eigenfunctions on $V_1$. (We only denote the non-zero values.)}
\end{center}
\end{figure}

\end{document}